\documentclass[acmsmall,screen,autherversion,final]{acmart}

\AtBeginDocument{%
  \providecommand\BibTeX{{%
    \normalfont B\kern-0.5em{\scshape i\kern-0.25em b}\kern-0.8em\TeX}}}

\setcopyright{acmcopyright}
\acmJournal{TELO}
\acmYear{2022} \acmVolume{1} \acmNumber{1} \acmArticle{1} \acmMonth{1} \acmPrice{15.00}\acmDOI{10.1145/3510425}

\usepackage{graphicx}
\graphicspath{{./}} 

\usepackage{subcaption}
\usepackage{booktabs} 
\usepackage{longtable}

\usepackage{latexsym}
\usepackage{amsmath} 
\usepackage{mathtools}       
\usepackage{stmaryrd}        
\usepackage{bm}              

\providecommand{\argmax}{\operatornamewithlimits{argmax}} 
\providecommand{\argmin}{\operatornamewithlimits{argmin}} 
\DeclareMathOperator{\Cond}{Cond} 
\DeclareMathOperator{\diag}{diag} 
\providecommand{\R}{\mathbb{R}} 
\providecommand{\E}{\mathbb{E}} 
\providecommand{\T}{\mathrm{T}} 
\renewcommand{\geq}{\geqslant} 
\renewcommand{\leq}{\leqslant} 
\DeclarePairedDelimiterX{\inner}[2]{\langle}{\rangle}{#1, #2}
\DeclarePairedDelimiter{\norm}{\lVert}{\rVert}
\DeclarePairedDelimiter{\abs}{\lvert}{\rvert}

\usepackage{amsthm}

\usepackage{xcolor}
\usepackage{ifthen}

\usepackage{cleveref}

\newtheorem{theorem}{Theorem}[section]
\newtheorem{proposition}[theorem]{Proposition}

\newtheorem{lemma}[theorem]{Lemma}
\theoremstyle{definition}
\newtheorem{definition}[theorem]{Definition}
\newtheorem{assumption}[theorem]{Assumption}

\usepackage{algorithm}
\usepackage{algpseudocode}

\providecommand{\X}{\mathbb{X}}
\providecommand{\Y}{\mathbb{Y}}
\providecommand{\Z}{\mathbb{Z}}
\providecommand{\XX}{\R^m}
\providecommand{\YY}{\R^n}
\providecommand{\ZZ}{\R^{\ell}}

\begin{document}

\title[Saddle Point Optimization with Approximate Minimization Oracle]{Saddle Point Optimization with Approximate Minimization Oracle and its Application to Robust Berthing Control}


\author{Youhei Akimoto}
\email{akimoto@cs.tsukuba.ac.jp}
\orcid{0000-0003-2760-8123}
\affiliation{%
  \institution{Faculty of Engineering, Information and Systems, University of Tsukuba \& RIKEN Center for Advanced Intelligence Project}
  \streetaddress{1-1-1 Tennodai}
  \city{Tsukuba}
  \state{Ibaraki}
  \country{Japan}
  \postcode{305-8573}
}

\author{Yoshiki Miyauchi}
\email{yoshiki_miyauchi@naoe.eng.osaka-u.ac.jp}
\author{Atsuo Maki}
\email{maki@naoe.eng.osaka-u.ac.jp}
\orcid{0000-0002-2819-1297}
\affiliation{%
  \institution{Department of Naval Architecture and Ocean Engineering, Graduate School of Engineering, Osaka University}
  \streetaddress{2-1 Yamadaoka}
  \city{Suita}
  \state{Osaka}
  \postcode{565-0971}
  \country{Japan}}



\begin{abstract}
We propose an approach to saddle point optimization relying only on oracles that solve minimization problems approximately. 
We analyze its convergence property on a strongly convex--concave problem and show its linear convergence toward the global min--max saddle point.
Based on the convergence analysis, we develop a heuristic approach to adapt the learning rate. 
An implementation of the developed approach using the (1+1)-CMA-ES as the minimization oracle, namely Adversarial-CMA-ES, is shown to outperform several existing approaches on test problems. 
Numerical evaluation confirms the tightness of the theoretical convergence rate bound as well as the efficiency of the learning rate adaptation mechanism.
As an example of real-world problems, the suggested optimization method is applied to automatic berthing control problems under model uncertainties, showing its usefulness in obtaining solutions robust to uncertainty.
\end{abstract}

\maketitle
\sloppy

\section{Introduction}

\emph{Simulation-based optimization} has recently received increasing attention from researchers. Here, the objective function $h: \X \to \R$, where $\X \subseteq \R^m$, is not explicitly given, but its value for each $x \in \X$ can be computed through computational simulation.
Solvers for simulation-based optimization problems have been widely developed. 
While some are domain-specific, others are general-purpose solvers. 
For a case where simulation-based optimization is required, we first need to design a simulator that models reality, for example, a physical equation, and computes the objective function value for each solution. 
Then, we apply a numerical solver to solve $\argmin_{x \in \X} h(x)$. 
However, owing to modeling errors and uncertainties, the optimal solution to $\argmin_{x \in \X} h(x)$ computed through a simulator is not necessarily optimal in real environments in which the obtained solution is used. 
This issue threatens the reliability of solutions obtained through simulation-based optimization.

An approach to obtain a solution that is robust against modeling errors and uncertainty is to formulate the problem as a min--max optimization 
\begin{equation}
  \min_{x \in \X} \max_{y \in \Y} f(x, y) ,\label{eq:spo}
\end{equation}
where $y \in \Y$ represents the model parameters and the uncertain parameters. In the following, $y$ is referred to as the \emph{uncertainty parameter}.
Assume that the real environment is represented by $y_\text{real} \in \Y$. 
The original objective $h(x)$ is equivalent to $f(x, y_\text{est})$ with an estimated parameter $y_\text{est} \in \Y$. 
Then, the solution $x_{y_\text{est}} = \argmin_{x \in \X} f(x, y_\text{est})$ obtained via simulation does not guarantee good performance in the real environment. That is, $f(x_{y_\text{est}}, y_\text{real})$ may be arbitrarily greater than $f(x_{y_\text{est}}, y_\text{est})$.
In contrast, the solution $x_{\Y} = \argmin_{x \in \X} \max_{y \in \Y}f(x, y)$ to \eqref{eq:spo} guarantees that $f(x_{\Y}, y_\text{real}) \leq \max_{y \in \Y}f(x_{\Y}, y)$. That is, by minimizing the worst-case objective value, one can guarantee performance in the real environment provided that $y_\text{real} \in \Y$. 

\paragraph{Robust Berthing Control}

As an important real-world application of the min--max optimization~\eqref{eq:spo}, we consider an automatic ship berthing task \cite{Maki2020b,maki2020}, which can be formulated as an optimization of the feedback controller of a ship.
Currently, the domestic shipping industry in Japan is facing a shortage of experienced on-board officers. Moreover, the existing workforce of officers is aging \cite{whitepaper2020}. This has generated considerable interest in autonomous ship operation to improve maritime safety, shipboard working environments, and productivity, and the technology is being actively developed.
Automatic berthing/docking requires fine control so that the ship can reach the target position located near the berth but avoid colliding with it. 
Therefore, automatic berthing is central to the realization of automatic ship operations. 
Because it is difficult to train the controller in a real environment owing to cost and safety issues, a typical approach first models the state equation of a ship, for example, using system identification techniques \cite{Abkowitz1980,Araki2012, Miyauchi20201SI, Wakita20201NN} and then optimizes the feedback controller in a simulator.
However, such an approach always suffers from modeling errors and uncertainties. For instance, the coefficients of a state equation model are often estimated based on captive model tests in towing tanks and regressions; hence, they may include errors. Moreover, the weather conditions at the time of operation could differ from those considered in the model. Optimization of the feedback controller on a simulator with an estimated model may result in a catastrophic accident, such as collision with the berth. Thus, to design a berthing control solution robust against modeling errors and uncertainties, we formulate the problem as a min--max optimization \eqref{eq:spo}, where $x$ is the parameter of the feedback controller and $y$ is the parameter representing the coefficients of the state equation model and weather conditions.

\paragraph{Saddle Point Optimization}

Here, we consider min--max continuous optimization \eqref{eq:spo}, where $f: \X \times \Y \to \R$ is the objective function and $\X \times \Y \subseteq \XX \times \YY$ is the search domain. 
In addition to the abovementioned situation, min--max optimization can be applied in many fields of engineering, including robust design \cite{Conn2012bilevel,mmde2018}, robust control \cite{Pinto2017icml,Shioya2018iclr}, constrained optimization \cite{Cherukuri2017sicon}, and generative adversarial networks (GANs) \cite{goodfellow2014generative,Salimans2016nips}. 
In particular, we are interested in the min--max optimization of a \emph{derivative-free} and \emph{black-box} objective $f$, where the gradient or higher-order information is unavailable (derivative-free) and no special structures such as convexity or the Lipschitz constant are available in advance (black-box) \cite{frazier2018tutorial}.\footnote{In this paper, \emph{simulation-based optimization} is used to refer to the problems described in the first paragraph above. The terms \emph{derivative-free} and \emph{black-box} are used to refer to the characteristics of the objective function.}

We aim to locate a local min--max saddle point of $f$, that is, a point $(x^*, y^*)$ satisfying $f(x, y^*) \geq f(x^*, y^*) \geq f(x^*, y)$ in a neighborhood of $(x^*, y^*)$. 
Generally, it is difficult to locate the global minimum of the worst-case objective $F(x) := \max_{y \in \Y} f(x, y)$. In a non-convex optimization context, the goal is often to locate a local minimum of an objective rather than the global minimum as a realistic target. However, in the min--max optimization context, it is still difficult to locate a local minimum of the worst-case objective $F(x)$ because doing so requires the maximization itself and there may exist local maxima of $f(x, y)$ unless $f(x, y)$ is concave in $y$ for all $x$.
A local min--max saddle point is considered as a local optimal solution in the min--max optimization context because it is a local minimum in $x$ and a local maximum in $y$.
Therefore, as a practical target, we focus on locating the local min--max saddle point of \eqref{eq:spo}. 

\paragraph{Related Works}

First-order approaches are often employed for \eqref{eq:spo} if gradients are available. A simultaneous gradient descent-ascent (GDA) approach
\begin{equation}
  (x_{t+1}, y_{t+1}) = (x_t, y_t) + \eta (- \nabla_x f(x_t, y_t), \nabla_y f(x_t, y_t))
  \label{eq:grad}
\end{equation}
has often been analyzed for its local and global convergence properties on twice continuously differentiable functions owing to its simplicity and popularity. 
A condition on the learning rate $\eta > 0$ for the dynamics \eqref{eq:grad} to be asymptotically stable at a local min--max saddle point has been studied \cite{Nagarajan2017nips,Mescheder2017nips}. Subsequently, \citet{pmlr-v89-adolphs19a} showed the existence of asymptotically stable points of \eqref{eq:grad} that are not local min--max saddle points. \citet{liang2019interaction} have derived a sufficient condition on $\eta$ for \eqref{eq:grad} to converge toward the global min--max saddle point on a locally strongly convex--concave function. Frank-Wolfe type approaches have also been analyzed for constrained situations \cite{gidel:hal-01403348,NIPS2019_9631}.
Although a convergence guarantee was not provided, \cite{Bertsimas.2010a,Bertsimas.2010b} have proposed a first-order approach targeting on $f$ that is non-concave in $y$.

Zero-order approaches for \eqref{eq:spo} include coevolutionary approaches \cite{Al-Dujaili2019lego,mmde2018,Jensen2004book,Branke2008ppsn,Zhou.2010}, surrogate-model--based approaches \cite{Picheny2019bo,Bogunovic2018nips,Conn2012bilevel}, and gradient approximation approaches \cite{Liu2020icml}.
Compared to first-order approaches, zero-order approaches have not been thoroughly analyzed in terms of their convergence guarantees and convergence rates. In particular, coevolutionary approaches are often designed heuristically and without convergence guarantees. Indeed, they fail to converge toward a min--max saddle point even on strongly convex--concave problems, as has been reported in \cite{akimoto2021} and as noted below in the experimental results.
Recently, \citet{Bogunovic2018nips} showed regret bounds for a Bayesian optimization approach and \citet{Liu2020icml} showed an error bound for a gradient approximation approach, where the error is measured by the square norm of the gradient. 
Both analyses show sublinear rates under possibly stochastic (i.e., noisy) versions of \eqref{eq:spo}.
However, compared to the first-order approach, which exhibits linear convergence, they show slower convergence. 

\paragraph{Contributions}

We propose an approach to saddle point optimization \eqref{eq:spo} that relies solely on numerical solvers that approximately solve $\argmin_{x' \in \X} f(x', y)$ for each $y \in \Y$ and $\argmin_{y' \in \Y} - f(x, y')$ for each $x \in \X$.
Given an initial solution $(x_0, y_0) \in \X\times\Y$, our approach iteratively locates the approximate solutions  $\tilde{x}_t \approx \argmin_{x' \in \X} f(x', y_t)$ and $\tilde{y}_t \approx \argmin_{y' \in \Y} -f(x_t, y')$ and updates the solution as
\begin{align}
    (x_{t+1}, y_{t+1}) = (x_t, y_t) + \eta \cdot (\tilde{x}_{t} - x_t, \tilde{y}_{t} - y_t) ,\label{eq:obspo}
\end{align}
where $\eta > 0$ is the learning rate.
This approach takes inspiration from the GDA method \eqref{eq:grad}, where we replace $-\nabla_x f(x_t, y_t)$ and $\nabla_y f(x_t, y_t)$ with $\tilde{x}_t - x_t$ and $\tilde{y}_t - y_t$.
However, unlike the GDA approach, the solvers need not be gradient-based. This is advantageous in the following situations: (1) there exists a well-developed numerical solver suitable for $\argmin_{x' \in \X} f(x', y)$ and/or $\argmin_{y' \in \Y} -f(x, y')$; (2) derivative-free approaches such as the covariance matrix adaptation evolution strategy (CMA-ES) \cite{hansen2001completely,hansen2003reducing,hansen2014principled,akimoto2020diagonal} are sought because gradient information is not available or gradient-based approaches are known to be sensitive to their initial search points. 

We analyze the proposed approach on strongly convex--concave problems, and prove its linear convergence in terms of the number of numerical solver calls.
In particular, we provide an upper bound on $\eta$ to guarantee linear convergence toward the global min--max saddle point and the convergence rate bound.
This corresponds to the known result for the GDA approach~\eqref{eq:grad}.
Compared to existing derivative-free approaches for saddle point optimization, this result is unique in that our convergence is linear, while the existing results show sublinear convergence \cite{Bogunovic2018nips,Liu2020icml}.
Although our motivational application is not necessarily a strongly convex--concave problem, the quantitative analysis helps to understand limitations of the approach---need for $\eta$ adaptation---and provide inspiration on how to improve the approach.

Moreover, we developed a heuristic adaptation mechanism for the learning rate in a black-box optimization setting.
In the black-box setting, we do not know in advance the characteristic constants of a problem that determines the upper bound for the learning rate to guarantee convergence.
Therefore, a learning rate adaptation mechanism is highly desired to avoid trial and error in tuning the learning rate.
We implemented two variants of the proposed approach, one using (1+1)-CMA-ES \cite{arnold2010}, a zero-order approach, as the minimization solver, and another using SLSQP \cite{kraft1988software}, a first-order approach. 
Empirical studies on test problems show that the learning rate adaptation achieved performance competitive to the proposed approach with the optimal static learning rate, while obviating the need for time-consuming parameter tuning.
We also demonstrate the limitations of existing coevolutionary approaches as well as the proposed approach. 

We apply our approach to robust berthing control optimization, as an example of a real-world application with a non-convex-concave objective. 
We consider the wind conditions and the coefficients of the state equation for the wind force as the uncertainty parameter $y$. 
Some related works address the wind force as an external disturbance when planning the trajectories \cite{miyauchi2021}; however, they treat the wind condition as an observable disturbance, and the control signal is selected according to the observed wind condition. In contrast, we optimize the on-line feedback controller under wind disturbance without considering the wind condition as an input to the controller. Moreover, among other studies on automatic berthing control, to the best of our knowledge, the present work is the first to address model uncertainty.
Compared to a naive baseline approach, the proposed approach located solutions with better worst-case performance.

This paper is an extension of a previous work \cite{akimoto2021}.
We have improved on the previous work in the following respects. 
First, we improved the convergence analysis in \Cref{sec:analysis}. We have removed unnecessary assumptions on the problem by refining the proof. Second, we have incorporated the covariance matrix adaptation into our proposed approach in \Cref{sec:acmaes}. Third, we have implemented a restart strategy and other ingenuity for practical use, summarized in \Cref{sec:restart}. Fourth, we have extended the comparison with existing approaches in \Cref{sec:comp}. Finally, we have evaluated the usefulness of the proposed approach in a real-world application in \Cref{sec:berth}.

Our implementation of the proposed approach in the Python programming language, Adversarial-CMA-ES, is publicly available at GitHub Gist.\footnote{%
\url{https://gist.github.com/youheiakimoto/ab51e88c73baf68effd95b750100aad0}
}%

\paragraph{Notation}

For a twice continuously differentiable function $f: \XX\times\YY \to \R$, that is, $f \in \mathcal{C}^2(\XX\times\YY, \R)$, let $H_{x,x}(x,y)$, $H_{x,y}(x,y)$, $H_{y,x}(x,y)$, and $H_{y,y}(x,y)$ be the blocks of the Hessian matrix $\nabla^2 f(x, y) = \begin{bmatrix} H_{x,x}(x,y) & H_{x,y}(x,y) \\ H_{y,x}(x, y) & H_{y,y}(x, y) \end{bmatrix}$ of $f$, whose $(i, j)$-th components are $\frac{\partial^2 f}{\partial x_i \partial x_j}$, $\frac{\partial^2 f}{\partial x_i \partial y_j}$, $\frac{\partial^2 f}{\partial y_i \partial x_j}$, and $\frac{\partial^2 f}{\partial y_i \partial y_j}$, respectively, evaluated at a given point $(x, y)$.

For symmetric matrices $A$ and $B$, by $A \succcurlyeq B$ and $A \succ B$, we mean that $A - B$ is non-negative and positive definite, respectively. 
For simplicity, we write $A \succcurlyeq a$ and $A \succ a$ for $a \in \R$ to mean $A \succcurlyeq a \cdot I$ and $A \succ a \cdot I$, respectively.
For a positive definite symmetric matrix $A$, let $\sqrt{A}$ denote the matrix square root, that is, $\sqrt{A}$ is a positive definite symmetric matrix such that $A = \sqrt{A}\cdot\sqrt{A}$.
Let $\norm{z}_{A} = [ z^\T A z ]^{1/2}$ for a positive definite symmetric $A$. 

Let $J_g(z)$ denote the Jacobian of a differentiable $g = (g_1, \dots, g_k): \ZZ \to \R^k$, where the ($i$, $j$)-th element is $\partial g_i / \partial z_j$ evaluated at $z = (z_1,\dots,z_\ell) \in \ZZ$. If $k = 1$, we write $J_g(z) = \nabla g(z)^\T$.

\section{Saddle Point Optimization}

Our objective is to locate the global or local min--max saddle point of the min--max optimization problem \eqref{eq:spo}.
In the following we first define the min--max saddle point.
We introduce the notion of the suboptimality error to measure the progress toward the global min--max saddle point.
Finally, we introduce a strongly convex--concave function as an important class of the objective function, on which we performed convergence analysis, which is described in the next section.

\subsection{Min--Max Saddle Point}

The min--max saddle point of a function $f:\X \times \Y \to \R$ is defined as follows. 
Here, $\mathcal{E}_z \subseteq \R^\ell$ is called a \emph{neighborhood} of $z^* \in \R^{\ell}$ if there exists an open ball $B(z^*, r) = \{z \in \R^\ell : \norm{z - z^*} < r\}$ such that $B(z^*, r) \subseteq \mathcal{E}_z$. A point $(x, y)$ is called a \emph{critical point} if $\nabla f(x, y) = (\nabla_x f(x, y), \nabla_y f(x, y)) = 0$.

\begin{definition}[Min--Max Saddle Point]\label{def:saddle}
A point $(x^*, y^*) \in \X \times \Y$ is a \emph{local min--max saddle point} of a function $f: \X \times \Y \to \R$ if there exists a neighborhood $\mathcal{E}_x \times \mathcal{E}_y \subseteq \X \times \Y$ including $(x^*, y^*)$ such that for any 
$(x, y) \in \mathcal{E}_x\times \mathcal{E}_y$, the condition 
$f(x, y^*) \geq f(x^*, y^*) \geq f(x^*, y)$ holds.
If $\mathcal{E}_x = \X$ and $\mathcal{E}_y = \Y$, the point $(x^*, y^*)$ is called the \emph{global min--max saddle point}. If the equality holds only if $(x, y) = (x^*, y^*)$, it is called a \emph{strict min--max saddle point}.
\end{definition}

For twice continuously differentiable function $f \in \mathcal{C}^2(\X\times\Y, \R)$, a point $(x^*, y^*)$ is a strict local min--max saddle point if it is a critical point and $H_{x,x}(x^*, y^*) \succ  0$ and $H_{y,y}(x^*, y^*) \prec 0$ hold. 
In general, the opposite does not hold. For example, a local min--max saddle point can be a boundary point of $\X\times\Y$ that is not a critical point. 

We comment on the relation between the min--max saddle point and the solutions to the worst-case objective function $F(x) := \max_{y \in \Y} f(x, y)$.
If there exists a global min--max saddle point $(x^*, y^*)$ of $f$, then $x^*$ is one of the global minimal point of the worst-case objective function $F$ and we have $F(x^*) = f(x^*, y^*)$.
However, a global minimal point $\bar{x} \in \X$ of $F(x)$ and its corresponding worst uncertainty parameter $\bar{y} \in \argmax_{y \in \Y} f(\bar{x}, y)$ do not necessarily form a global min--max saddle point in general. 
An example case is $f(x, y) = (x + \sin(\pi y))^2$, where the worst-case objective function is $F(x) = (\abs{x} + 1)^2$ and its global minimal point is $\bar{x} = 0$. The corresponding worst uncertainty parameters are $\frac12 + i$ for all integer $i$, but $(\bar{x}, \bar{y})$ is not even a local min--max saddle point for any $\bar{y} = \frac12 + i$. Moreover, if $(x^*, y^*)$ is a local min--max saddle point of $f$, the point $x^*$ is not necessarily a local minimal point of $F$.

\subsection{Suboptimality Error}

The suboptimality error \cite{gidel:hal-01403348} is a quantity that measures the progress toward the global min--max saddle point, defined as follows.
\begin{definition}[Suboptimality Error]\label{def:subopt}
  For function $f:\X\times\Y \to \R$, the suboptimality error $G_x:\X\times\Y\to[0, \infty)$ in $x$ and the suboptimality error $G_y:\X\times\Y\to[0, \infty)$ in $y$ are defined as 
\begin{align}
  G_x(x, y) &:= f(x, y) - \min_{x' \in \X} f(x', y) ,\\
  G_y(x, y) &:= \max_{y' \in \Y} f(x, y') - f(x, y) ,
    \label{eq:subopt-xy}
\end{align}
and the \emph{suboptimality error} is 
\begin{align}
  G(x, y) &:= G_x(x, y) + G_y(x, y) = \max_{y' \in \Y} f(x, y') - \min_{x' \in \X} f(x', y)
    . \label{eq:subopt-global}
\end{align}
\end{definition}

The suboptimality error $G(x, y)$ is zero if and only if $(x, y)$ is the global min--max saddle point of~$f$. 
Moreover, the local min--max saddle points of $f$ are characterized by suboptimality errors.
This is summarized in the following proposition, whose proof is given in \Cref{apdx:proof}.
\begin{proposition}\label{prop:saddle-error}
The point
  $(x^*, y^*)$ is the global min--max saddle point of $f$ if and only if it is the 
  global minimal point of $G$, that is, $G(x, y) \geq 0$ for any $(x, y) \in \X\times\Y$.
  The point $(x^*, y^*)$ is a local min--max saddle point of $f$ if and only if $x^*$ and $y^*$ are  
  local minimal points of $G_x(\cdot, y^*)$ and $G_y(x^*, \cdot)$, respectively, that is, there exists a neighborhood $\mathcal{E}_x\times\mathcal{E}_y$ of $(x^*, y^*)$ such that $G_x(x, y^*) \geq G_x(x^*, y^*)$ and $G_y(x^*, y) \geq G_y(x^*, y^*)$ for any $(x, y) \in \mathcal{E}_x\times\mathcal{E}_y$.
\end{proposition}

\subsection{Strongly Convex--Concave Function}

A strongly convex--concave function is often considered for the theoretical investigation of first-order min--max optimization approaches.

\begin{definition}
  A twice continuously differentiable function $f \in \mathcal{C}^2(\XX\times\YY, \R)$ is \emph{locally $\mu$-strongly convex--concave} around a critical point $(x^*, y^*)$ for some $\mu > 0$ if there exist open sets $\mathcal{E}_x \subseteq \XX$ including $x^*$ and $\mathcal{E}_y \subseteq \YY$ including $y^*$ such that $H_{x,x}(x, y) \succcurlyeq \mu$ and $- H_{y,y}(x, y) \succcurlyeq \mu$ for all $(x, y) \in \mathcal{E}_x \times \mathcal{E}_y$. The function $f$ is \emph{globally $\mu$-strongly convex--concave} if $\mathcal{E}_x = \XX$ and $\mathcal{E}_y = \YY$.
  We say that $f$ is locally or globally strongly convex--concave if $f$ is locally or globally $\mu$-strongly convex--concave for some $\mu > 0$. 
\end{definition}

If the objective function $f$ is a globally strongly convex--concave, the global minimal point of the worst-case objective function $F(x)$ is the global min--max saddle point, and it is the only local min--max saddle point. 

\providecommand{\opty}{\hat{y}}
\providecommand{\optx}{\hat{x}}

The implicit function theorem, for example, Theorem~5 of \cite{deoliveira2013}, provides important characteristics of strongly convex--concave functions. 
\begin{proposition}[Implicit Function Theorem]\label{prop:implicit}
  Let $(x^*, y^*)$ be a min--max saddle point of $f \in C^2(\XX\times\YY, \R)$ and $f$ be (at least) locally strongly convex--concave around $(x^*, y^*)$ in $\mathcal{E}_{x} \times \mathcal{E}_{y} \subseteq \XX \times \YY$.
  
  There exist open sets $\mathcal{E}_{x,x} \subseteq \mathcal{E}_x$ including $x^*$ and $\mathcal{E}_{x,y} \subseteq \mathcal{E}_y$ including $y^*$,
  such that there exists a unique $\opty: \mathcal{E}_{x,x} \to \mathcal{E}_{x,y}$ such that $\nabla_y f(x, \opty(x)) = 0$. Moreover, $y^* = \opty(x^*)$ and $J_{\opty}(x) = -(H_{y,y}(x, \opty(x)))^{-1} H_{y,x}(x, \opty(x))$ for all $x \in \mathcal{E}_{x,x}$.

  Analogously, 
  there exist open sets $\mathcal{E}_{y,y} \subseteq \mathcal{E}_y$ including $y^*$ and $\mathcal{E}_{y,x} \subseteq \mathcal{E}_x$ including $x^*$,
  such that there exists a unique $\optx: \mathcal{E}_{y,y} \to \mathcal{E}_{y,x}$ such that $\nabla_x f(\optx(y), y) = 0$. Moreover, $x^* = \optx(y^*)$ and $J_{\optx}(y) = -(H_{x,x}(\optx(y), y))^{-1} H_{x,y}(\optx(y), y)$ for all $y \in \mathcal{E}_{y,y}$.

  If $f$ is globally strongly convex--concave, one can take $\mathcal{E}_{x,x} = \mathcal{E}_{y,x} = \XX$ and $\mathcal{E}_{y,y} = \mathcal{E}_{x,y} = \YY$ in the above statements. 
\end{proposition}

\Cref{prop:implicit} states that for a globally strongly convex-concave $f \in C^2(\X\times\Y, \R)$, for each $x \in \XX$ there exists a unique global maximal point $\opty(x) \in \YY$ such that $\opty(x) = \argmax_{y \in \YY} f(x, y)$, and for each $y \in \YY$ there exists a unique global minimal point $\optx(y) \in \XX$ such that $\optx(y) = \argmin_{x \in \XX} f(x, y)$.

The following lemma shows the positivity of the Hessian of the suboptimality error $G$, which implies that the suboptimality error $G$ is a globally strongly convex function. The proof is provided in \Cref{apdx:proof}.

\begin{lemma}\label{lem:hess-gap}
  Suppose that $f \in \mathcal{C}^2(\XX\times\YY, \R)$ is globally $\mu$-strongly convex--concave for some $\mu > 0$.
  The Hessian matrix of the suboptimality error $G$ is $\nabla^2 G(x, y) = \diag(G_{x,x}(x, \opty(x)), G_{y,y}(\optx(y), y))$, where 
  \begin{align*}
    G_{x,x}(x, y) &= H_{x,x}(x, y) + H_{x,y}(x, y)(- H_{y,y}(x, y) )^{-1}  H_{y,x}(x, y) \\
    G_{y,y}(x, y) &= - H_{y,y}(x, y ) + H_{y,x}(x, y)(H_{x,x}(x, y))^{-1}  H_{x,y}(x, y) 
  \end{align*}
  and they are symmetric, and $G_{x,x}(x, y) \succcurlyeq \mu$ and $G_{y,y}(x, y) \succcurlyeq \mu$.
\end{lemma}

\section{Oracle-based Saddle Point Optimization}

We now analyze saddle point optimization based on the approximate minimization oracle outlined in \eqref{eq:obspo}. 
In the following, we formally state the condition for the approximate minimization oracle. Then, we show the global convergence of \eqref{eq:obspo} on strongly convex--concave functions. 

\subsection{Approximate Minimization Oracle}


First, we formally define the requirement for the minimization problem solvers.
\begin{definition}[Approximate Minimization Oracle]
  Given an objective function $h:\Z \to \R$ to be minimized and a reference solution $\bar{z} \in \Z$, an \emph{approximate minimization oracle} $\mathcal{M}$ outputs a solution $\tilde{z} = \mathcal{M}(h, \bar{z})$ satisfying $h(\tilde{z}) < h(\bar{z})$ unless $\bar{z}$ is a local minimal point of $h$. 
\end{definition}

We now reformulate the saddle point optimization with approximate minimization oracles.
Suppose that we have an approximate minimization oracle $\mathcal{M}_{x}$ solving $\argmin_{x'\in\X} f(x', y)$ for any $y \in \Y$ and an approximate minimization oracle $\mathcal{M}_{y}$ solving $\argmin_{y'\in\Y} - f(x, y')$ for any $x \in \X$.
At each iteration, the algorithm asks the approximate minimization oracles to output the approximate solutions to $\argmin_{x' \in \X} f(x', y_t)$ and $\argmin_{y' \in \X} f(x_t, y')$ with the current solution $(x_t, y_t)$ as their reference point.
Let $\tilde{x}_t = \mathcal{M}_{x}(f(\cdot, y_t), x_t)$ and $\tilde{y}_t = \mathcal{M}_{y}(-f(x_t, \cdot), y_t)$.
The update follows
\begin{equation}
  \begin{split}
  x_{t+1} &= x_t + \eta \cdot (\tilde{x}_t - x_t) ,\\
  y_{t+1} &= y_t + \eta \cdot (\tilde{y}_t - y_t) .
\end{split}\label{eq:algo}
\end{equation}

A point $(x, y) \in \X \times \Y$ is a stationary point of the dynamics of \eqref{eq:algo} only if it is a local min--max saddle point of $f$. 
Moreover, if $(x, y)$ is a strict local min--max saddle point of $f$, it is a stationary point of the dynamics of \eqref{eq:algo}.
Therefore, if it converges, the final solution is guaranteed to be a local min--max saddle point of $f$. To guarantee its convergence, we further assume the following requirement.
\begin{assumption}\label{def:amo}
  Given an objective function $h:\Z \to \R$ to be minimized and a reference solution $\bar{z} \in \Z$, an approximate minimization oracle $\mathcal{M}$ with an approximation precision parameter $\epsilon \in [0, 1)$ outputs a solution $\tilde{z} = \mathcal{M}(h, \bar{z})$ satisfying
  \begin{equation}
    h(\tilde{z}) - \min_{z \in \Z} h(z) \leq \epsilon \cdot (h(\bar{z}) - \min_{z \in \Z} h(z)) .
    \label{eq:amo-cond}
  \end{equation}
\end{assumption}

We are particularly interested in algorithms that decrease the objective function value at a geometric rate on (at least locally) strongly convex objective $h$ as instances of the approximate minimization oracle $\mathcal{M}$. That is, the runtime --- number of $h$ calls or $\nabla h$ calls --- to decrease the objective function difference $h(z) - h(z^*)$ from a local minimum by the factor $\epsilon$ is $O(\log(1/\epsilon))$.
For example, the gradient descent method is well known to exhibit a geometric decrease in the objective function value on strongly convex functions with Lipschitz continuous gradients \cite{polyakcondition2016,boyd2004convex}.
The (1+1)-ES also exhibits a geometric decrease on strongly convex functions with Lipschitz continuous gradients \cite{morinaga2019generalized}.
We can satisfy the oracle requirement \eqref{eq:amo-cond} by performing $O(\log(1/\epsilon))$ iterations of such algorithms.
The condition can also be satisfied by algorithms that exhibit slower convergence, that is, sublinear convergence.
However, for such algorithms, the runtime increases as a candidate solution becomes closer to a local optimum. Therefore, the stopping condition for the internal algorithm to satisfy \eqref{eq:amo-cond} needs to be carefully designed.

\subsection{Analysis on Strongly Convex--Concave Functions}\label{sec:analysis}

\providecommand{\bx}{\bar{x}}
\providecommand{\by}{\bar{y}}
\providecommand{\bbx}{\bar{\bar{x}}}
\providecommand{\bby}{\bar{\bar{y}}}

Next, we investigate the convergence property of the oracle-based saddle point optimization \eqref{eq:algo} on strongly convex--concave functions.
In particular, we are interested in knowing how small the learning rate $\eta$ needs to be to guarantee convergence and how fast it converges.
The following theorem provides an upper bound of the suboptimality error at iteration $t+1$ given the solution at iteration $t$. The proof is provided in \Cref{apdx:proof}. 

\begin{theorem}\label{thm:conv}
  Suppose that $f \in \mathcal{C}^2(\XX\times\YY, \R)$ is globally strongly convex--concave, and there exist $\beta_G \geq \alpha_G > 0$ and $\beta_H \geq \alpha_H > 0$ such that
  \begin{enumerate}
    \item $\beta_H \succcurlyeq \sqrt{H_{x,x}^*}^{-1} H_{x,x}(x, y) \sqrt{H_{x,x}^*}^{-1} \succcurlyeq \alpha_H$;
    \item $\beta_H \succcurlyeq \sqrt{-H_{y,y}^*}^{-1} (-H_{y,y}(x, y)) \sqrt{-H_{y,y}^*}^{-1} \succcurlyeq \alpha_H$;  
    \item $\beta_G \succcurlyeq \sqrt{H_{x,x}^*}^{-1} G_{x,x}(x, y) \sqrt{H_{x,x}^*}^{-1} \succcurlyeq \alpha_G$;
    \item $\beta_G \succcurlyeq \sqrt{-H_{y,y}^*}^{-1} G_{y,y}(x, y) \sqrt{-H_{y,y}^*}^{-1} \succcurlyeq \alpha_G$,
  \end{enumerate}
where $H_{x,x}^* = H_{x,x}(x^*, y^*)$, $H_{y,y}^* = H_{y,y}(x^*, y^*)$ and $(x^*, y^*)$ is the global min--max saddle point of~$f$. 
  Consider approach \eqref{eq:algo} with approximate minimization oracles $\mathcal{M}_x$ and $\mathcal{M}_y$ satisfying \Cref{def:amo} with approximate precision $\epsilon < \frac{\alpha_H^5}{\beta_H^4 \beta_G}$.
  Let
  \begin{align}
    \eta^* &= \frac{\alpha_H}{\beta_H}\frac{\alpha_H}{\beta_G} \frac{1 - (\beta_H/\alpha_H)^2 \sqrt{(\beta_G / \alpha_H) \cdot \epsilon}  }{(1 + \sqrt{\epsilon})^2} ,
     \label{eq:eta}\\
    \gamma &= - 2 \eta \frac{\alpha_{H}}{\beta_{H}} \left( 1 - \frac{\beta_H^2}{\alpha_H^2} \sqrt{\frac{\beta_G}{\alpha_{H}} \cdot \epsilon } \right)
    + \eta^2 \cdot (1 + \sqrt{\epsilon})^2 \cdot \frac{\beta_{G}}{\alpha_H} .\label{eq:cr}
  \end{align}
  Then, for any $\eta < 2 \cdot \eta^*$, we have $\gamma < 0$ and $\log\left( G(x_{t+1}, y_{t+1}) \right) - \log\left( G(x_{t}, y_{t}) \right) < \gamma$.
  In other words, the runtime $T_{\zeta}$ to reach $\{(x, y) \in \XX \times \YY : G(x, y) \leq \zeta \cdot G(x_{0}, y_{0}) \}$ for $\zeta \in (0, 1)$ is
  $T_{\zeta} \leq \left\lceil \frac{1}{\abs{\gamma}} \log\left(\frac{1}{\zeta}\right)  \right\rceil$ for any initial point $(x_0, y_0) \in \XX\times\YY$ .
  Moreover, $G(x_{t+1}, y_{t+1}) > G(x_t, y_t)$ for all $(x_t, y_t)$ if $\eta > 2 \cdot \bar{\eta}$, where
  \begin{equation}
    \bar{\eta} = \frac{\beta_H}{\alpha_G}\frac{\beta_H}{\alpha_H}  \frac{1 + \sqrt{(\beta_G / \alpha_H) \cdot \epsilon}}{  ( 1 - \sqrt{(\beta_H/ \alpha_H) \cdot \epsilon})^2 }  . \label{eq:eta-negative}
  \end{equation}
\end{theorem}

\paragraph{Linear Convergence}
The proposed approach \eqref{eq:algo} satisfying \Cref{def:amo} converges linearly toward the global min--max saddle point on a strongly convex--concave objective function if $\eta < 2 \eta^*$. 
If $\mathcal{M}_{x}$ and $\mathcal{M}_{y}$ are implemented with algorithms that exhibit linear convergence, we can conclude that the runtime in terms of $f$-calls and/or $\nabla f$-calls is 
\begin{equation}\label{eq:runtime-f}
\mathcal{O}\left( \frac{1}{\abs{\gamma}} \log\left(\frac{1}{\zeta}\right) \log\left(\frac{1}{\epsilon}\right) \right).
\end{equation}

\paragraph{Necessary Condition} To exhibit convergence, shrinking the learning rate $\eta$ is not only sufficient but also necessary.
To determine the closeness of the upper bound $2 \cdot \eta^*$ in the sufficient condition and the lower bound $2 \cdot \bar{\eta}$ in the necessary condition, consider a convex--concave quadratic function $f(x, y) = \frac{a}{2} x^2 + b x y - \frac{c}{2} y^2$ for instance, where $a > 0$, $b \in \R$ and $c > 0$. Then, we have $\alpha_H = \beta_H = 1$ and $\alpha_G = \beta_G = \frac{ac}{ac + b^2}$. Ignoring the effect of $\epsilon$, we have $\eta^* = \bar{\eta} = \frac{ac}{ac + b^2}$. This implies that the sufficient condition for linear convergence, $\eta < 2\cdot \eta^*$, is indeed the necessary condition for the convergence itself in this example situation.
This reveals a limitation of existing approaches \cite{Al-Dujaili2019lego,Pinto2017icml}, which corresponds to \eqref{eq:algo} with $\eta = 1$. 

\paragraph{Runtime Bound}
The runtime bound $T_\zeta$ is proportional to $\frac{1}{\abs{\gamma}}$ in \eqref{eq:cr}.
The log-progress bound~$\abs{\gamma}$ is roughly proportional to $2 \cdot \eta$ if $\eta \ll 1$.
That is, the runtime is proportional to $\frac{1}{2\cdot\eta}$.
The minimal runtime bound is obtained when $\eta = \eta^*$, where
\begin{align}
  \gamma = \gamma^* := - \frac{\alpha_H}{\beta_G} \left(\frac{\alpha_H}{\beta_H}\right)^2 \left(\frac{1 - (\beta_H/\alpha_H)^2 \sqrt{(\beta_G / \alpha_H) \cdot \epsilon}}{1 + \sqrt{\epsilon}} \right)^2.\label{eq:crstar}
\end{align}
Provided that $\epsilon \ll 1$, we have $\eta^* \approx \frac{\alpha_H}{\beta_G} \frac{\alpha_H}{\beta_H}$ and $\gamma^* \approx - \frac{\alpha_H}{\beta_G} \left(\frac{\alpha_H}{\beta_H}\right)^2$.
The main factor that limits $\eta^*$ and $\gamma^*$ is $\frac{\alpha_H}{\beta_G}$. As noted above in the above example of a convex--concave quadratic function, the ratio $\frac{\alpha_H}{\beta_G}$ is smaller as the influence of the interaction term between $x$ and $y$ on the objective function value is greater than that to the other terms, that is, as $b^2 / a c$ is greater.
The other factor, $\frac{\alpha_H}{\beta_H}$, is smaller as the condition number $\Cond(H_{x,x}(x, y) (H_{x,x}^*)^{-1})$ or $\Cond(H_{y,y}(x, y) (H_{y,y}^*)^{-1})$ is higher.
This depends on the change in the Hessian matrix over the search space $\XX\times\YY$. If the objective function is convex--concave quadratic, that is, $f(x, y) = \frac12 x^\T H_{x,x} x + x^\T H_{x,y} y + \frac12 y^\T H_{y,y} y$, the Hessian matrix is constant over the search space, and we have $\alpha_H / \beta_H = 1$, whereas $\beta_G = 1 + \sigma_{\max}^2(\sqrt{H_{x,x}}^{-1} H_{x,y} \sqrt{-H_{y,y}}^{-1})$ and $\alpha_G = 1 + \sigma_{\min}^2(\sqrt{H_{x,x}}^{-1} H_{x,y} \sqrt{-H_{y,y}}^{-1})$, where $\sigma_{\min}$ and $\sigma_{\max}$ denote the smallest and greatest singular values.
Therefore, we have
\begin{equation}
   \frac{1}{\abs{\gamma^*}} = \left( 1 + \sigma_{\max}^2(\sqrt{H_{x,x}}^{-1} H_{x,y} \sqrt{-H_{y,y}}^{-1}) \right) \left(\frac{1 + \sqrt{\epsilon}}{1 - \sqrt{\epsilon}} \right)^2.\label{eq:cr-quad}
\end{equation}

\paragraph{Comparison with GDA}
Theorem~1 of \cite{liang2019interaction} shows that the runtime $T_\zeta$ of the GDA \eqref{eq:grad} is $\mathcal{O}\left(\frac{1}{\gamma_\text{GDA}} \log\left(\frac{1}{\zeta}\right)\right)$, where
\begin{align}
    \frac{1}{\gamma_\text{GDA}} = 
  \frac{\max_{(x,y) \in \XX\times\YY}\lambda_{\max}(K(x,y))}{ \min_{(x,y) \in \XX\times\YY}\lambda_{\min}(\diag(H_{x,x}(x,y)^2, (-H_{y,y}(x,y))^2)) } ,
  \label{eq:grad-rate}
\end{align}
$\lambda_{\min}$ and $\lambda_{\max}$ denote the smallest and greatest eigenvalues, and
\begin{align}
  K(x, y) = \begin{bmatrix}
    H_{x,x}^2 + H_{x,y}H_{y,x} & - H_{x,x}H_{x,y} + H_{x,y} (-H_{y,y}) \\
   -H_{y,x} H_{x,x} + (-H_{y,y}) H_{y,x} & (-H_{y,y})^2 + H_{y,x}H_{x,y}
  \end{bmatrix},
\end{align}
where we drop $(x,y)$ from $H_{x,x}(x, y)$, $H_{y,y}(x, y)$, $H_{x,y}(x, y)$, and $H_{y,x}(x, y)$ for compact expression.
To compare this with our result, consider the pre-conditioned convex--concave quadratic function
$f(x, y) = \frac12 x^\T I x + x^\T \sqrt{H_{x,x}}^{-1} H_{x,y} \sqrt{-H_{y,y}}^{-1} y + \frac12 y^\T(- I) y$. 
Then, it may be easily observed that $\lambda_{\min}(\diag(H_{x,x}(x,y)^2, (-H_{y,y}(x,y))^2))  = 1 = \alpha_H = \beta_H$ and $\lambda_{\max}(K(x,y)) = 1 + \sigma_{\max}^2(\sqrt{H_{x,x}}^{-1} H_{x,y} \sqrt{-H_{y,y}}^{-1}) = \beta_G$. 
Therefore, $\frac{1}{\gamma_\text{GDA}} = 1 + \sigma_{\max}^2(\sqrt{H_{x,x}}^{-1} H_{x,y} \sqrt{-H_{y,y}}^{-1})$. Because the GDA requires one $\nabla f$ call per iteration, $T_\zeta$ is the runtime w.r.t.\ $\nabla f$ calls as well. 
It indicates that the runtime bound w.r.t.\ $f$ and/or $\nabla f$ calls are the same for the GDA and the oracle-based saddle point optimization \eqref{eq:algo} with linearly convergent oracles, whose runtime is obtained by substituting \eqref{eq:cr-quad} into \eqref{eq:runtime-f}, by ignoring the effect of $\epsilon$. 
Note, however that the runtime of the GDA depends on the pre-conditioning, as it is a first-order approach. The number of oracle calls of the oracle-based saddle point optimization is independent of the pre-conditioning, but the number of $f$ and/or $\nabla f$ calls in each oracle call may depend on the pre-conditioning. 


\section{Saddle Point Optimization with Learning Rate Adaptation}

In this section, we propose practical implementations of the saddle point optimization approach \eqref{eq:algo} with a heuristic mechanism to adapt the learning rate $\eta$. We implement the proposed approach using two minimization routines.
The first is (1+1)-CMA-ES, which is a zero-order randomized hill-climbing approach.
The second is SLSQP, which is a first-order deterministic hill-climbing approach.

\subsection{Learning Rate Adaptation}\label{sec:lra}

The main limitation of oracle-based saddle point optimization when it is applied to a simulation-based optimization task is that we rarely know the right $\eta$ value in advance.
As we see in \Cref{thm:conv}, $\eta < 2 \cdot \eta^*$ must be selected to guarantee convergence on a convex--concave function. However, the optimal value, $\eta^*$, depends on the problem characteristics and is unknown in advance when considering a black-box setting. 
In practice, it is a tedious task to find a reasonable $\eta$. 

To address this issue, we propose adapting $\eta$ during the optimization process.
The overall framework is presented in \Cref{alg:lr}, where we assume $f_\mathcal{Y} = f$ for the moment, that is, $\mathcal{Y} = \emptyset$ to simplify the main idea. 

The main idea is to estimate the convergence speed in terms of the suboptimality error by running $N_\text{step}$ iterations of algorithm \eqref{eq:algo} with a candidate learning rate $\eta_c$ (lines~\ref{l:xcyc}--\ref{l:slope}). If the estimated convergence speed $\tilde{\gamma}_c$ associated with $\eta_c$ is better (greater absolute value with a negative sign) than the estimated convergence speed $\tilde{\gamma}$ associated with the base learning rate $\eta$, we replace $\eta$ with $\eta_c$ (lines~\ref{l:gammacheck}--\ref{l:revert}). 
The next candidate learning rate is chosen randomly from $\min(\eta \cdot c_\eta, 1)$ (greater learning rate), $\eta$ (current learning rate), and $\min(\eta / c_\eta, \eta_{\min})$ (smaller learning rate) with equal probability, where $\eta_{\min}$ is the minimal learning rate value and $c_\eta > 1$ is the hyperparameter that determines the granularity of the $\eta$ update. A smaller $c_\eta$ results in a smoother $\eta$ change, but it may require more time to adapt $\eta$. It is advised to set $c_\eta < 2$ because \Cref{thm:conv} indicates that the upper bound on $\eta$ for convergence is $2 \cdot \eta^*$, where $\eta^*$ is the optimal value.  

\begin{algorithm}[t]
  \caption{Saddle Point Optimization with Learning Rate Adaptation}\label{alg:lr}
  \begin{algorithmic}[1]
    \Require $x\in\X$, $y\in\Y$, $\theta^x$, $\theta^y$, $a_\eta > 0$, $b_\eta \geq 0$, $c_\eta > 1$
    \Require (optional) $P_x$, $P_y$, $\eta_{\min} \geq 0$, $G_\text{tol}\geq 0$, $d^y_{\min} \geq 0$
    \State $\eta \leftarrow 1$, $\tilde\gamma \leftarrow 0$, $\mathcal{Y} \leftarrow \emptyset$, $\mathcal{X} \leftarrow \emptyset$
    \For{$t = 1, \cdots, T$}
    \State $(x_t, \tilde{x}, \theta_t^x, y_t, \tilde{y}, \theta_t^y) \leftarrow (x, x, \theta^x, y, y, \theta^y)$\label{l:keep}
    \State $\eta_c \leftarrow \{ \min(\eta \cdot c_\eta, 1), \eta, \max(\eta / c_\eta, \eta_{\min}) \}$ w.p.\ $1/3$ for each\label{l:etac}
    \State $N_\text{step} \leftarrow \lfloor b_\eta +  a_\eta / \eta_c \rfloor$
    \For{$s = 1, \cdots, N_\text{step}$}\label{l:xcyc}
    \State Let $f_{\mathcal{Y}}(x, y) := \max_{y' \in \mathcal{Y}\cup\{y\}} f(x, y')$\label{l:fy}
    \State $(\tilde{x}, \theta^x) \leftarrow (x, \theta^x_t)$ \textbf{ if } $f_{\mathcal{Y}}(\tilde{x}, y) > f_{\mathcal{Y}}(x, y)$
    \State $(\tilde{y}, \theta^y) \leftarrow (y, \theta^y_t)$ \textbf{ if } $f(x, \tilde{y}) < f(x, y)$
    \State $\tilde{x}, \theta^x \leftarrow \mathcal{M}_{x}(f_{\mathcal{Y}}(\cdot, y), \tilde{x}; \theta^x)$\label{l:xhat}
    \State $\tilde{y}, \theta^y \leftarrow \mathcal{M}_{y}(-f(x, \cdot), \tilde{y}; \theta^y)$\label{l:yhat}
    \State $(\tilde{x}, \theta^x) \leftarrow (x', \theta^x_{t})$ \textbf{ if } $f_{\mathcal{Y}}(x', y)< f_{\mathcal{Y}}(\tilde{x}, y)$ for $x' \sim P_x$\label{l:px}
    \If{$f(x, y') > f(x, \tilde{y})$ for $y' \sim P_y$}\label{l:py}
    \State $\mathcal{Y} \leftarrow \mathcal{Y} \cup \{\tilde{y}\}$ \textbf{if} $f(x, \tilde{y}) \geq f_\mathcal{Y}(x, y)$ and $\norm{\bar{y} - \tilde{y}} > d^y_{\min}$ for all $\bar{y} \in \mathcal{Y}$
    \State $(\tilde{y}, \theta^y) \leftarrow (y', \theta^y_{t})$
    \EndIf
\State $F_{s} \leftarrow  f_{\mathcal{Y}}(x, \tilde{y}) - f_{\mathcal{Y}}(\tilde{x}, y)$\label{l:gap}
    \State $(x, y) \leftarrow (x, y) + \eta_c (\tilde{x} - x, \tilde{y} - y)$\label{l:xyupdate}
    \State \textbf{break} \textbf{if} $s \geq b_\eta$ and $F_{s} > \cdots > F_{s-b_\eta+1}$
    \EndFor
    \State $\tilde\gamma_c,\ \sigma_{\tilde\gamma_c} \leftarrow  \textsc{slope}(\log(F_1), \dots, \log(F_{s}))$\label{l:slope}
    \If{$\tilde\gamma \geq 0$ and $\tilde\gamma_c \geq 0$} \label{l:gammacheck}
    \State $\eta \leftarrow \eta / c_\eta^3$\label{l:etadec}
    \ElsIf{$\tilde\gamma_c \leq \tilde\gamma$ or $\eta = \eta_c$}
    \State $\eta \leftarrow \eta_c$, $\tilde\gamma \leftarrow \tilde\gamma_c$\label{l:etaup}
    \EndIf
    \State $(x, y, \theta^x, \theta^y) \leftarrow (x_t, y_t, \theta_t^x, \theta_t^y)$ \textbf{if} $\tilde\gamma_c - 2 \sigma_{\tilde\gamma_c} > 0$\label{l:revert}
    \If{$F_{s} \leq G_\text{tol}$}\label{l:restart}
    \State $\mathcal{X} \leftarrow \mathcal{X} \cup \{x\}$
    \State $\mathcal{Y} \leftarrow \mathcal{Y} \cup \{y\}$ \textbf{if} $\norm{\bar{y} - \tilde{y}} > d^y_{\min}$ for all $\bar{y} \in \mathcal{Y}$
    \State Re-initialize $x$, $y$, $\theta^x$, $\theta^y$ and reset $\eta \leftarrow 1$ and $\tilde\gamma \leftarrow 0$
    \EndIf\label{l:restart-end}
    \EndFor
    \State \Return $\argmin_{x' \in \mathcal{X}\cup\{x\}} f_\mathcal{Y}(x', y)$
  \end{algorithmic}
\end{algorithm}

We estimate the convergence speed by running the algorithm for $N_\text{step}$ iterations. The suboptimality error $G(x, y)$ is approximated by $F_s$ in line~\ref{l:gap}. Because of the oracle condition \eqref{eq:amo-cond}, if there exists a unique (hence, global) min--max saddle point, we have $(1 - \max(\epsilon_x, \epsilon_y)) \cdot G(x, y) \leq F_s \leq G(x, y)$. Then, we have 
\begin{equation*}
 \frac{ 1 }{ N_{\mathrm{step}} - 1} \abs*{ \log\left( \frac{G(x_{N_\mathrm{step}}, y_{N_\mathrm{step}})}{ G(x_1, y_1) } \right)
    - \log\left( \frac{F_{N_\mathrm{step}} }{ F_1 } \right) }
  \leq \frac{ \abs*{ \log(1 - \max(\epsilon_x, \epsilon_y)) } }{ N_{\mathrm{step}} - 1}
  .
\end{equation*}
Based on \Cref{thm:conv}, if the objective function is strongly convex--concave, the convergence speed will be proportional to $1/\eta$.
Then, to approximate the convergence speed in line~\ref{l:gap}, $N_\mathrm{step} \in \Omega(1/\eta)$ must be set to alleviate the estimation error, that is, the right-hand side of the above inequality.
Therefore, we set $N_\text{step} = \lfloor b_\eta +  a_\eta / \eta_c \rfloor$, where $a_\eta > 0$ and $b_\eta \geq 0$ are the hyperparameters. The greater they are, the more accurate the estimated convergence speed, but the slower the speed of adaptation of $\eta$.
If the objective function is not strongly convex--concave, the above argument may not hold, yet we optimistically expect it to reflect the convergence speed of the algorithm toward a local min--max saddle point.

After estimating the convergence speed $\tilde{\gamma}_c$ for the candidate learning rate $\eta_c$, we replace $\eta$ and $\tilde\gamma$ with $\eta_c$ and $\tilde\gamma_c$ if $\tilde{\gamma}_c$ is equal to or better than the convergence speed $\tilde{\gamma}$ for the current learning rate~$\eta$~(line~\ref{l:etaup}). We also update $\tilde{\gamma}$ when $\eta = \eta_c$. If both $\tilde\gamma$ and $\tilde\gamma_c$ are non-negative, the learning rate is too high, and we reduce $\eta$ by multiplying $1/c_\eta^3$. If $\tilde\gamma_c - 2 \sigma_{\tilde\gamma_c} > 0$, where $\sigma_{\tilde\gamma}$ is the estimated standard deviation of $\tilde\gamma$, we revert the solutions and other strategy parameters $\theta^x$ and $\theta^y$.

Based on our preliminary experiments and the above argument, we set $a_\eta = 1$, $b_\eta = 5$, and $c_\eta = 1.1$ as the default values.

\subsection{Ingenuity for practical use}\label{sec:restart}

Our approach is designed to locate a min--max saddle point. 
However, in practice, we often cannot guarantee the existence of min--max saddle points. 
In such a situation, $x$ and $y$ may not converge and oscillate. 
For example, consider $f(x, y) = x^\T y$ on $[-1, 1] \times [-1, 1]$. 
The worst $y$ is $-1$ if $x < 0$ and $1$ if $x > 0$, and the best $x$ is $1$ if $y < 0$ and $-1$ if $y > 0$. 
This causes a cyclic behavior: $x$ is positive, then $y$ becomes positive, then $x$ becomes negative, then $y$ becomes negative, and so on. 
To stabilize the algorithm in such situations, we maintain a set $\mathcal{Y}$ of $y \in \Y$ and replace $f$ with $f_{\mathcal{Y}}(x, y) := \max_{y' \in \mathcal{Y}\cup\{y\}} f(x, y')$ using the approach described in \Cref{sec:lra}. In the above example, provided that there are points $y_1 < 0$ and $y_2 > 0$ in $\mathcal{Y}$, the optimal $x$ of $f_{\mathcal{Y}}$ is zero regardless of $y$. This is the optimal solution for $\min_{-1 \leq x \leq 1} \max_{-1 \leq y \leq 1} f(x, y)$. However, if we replace $-f$ with $-f_{\mathcal{Y}}$ for the objective function of $\mathcal{M}_y$, the optimization is likely to fail because $f_{\mathcal{Y}}(x, y)$ is constant with respect to $y$ over $\{y \in \Y : f(x, y) \geq f(x, y') \text{ for some } y' \in \mathcal{Y}\}$. 
Therefore, we replace $f$ with $f_{\mathcal{Y}}$ only for the parts regarding $x$ optimization.
We initialize $\mathcal{Y}$ with the empty set; hence, $f_\mathcal{Y} = f$ at the beginning. 
The output $\tilde{y}$ of $\mathcal{M}_y$ is registered to $\mathcal{Y}$ if a random sample $y' \sim P_y$ provides a worse objective value than $\tilde{y}$, $f(x, \tilde{y}) \geq f_{\mathcal{Y}}(x, y)$, and none of the registered points $\bar{y} \in \mathcal{Y}$ is in the closed ball centered at $\tilde{y}$ with radius $d^{y}_{\min}$, which is a hyperparameter. 

The existence of multiple local min--max saddle points is another difficulty that is often encountered in practice. 
For such problems, we would like to locate a local min--max saddle point whose worst-case objective value is as small as possible. 
To tackle this difficulty, we implement a restart strategy in lines~\ref{l:restart}--\ref{l:restart-end} of \Cref{alg:lr}. 
First, we check whether the current solution is nearly a local min--max saddle point by checking $F_{s} \leq G_\text{tol}$, where $G_\text{tol}$ is a user-defined threshold parameter. 
Note that $F_s$ can be close to zero at a local min--max saddle point even if the true suboptimality error is nonzero because $F_s$ is computed using the outputs of $\mathcal{M}_x$ and $\mathcal{M}_y$. Therefore, a small $F_{s}$ value is a sign of a local min--max saddle point. 
If this restart criterion is satisfied, we register the current solution $x$ as a candidate for the final solution and append the current $y$ to $\mathcal{Y}$ unless it is sufficiently close to the already registered points in $\mathcal{Y}$. Then, we re-initialize the solutions $x$ and $y$ and the internal parameters $\theta^x$ and $\theta^y$, and restart the search with $\eta = 1$. 

The other details are described as follows.
First, we allow the sharing of the internal parameters $\theta^x$ and $\theta^y$ of $\mathcal{M}_{x}$ and $\mathcal{M}_{y}$ over oracle calls. Second, we feed the last outputs $\tilde{x}$ and $\tilde{y}$ to $\mathcal{M}_{x}$ and $\mathcal{M}_{y}$ as the reference points instead of the current solutions $x$ and $y$ if the former is better. This contributes to realizing smaller approximation errors $\epsilon$. Third, we optionally try random samples $x' \sim P_x$ and $y' \sim P_y$ and check if they are better than the outputs of the oracles if $P_x$ and $P_y$ are given. A typical choice for $P_x$ and $P_y$ is the uniform distribution on $\X$ and $\Y$ if they are bounded. Fourth, we optionally introduce the minimal learning rate $\eta_{\min}$. Because a small $\eta$ slows down the optimization speed, it is not practical to set an extremely small $\eta$, even though it is necessary for convergence. 

\subsection{Adversarial-CMA-ES}\label{sec:acmaes}

We implemented the proposed approach with (1+1)-CMA-ES as $\mathcal{M}_{x}$ and $\mathcal{M}_{y}$.
The (1+1)-CMA-ES is a derivative-free randomized hill-climbing approach with step-size adaptation and covariance matrix adaptation. 
It samples a candidate solution $z' \sim \mathcal{N}(z, (\sigma A)(\sigma A)^\T)$, where $\sigma$ is the step size and $A \cdot A^\T$ is the covariance matrix. 
The step size is adapted with the so-called 1/5-success rule \cite{Devroye:72,Schumer:Steiglitz:68,rechenberg:1973}, which maintains $\sigma$ such that the probability of generating a better solution is approximately 1/5.
We implemented a simplified 1/5-success rule proposed by \cite{Kern2004}.
The covariance matrix was adapted with the active covariance matrix update \cite{arnold2010}.
The results show empirically that the covariance matrix learned the inverse Hessian matrix on a convex quadratic function.
The algorithm is summarized in \Cref{alg:es}.
We call \Cref{alg:lr} with \Cref{alg:es} as $\mathcal{M}_x$ and $\mathcal{M}_x$ \emph{Adversarial-CMA-ES}.

\begin{algorithm}[t]
  \caption{(1+1)-CMA-ES as Minimization Oracle}\label{alg:es}
  \begin{algorithmic}[1]
    \Require $h: \R^\ell \to \R$, $z \in \R^\ell$, $\sigma > 0$, $A \in \R^{\ell \times \ell}$, $h_{z} = h(z)$, $\tau_\text{es},\ \tau_\text{es}' \in \mathbb{N}$
    \Require (optional) $\bar{\sigma}_{\min} \geq 0$
    \State $c = e^{\frac{2}{2 + \ell}}$, $c_p = \frac{1}{12}$, $c_c = \frac{2}{\ell + 2}$, $c_\text{cov+} = \frac{2}{\ell^2 + 6}$, $c_\text{cov-} = \frac{0.4}{\ell^{1.6} + 1}$, $p_\text{thre}=0.44$
    \State $p \leftarrow 0 \in \R^{\ell}$,
    $p_\text{succ} \leftarrow 0.5 \in [0, 1]$,
    $n_\text{succ} = 0$
    \State Initialize $H \in \R^5$ with $H_1 = h_z$ and $H_2 = H_3 = H_4 = H_5 = \infty$
    \While{$n_\text{succ} < \tau_\text{es} \cdot \ell + \tau_\text{es}'$}
    \State $z' \leftarrow z + \sigma A \mathcal{N}(0, I)$\label{alg:aes:esstart} 
    \State $h_{z'} = h(z')$
    \If{$h_{z'} \leq H_1$}
    \State $H \leftarrow (h_{z'}, H_1, H_2, H_3, H_4)$
    \State $p_\text{succ} \leftarrow (1 - c_p) \cdot p_\text{succ} + c_p$
    \If{$p_\text{succ} > p_\text{thre}$}\label{l:bad1}
    \State $p \leftarrow (1 - c_c) \cdot p$, $c_\text{cov} = c_\text{cov+} (1 - c_c \cdot (2 - c_c))$\label{l:bad2}
    \Else
    \State $p \leftarrow (1 - c_c) \cdot p + \sqrt{c_c \cdot (2 - c_c)} \frac{z' - z}{\sigma}$, 
    $c_\text{cov} = c_\text{cov+}$
    \EndIf
    \State $w = A_\text{inv} \cdot p$
    \State $a = \sqrt{1 - c_\text{cov}}$,
    $b = \frac{\sqrt{1 - c_\text{cov}}}{\norm{w}^2} \left( \sqrt{1 + \frac{c_\text{cov} }{1 - c_\text{cov}} \norm{w}^2} - 1 \right)$
    \State $A \leftarrow a \cdot A + b \cdot (A \cdot w) \cdot w^\mathrm{T}$,
    $A_\text{inv} \leftarrow \frac{1}{a} \cdot A_\text{inv} - \frac{b}{a^2 + a\cdot b\cdot \norm{w}^2} \cdot w \cdot (w^\mathrm{T} A_\text{inv})$
    \State $\sigma \leftarrow \sigma \cdot c$,
    $z \leftarrow z'$,
    $n_\text{succ} \leftarrow n_\text{succ} + 1$
    \Else
    \State $p_\text{succ} \leftarrow (1 - c_p) \cdot p_\text{succ}$    
    \If{$h_{z'} > H_5$ and $p_\text{succ} \leq p_\text{thre}$}
    \State $w = A_\text{inv} \cdot \frac{z' - z}{\sigma}$
    \State $c_\text{cov} = c_\text{cov-}$ \textbf{if} $c_\text{cov-} (2 \cdot \norm{w}^2 - 1) \leq 1$ \textbf{else} $c_\text{cov} = \frac{1}{2 \cdot \norm{w}^2 - 1}$
    \State $a = \sqrt{1 + c_\text{cov}}$,
    $b = \frac{\sqrt{1 + c_\text{cov}}}{\norm{w}^2} \left( \sqrt{1 - \frac{c_\text{cov} }{1 + c_\text{cov}} \norm{w}^2} - 1 \right)$
    \State $A \leftarrow a \cdot A + b \cdot (A \cdot w) \cdot w^\mathrm{T}$,
    $A_\text{inv} \leftarrow \frac{1}{a} \cdot A_\text{inv} - \frac{b}{a^2 + a\cdot b\cdot \norm{w}^2} \cdot w \cdot (w^\mathrm{T} A_\text{inv})$
    \EndIf
    \State $\sigma \leftarrow \sigma \cdot c^{-1/4}$
    \EndIf\label{alg:aes:esend}
    \State $\sigma \leftarrow \sigma \cdot \frac{\norm{A}_F}{\sqrt{\ell}}, A \leftarrow A \cdot \frac{\sqrt{\ell}}{\norm{A}_F}, A_\text{inv} \leftarrow A_\text{inv} \cdot \frac{\norm{A}_F}{\sqrt{\ell}}, p \leftarrow p \cdot \frac{\sqrt{\ell}}{\norm{A}_F}$ every $\ell$ iterations\label{l:balance}
    \State \textbf{break} \textbf{if} $\sigma < \bar{\sigma}_{\min}$
    \EndWhile
    \State \Return $z$, $\max(\sigma, \bar{\sigma}_{\min})$, $A$
  \end{algorithmic}
\end{algorithm}

We shared the strategy parameter $\theta = (\sigma, A)$ over oracle calls.
Here, we implicitly assumed that the objective function $h$ of the current oracle call and that of the last oracle call are similar because the changes in $x_t$ and $y_t$ are small if $\eta$ is small. Then, reusing the strategy parameter of the last oracle call reduced the need for its adaptation time.

We ran (1+1)-CMA-ES until it improved the solution $\tau_\text{es} \cdot \ell + \tau_\text{es}'$ times. 
The reason for this procedure is described below.
Because the step size is maintained such that the probability of generating a successful solution is approximately 1/5, the algorithm runs approximately $T = 5 \cdot (\tau_\text{es} \cdot \ell + \tau_\text{es}')$ iterations. 
It was shown in \cite{morinaga2019generalized} that the expected runtime $\E[T_\epsilon]$ of (1+1)-ES with the simplified 1/5-success rule is $\Theta(\log(1/\epsilon))$ on strongly convex functions with Lipschitz continuous gradients and their strictly increasing transformations. Moreover, the scaling of the runtime with respect to dimension $\ell$ is $\Theta(\ell)$ on general convex quadratic functions \cite{morinaga2021}.
Therefore, we expect that $T$ iterations of (1+1)-CMA-ES approximates $\mathcal{M}$ with $\epsilon \in \exp(-\Theta(T / \ell)) = \exp(-\Theta(1))$.
The reason that we count the number of successful iterations instead of the number of total iterations is to avoid producing no progress because of a bad initialization of each oracle call. 

Another optional stopping condition is $\sigma < \bar{\sigma}_{\min}$ for a given minimal step size $\bar{\sigma}_{\min} \geq 0$. Once $\sigma$ reaches $\bar{\sigma}_{\min}$, \Cref{alg:es} returns $\sigma = \bar{\sigma}_{\min}$. Then, the next $\mathcal{M}$ call starts with $\sigma = \bar{\sigma}_{\min}$ and it is expected to stop after a few iterations. That is, if $\sigma$ for $\mathcal{M}_x$ reaches $\bar{\sigma}_{\min}$ while $\sigma > \bar{\sigma}_{\min}$ for $\mathcal{M}_y$, \Cref{alg:lr} spends more $f$-calls for $\mathcal{M}_y$ than for $\mathcal{M}_x$, and vice versa. 

Based on our preliminary experiments, we set $\tau_\text{es} = \tau_\text{es}' = 5$ as their default values. If they are set greater, we expect that (1+1)-CMA-ES approximates condition \eqref{eq:amo-cond} with a smaller $\epsilon$. 

\subsection{Adversarial-SLSQP}

We also implemented the algorithm with a sequential least squares quadratic programming (SLSQP) subroutine \cite{kraft1988software} to demonstrate the applicability of the proposed $\eta$ adaptation mechanism.
This was a first-order approach, which required access to $\nabla f$.
Unlike Adversarial-CMA-ES, no strategy parameter for SLSQP is shared over oracle calls. 
The maximum number of iterations is set to $\tau_\mathrm{slsqp} = 5$. 
We used the \texttt{scipy} implementation of SLSQP as $\mathcal{M}$ in \Cref{alg:lr}. 
We call this first-order approach \emph{Adversarial-SLSQP (ASLSQP)}.

\section{Numerical Evaluation}

Through experiments on test problems as described below, we confirmed the following hypotheses.
(A) Our implementations of the proposed approach, Adversarial-CMA-ES and Adversarial-SLSQP, performed as well as the theory implies.
(B) Our learning rate adaptation located a nearly optimal learning rate with little compromise of the objective function calls.
(C) Local strong convexity--concavity of the objective function is necessary for good min--max performance of the proposed approach.
(D) Existing coevolutionary approaches fail to converge even on a convex--concave quadratic problem.

\subsection{On Convex--Concave Quadratic Functions}\label{sec:ccqf}

To confirm (A) and (B), we ran Adversarial-CMA-ES and Adversarial-SLSQP on the following convex-concave quadratic function $f_1:\XX\times\YY\to\R$ with $n = m$:
\begin{equation}
  f_1(x, y) = \frac{a}{2} \norm{x}^2 + b \inner{x}{y} - \frac{c}{2} \norm{y}^2 ,
\end{equation}
where $a, c > 0$ and $b \in \R$.
The global min--max saddle point was located at $(x^*, y^*) = (0, 0)$.
The suboptimality error is $G_1(x, y) = \frac{a\cdot c + b^2}{2 a \cdot c} (\norm{x}^2 + \norm{y}^2)$. In this problem, we have $\alpha_H = \beta_H = 1$ and $\alpha_G = \beta_G = 1 + \frac{b^2}{a \cdot c}$; hence, for $\epsilon \ll 1$, we have $\eta^* \approx \bar{\eta} \approx \frac{a\cdot c}{a \cdot c + b^2}$. Moreover, for $\eta = \delta \cdot \eta^*$ for $\delta \in (0, 2)$, $\gamma = - \frac{a\cdot c}{a\cdot c+b^2} \delta ( 2 - \delta)$.
That is, \Cref{thm:conv} indicates that the runtime of the proposed approach with a fixed learning rate was proportional to $\left(1 + \frac{b^2}{a\cdot c}\right) \frac{1}{\delta (2 - \delta)}$.

The experimental settings were as follows.
We draw the initial solution $(x, y)$ uniform-randomly from $[-1, 5]^{m} \times [-1, 5]^{n}$.
The strategy parameters for Adversarial-CMA-ES are $\theta^x = (\sigma^x, A^x)$ and $\theta^y = (\sigma^y, A^y)$.
The step sizes $\sigma^x$ and $\sigma^y$ are initialized as one-fourth of the length of the initialization interval, that is, $\sigma^x = \sigma^y = 1.5$. The factors $A^x$ and $A^y$ are initialized by the identity matrix.
We used the default hyperparameter values described in the previous section. 
We omitted lines~\ref{l:px}--\ref{l:py} and lines~\ref{l:restart}--\ref{l:restart-end} of \Cref{alg:lr} (i.e., neither $P_x$ nor $P_y$ are given and $G_\text{tol} = 0$) in this experiment. The minimal learning rate was set to $\eta_{\min} = 10^{-4}$.
The minimal step sizes are set to $\bar{\sigma}_{\min}^x = \bar{\sigma}_{\min}^y = 0$.
We run $50$ independent trials for each setting, with the maximum number of $f$-calls of $10^{7}$.

\begin{figure}[t]
  \begin{subfigure}{0.5\hsize}%
    \includegraphics[width=\hsize]{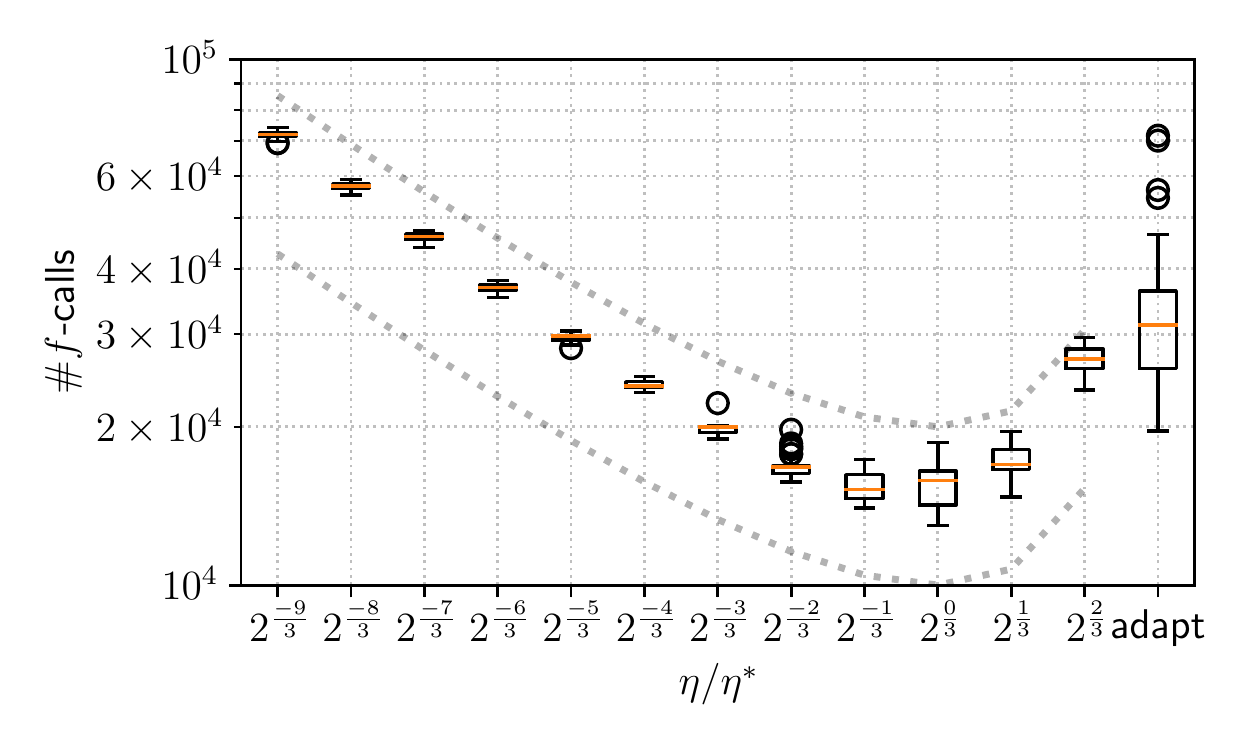}%
    \caption{Adversarial-CMA-ES}
  \end{subfigure}%
  \begin{subfigure}{0.5\hsize}%
    \includegraphics[width=\hsize]{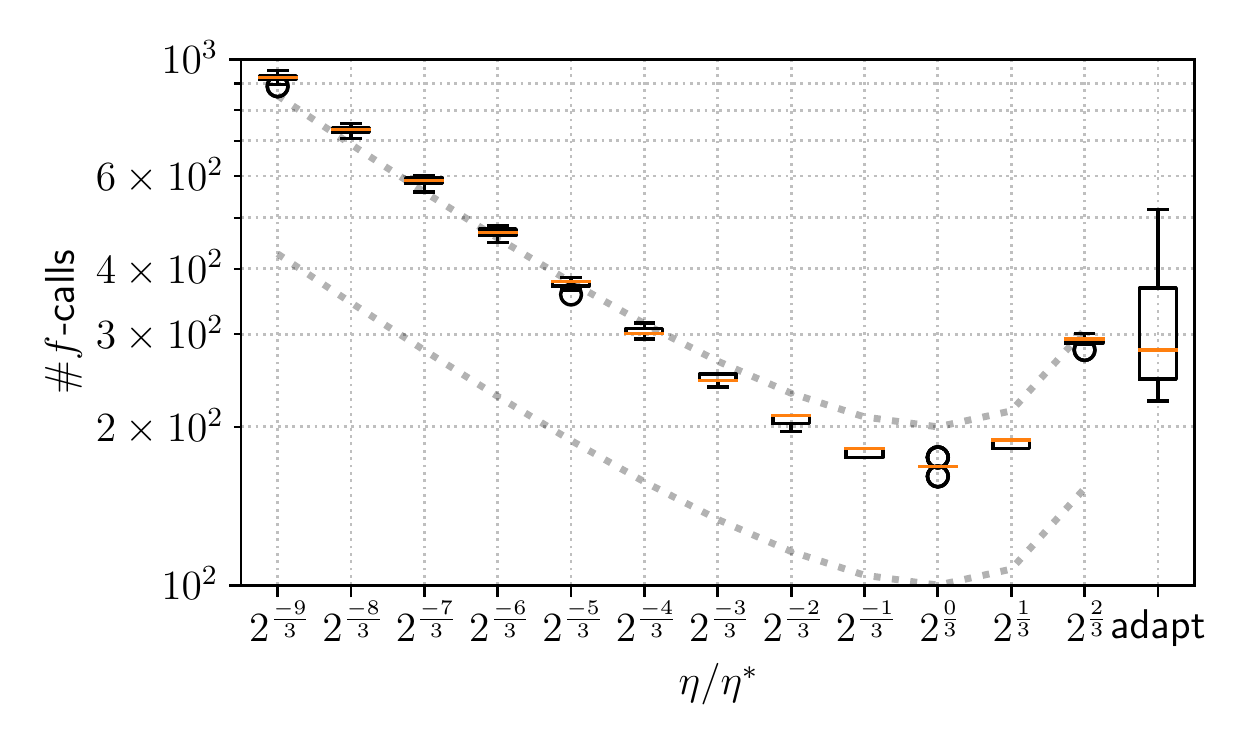}%
    \caption{Adversarial-SLSQP}
  \end{subfigure}%
  \caption{The number of $f$-calls until $G(x, y) \leq 10^{-5}$ is reached on $f_1$ with $n = m = 10$ and $a = b = c = 1$. (a) Adversarial-CMA-ES with $\eta$-adaptation (adapt) and fixed $\eta = \eta^* \times 2^{\frac{3-k}{3}}$ for $k = 1,\dots, 12$. (b) Adversarial-SLSQP with $\eta$-adaptation (adapt) and fixed $\eta = \eta^* \times 2^{\frac{3-k}{3}}$ for $k = 1,\dots, 12$. The dashed lines are proportional to $\frac{1}{(\eta/\eta^*)(2 - \eta/\eta^*)}$.}
  \label{fig:eta}
\end{figure}

\Cref{fig:eta} compares the proposed approaches with and without $\eta$-adaptation mechanism.
For fixed $\eta$ cases, we set $\eta$ to $\delta \cdot \eta^*$ with $\delta \in \{2^\frac{3-k}{3} : k = 1, \dots, 12\}$. We remark that for both algorithms, all the trials with $\eta = 2 \times \eta^*$ fail to converge, as implied by \Cref{thm:conv}. As expected, the runtimes of both algorithms with fixed $\eta$ were nearly proportional to $\frac{1}{(\eta/\eta^*)(2 - \eta/\eta^*)}$. The best $\eta$ is approximately $\eta^*$.
We conclude that our implementations closely approximate the oracle condition \eqref{eq:amo-cond} and that the proposed approach works as the theory implies.

The proposed approach with the $\eta$-adaptation mechanism succeed in converging toward the global min--max saddle point. 
Comparing the runtime of the $\eta$-adaptation mechanism and the best fixed $\eta = \eta^*$, we compromise the number of $f$-calls at most three times in the median case for both Adversarial-CMA-ES and Adversarial-SLSQP to adapt $\eta$.
There are also trials that required a few times more runtime than the median case. 
However, considering the difficulty in tuning $\eta$ in advance, we conclude that this $\eta$-adaptation mechanism is promising to waive the need for $\eta$ tuning in advance. 

\begin{figure}[t]
  \centering%
  \begin{subfigure}{0.33\hsize}%
    \centering%
    \includegraphics[width=\hsize]{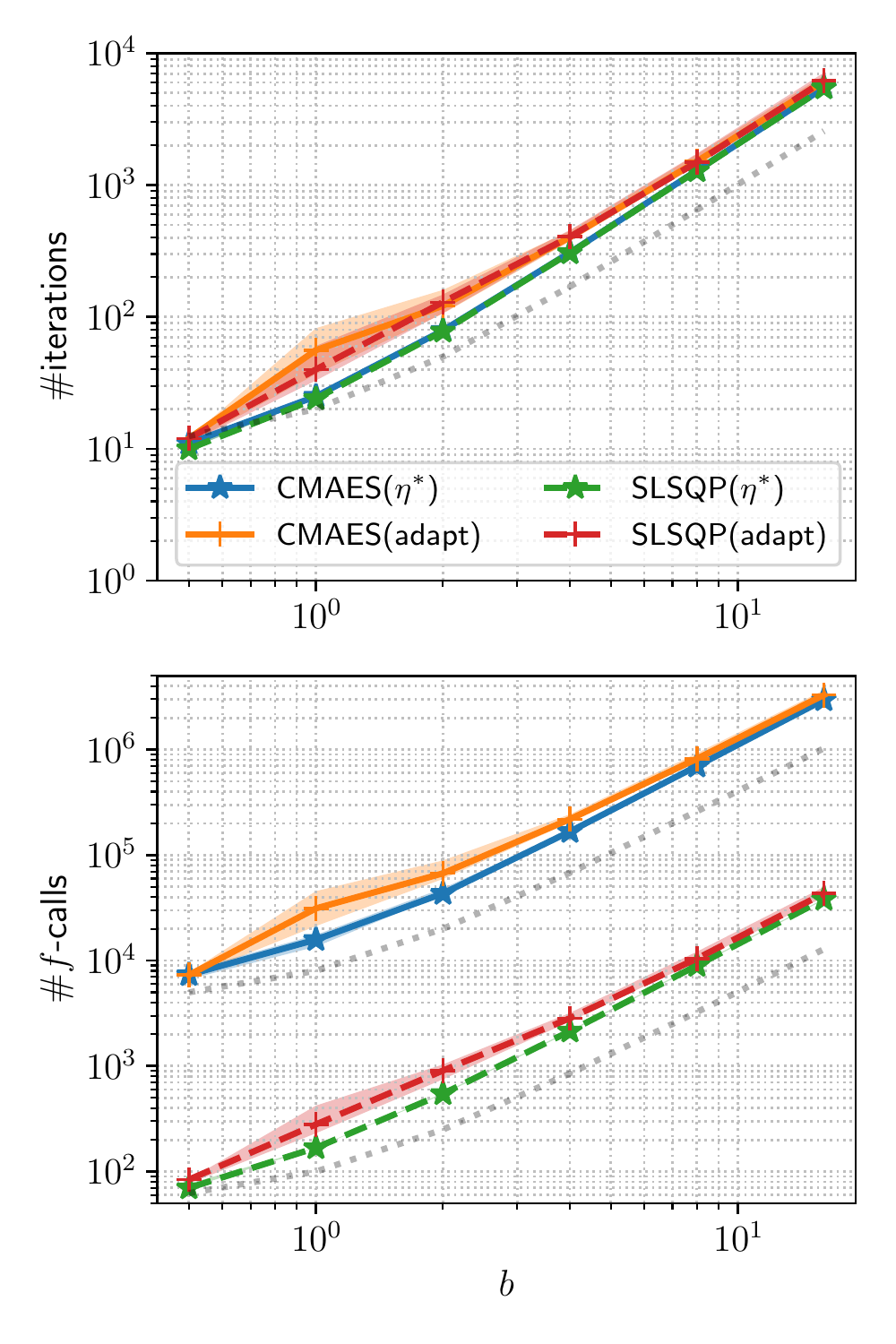}%
    \caption{Varying $b$\\
      $a = c = 1$ and $n = m = 10$}\label{fig:b1}%
  \end{subfigure}%
  \begin{subfigure}{0.33\hsize}%
    \centering%
    \includegraphics[width=\hsize]{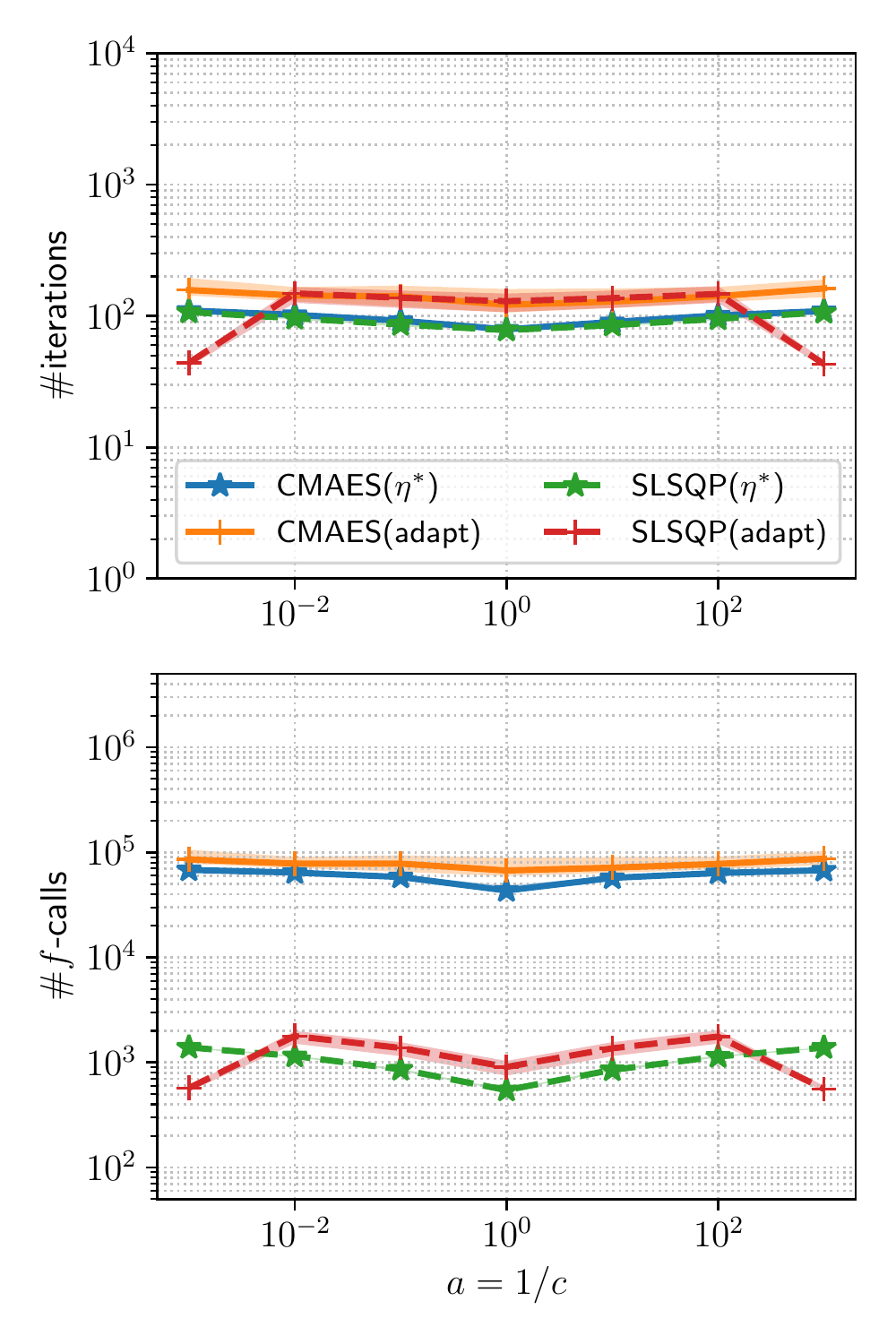}%
    \caption{Varying $\frac{a}{c}$\\
      $b = ac = 1$ and $n = m = 10$}\label{fig:b2}%
  \end{subfigure}%
  \begin{subfigure}{0.33\hsize}%
    \centering%
    \includegraphics[width=\hsize]{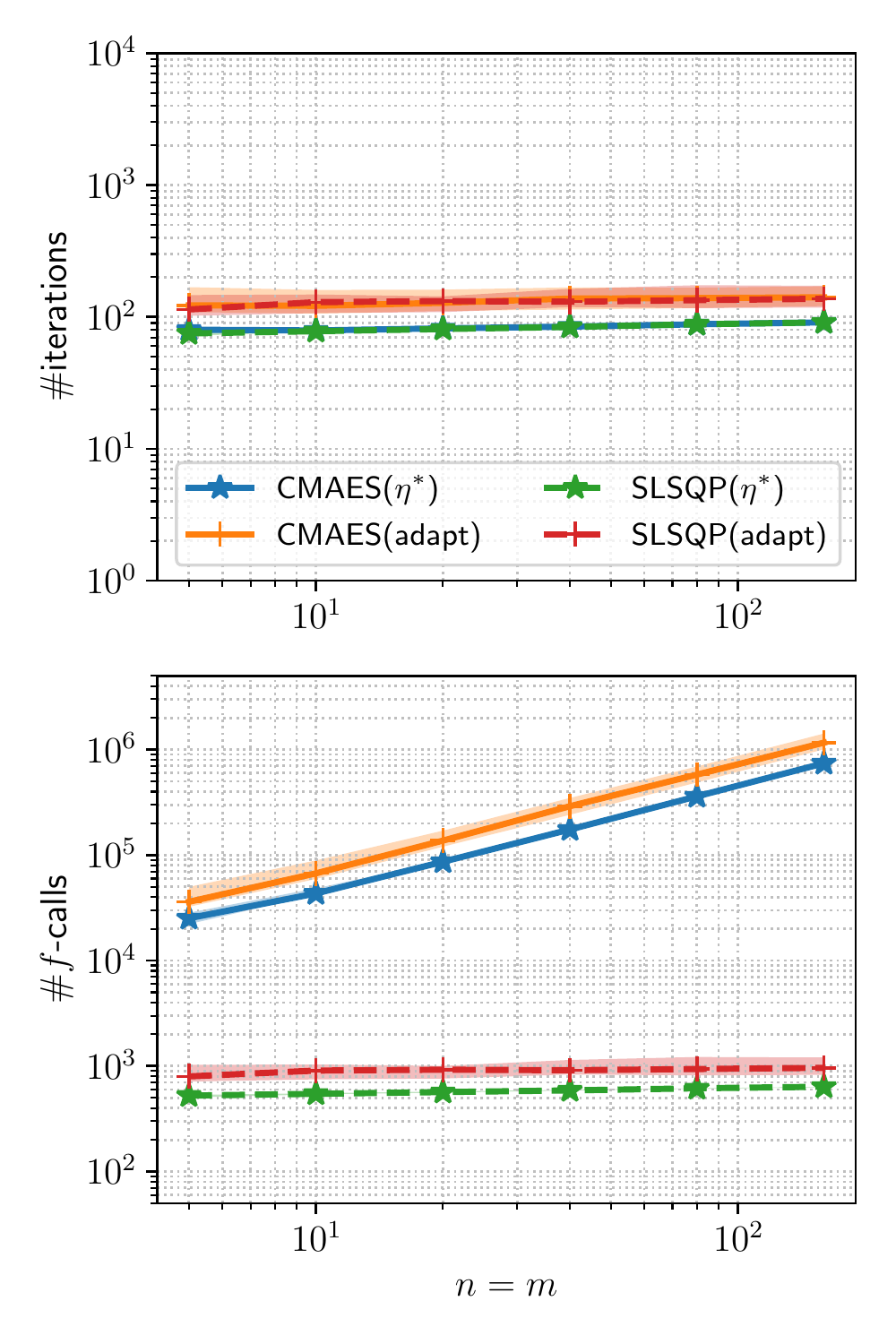}%
    \caption{Varying $n = m$\\
      fixed $a = c = 1$ and $b = 2$}\label{fig:dim}%
  \end{subfigure}%
  \caption{The number of iterations and the number of $f$-calls until $G_1(x, y) \leq 10^{-5}$ is reached on $f_1$.
    The solid lines indicate the median and the shaded areas indicate the $10$--$90$ percentile ranges.
    Dashed lines are proportional to $1 + \frac{b^2}{a\cdot c}$.}
\end{figure}

\Cref{fig:b1,fig:b2} show the runtime of the proposed approaches with and without $\eta$-adaptation for varying $b$ and for varying $a/c$. For the fixed $\eta$ case, we set $\eta = \eta^*$. It may be observed that the runtimes in terms of the number of iterations are proportional to $1 + \frac{b^2}{a \cdot c}$, as expected from \Cref{thm:conv}. Moreover, the number of iterations was largely the same for all algorithms, as they all approximate \eqref{eq:amo-cond} with $\epsilon \ll 1$. In contrast, the number of $f$-calls was different for the two algorithms. This is because Adversarial-CMA-ES is expected to spend approximately $5 (\tau_\mathrm{es} \times \ell + \tau_\mathrm{es}') = 275$ $f$-calls per oracle call, whereas Adversarial-SLSQP spends $\tau_\mathrm{slsqp} = 5$ $f$-calls. We remark that if one of the CMA-ES in Adversarial-CMA-ES (i.e., either $\mathcal{M}_{x}$ or $\mathcal{M}_{y}$) is replaced with SLSQP, the number of $f$-calls will be approximately halved. Therefore, it is advisable to use SLSQP, or another first-order approach, as an approximate minimization oracle if $\nabla f$ is available and cheap to compute.
\Cref{fig:dim} shows the scaling of the runtime with respect to the dimension $n = m$. The number of iterations did not depend on the search space dimension. The number of $f$-calls was also constant over varying $n = m$ for Adversarial-SLSQP. However, it was proportional to $n + m$ for Adversarial-CMA-ES.\footnote{We comment on the computational complexity of the algorithm. The bottleneck of the execution time of each iteration of \Cref{alg:lr} is an $\mathcal{M}_x$-call and an $\mathcal{M}_y$-call. The execution time for the $\eta$-adaptation was negligible. The time and space complexity of \Cref{alg:es} per $f$-call is $O(\ell^2)$, where $\ell$ is the search space dimension. Therefore, if the number of $f$-calls scales linearly in $n + m$, the execution time of Adversarial-CMA-ES scales as $O(n^3+m^3)$.} This is because the runtime of (1+1)-CMA-ES is proportional to the dimension, and iterations must be run proportional to the search space dimension to approximate \eqref{eq:amo-cond}. 

\subsection{Comparison with Baseline Approach}\label{sec:comp}

To confirm (C) and (D), we ran Adversarial-CMA-ES on the six test problems summarized in \Cref{tbl:func}. 
In all cases, the domain of the objective function is $\X \times \Y = [-1, 5]^m \times [-1, 5]^n$.
The function $f_2$ is globally strongly convex--concave, while $f_3$ is locally strongly convex--concave. The function $f_4$ is globally convex--concave but not strongly convex--concave. These functions exhibited a global min--max saddle point at $(x^*, y^*) = (0, 0)$ and $x^*$ was the global optimal solution to the worst-case objective $F(x) = f(x, \opty(x))$. The function $f_5$ was not strongly convex--concave, but the worst case $y$ is independent of $x$, and the optimal $x$ is constant over $y$ such that $\sum_{j=1}^{n} y_j > 0$. The optimal solutions $x^* = 0$ to the worst-case objective functions for $f_6$ and $f_7$ were not min--max saddle points. 

\begin{table}[t]
\centering
\caption{Definition of the test functions $f(x, y)$ and their worst-case variable $\opty(x)= \argmax_{y \in \Y} f(x, y)$}\label{tbl:func}%
\begin{tabular}{lll}
\toprule
& $f(x, y)$ & $\opty(x)$
\\
\midrule
$f_2$ & $\frac{1}{2} \norm{x}^2 + \frac{1}{m}\sum_{i=1}^{m}  x_i \sum_{j=1}^{n}y_j - \frac{1}{2}\norm{y}^2$
& 
$\left(\frac{1}{m}\sum_{i=1}^{m} x_i \right) \mathbf{1}$
\\
$f_3$ & $\frac{1}{2} \min\left[\norm{x}^2,\ \norm{x - 4\cdot\mathbf{1}}^2\right] + \frac{1}{m}\sum_{i=1}^{m}  x_i \sum_{j=1}^{n}y_j - \frac{1}{2}\norm{y}^2$ 
& 
$\left(\frac{1}{m}\sum_{i=1}^{m}  x_i \right) \mathbf{1}$
\\
$f_4$ & $\frac{1}{2} \norm{x}^2 + \frac{1}{m}\sum_{i=1}^{m}  x_i \sum_{j=1}^{n}y_j - \left(\frac{1}{2}\norm{y}^2\right)^2 $
&
$\left(\frac{1}{m\cdot n}\sum_{i=1}^{m} x_i \right)^\frac{1}{3} \mathbf{1}$
\\
$f_5$ & $\frac{1}{m}\norm{x_i} \sum_{j=1}^{n}y_j$
&
$5 \cdot \mathbf{1}$
\\
$f_6$ & $\frac{1}{m}\sum_{i=1}^{m}x_i \sum_{j=1}^{n}y_j - \frac{1}{2}\norm{y}^2 $
&
$\left(\frac{1}{m}\sum_{i=1}^{m} x_i \right) \mathbf{1}$
\\
$f_7$ & $\frac{1}{2} \norm{x}^2 + \frac{1}{m}\sum_{i=1}^{m}  x_i \sum_{j=1}^{n}y_j $
&
$\begin{cases} 5 \cdot \mathbf{1} & \sum_{i=1}^{m}x_i \geq 0 \\ - 1 \cdot \mathbf{1} & \sum_{i=1}^{m}x_i < 0\end{cases}$
\\
\bottomrule
\end{tabular}
\end{table}

The experimental settings were as follows. 
We ran Adversarial-CMA-ES with and without sampling distributions $P_x$ and $P_y$.
For the distributions $P_x$ and $P_y$, uniform distributions over $\X$ and $\Y$ are used. 
Moreover, we use the same initialization as in \Cref{sec:ccqf}. The minimal learning rate is $\eta_{\min} = 10^{-4}$. 
The minimal step sizes were set to $\bar{\sigma}_{\min}^x = \bar{\sigma}_{\min}^y = 0$.
The restart was not performed, that is, $G_\text{tol} = 0$. 
The boundary constraint was handled using the mirroring technique, that is, the domain was virtually extended to $\XX\times\YY$ by defining the function value $f(x, y)$ for $(x, y) \notin \X\times\Y$ by $f(T_\X(x), T_\Y(y))$, where $T_\X$ and $T_\Y$ map each coordinate to $U - \abs{\text{mod}(x - L, 2(U-L)) - (U-L)}$, where $U = 5$ and $L = -1$ denote the upper and lower bounds of each coordinate.
The output of (1+1)-CMA-ES ($\mathcal{M}_x$ and $\mathcal{M}_y$) is repaired into the feasible domain by applying $T_\X$ and $T_\Y$.
We compare the results with those of the naive baseline approach, referred to as CMA-ES($N_y$). 
We sampled $N_y = 10$ or $100$ points uniform-randomly in $\Y$, and they were denoted as $y^{k}$ for $k = 1, \dots, N_y$.
The approximate worst-case objective was defined as $F_{N_y}(x) = \max_{1 \leq k \leq N_y} f(x, y^k)$. 
Then, we solve $F_{N_y}$ with (1+1)-CMA-ES (\Cref{alg:es}) using mirroring boundary constraint handling.
These algorithms are run $10$ times with different initial solutions.
We also compared two coevolutionary approaches, MMDE \cite{mmde2018} and COEVA \cite{Al-Dujaili2019lego}. 
These approaches were implemented based on the Python code provided by the authors of \cite{Al-Dujaili2019lego}.

\begin{figure}[t]
  \centering%
  \begin{subfigure}{0.33\hsize}%
    \centering%
    \includegraphics[width=\hsize]{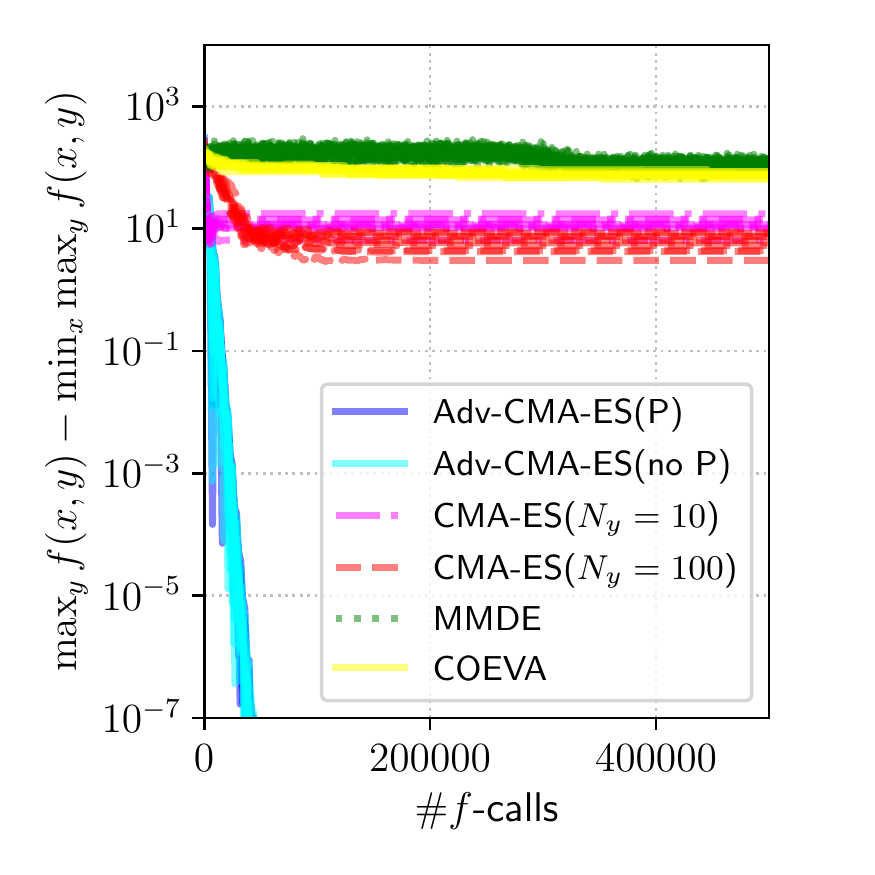}%
    \caption{$f_2$}%
  \end{subfigure}%
  \begin{subfigure}{0.33\hsize}%
    \centering%
    \includegraphics[width=\hsize]{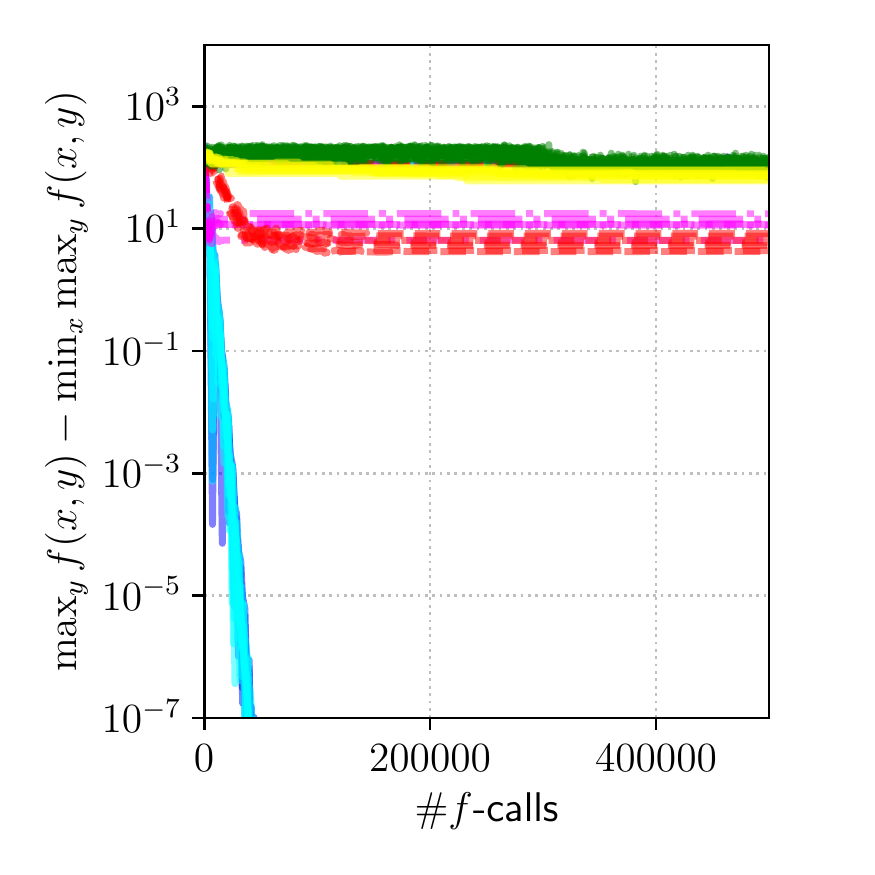}%
    \caption{$f_3$}%
  \end{subfigure}%
  \begin{subfigure}{0.33\hsize}%
    \centering%
    \includegraphics[width=\hsize]{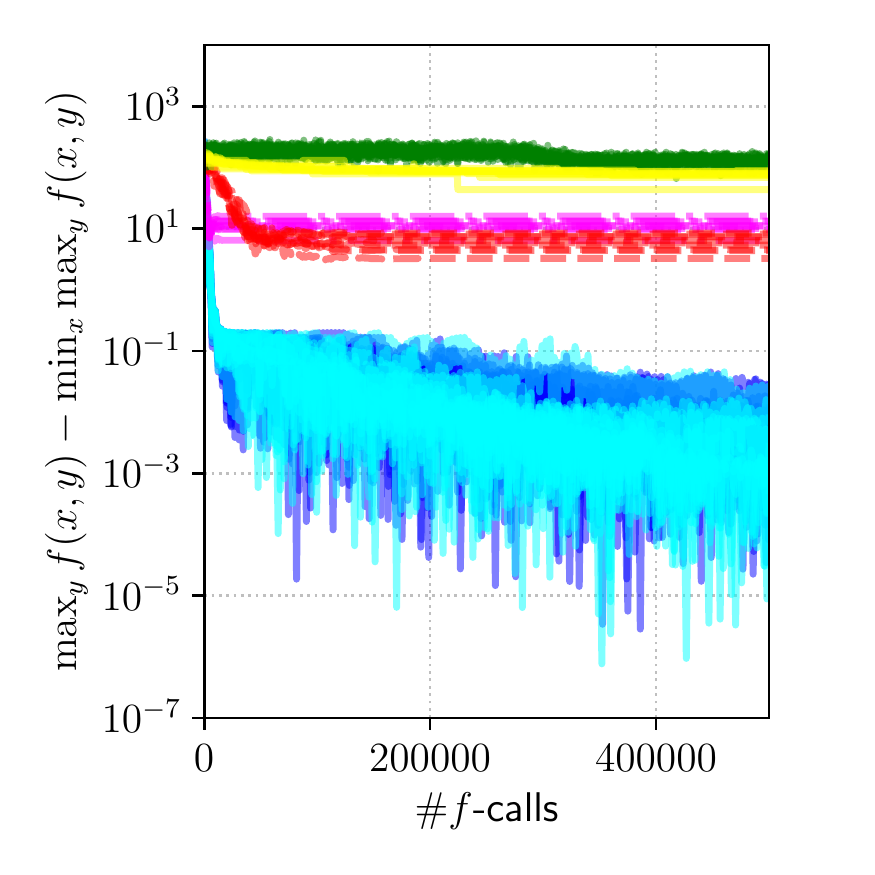}%
    \caption{$f_4$}%
  \end{subfigure}%
  \\
  \begin{subfigure}{0.33\hsize}%
    \centering%
    \includegraphics[width=\hsize]{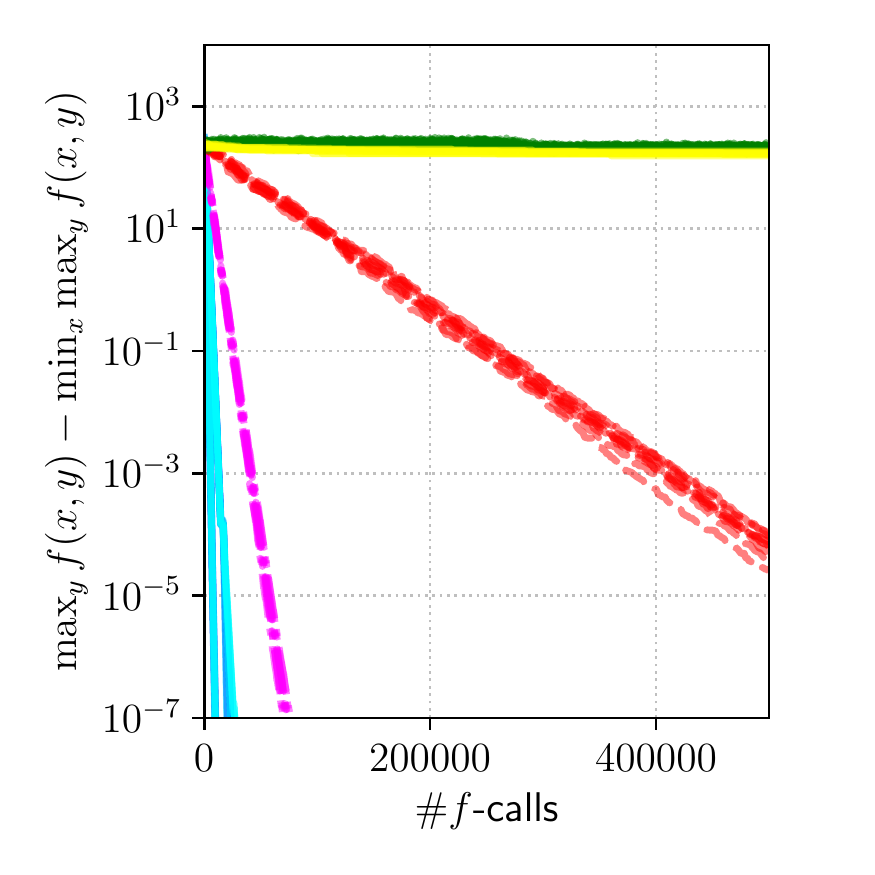}%
    \caption{$f_5$}%
  \end{subfigure}%
  \begin{subfigure}{0.33\hsize}%
    \centering%
    \includegraphics[width=\hsize]{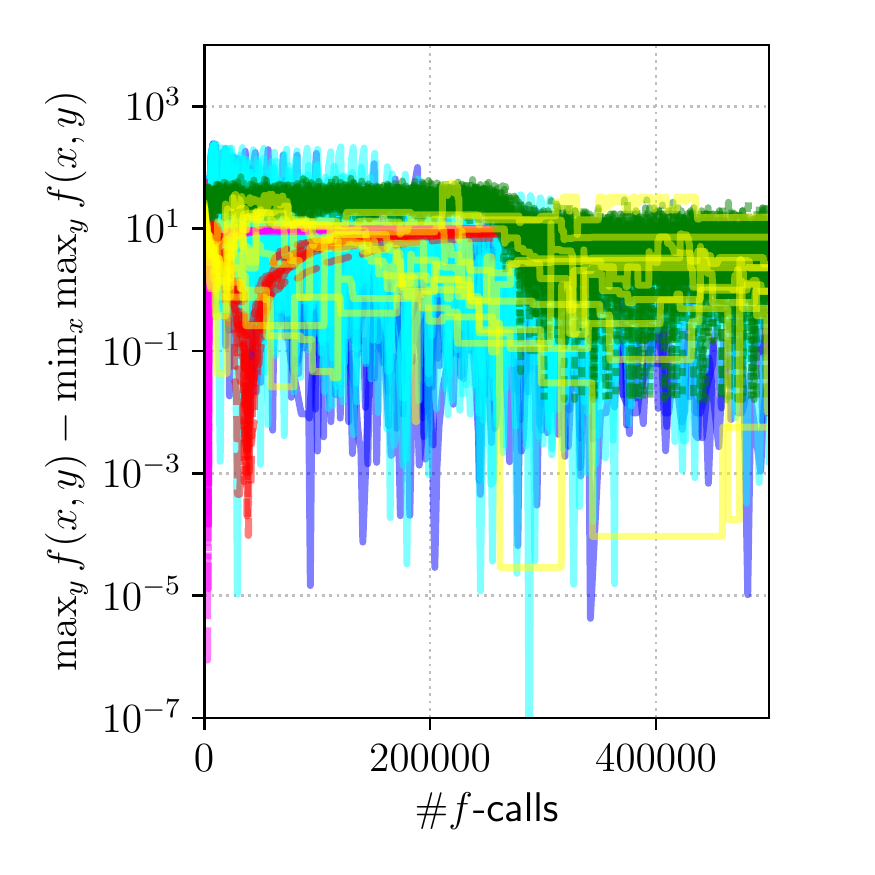}%
    \caption{$f_6$}%
  \end{subfigure}%
   \begin{subfigure}{0.33\hsize}%
    \centering%
    \includegraphics[width=\hsize]{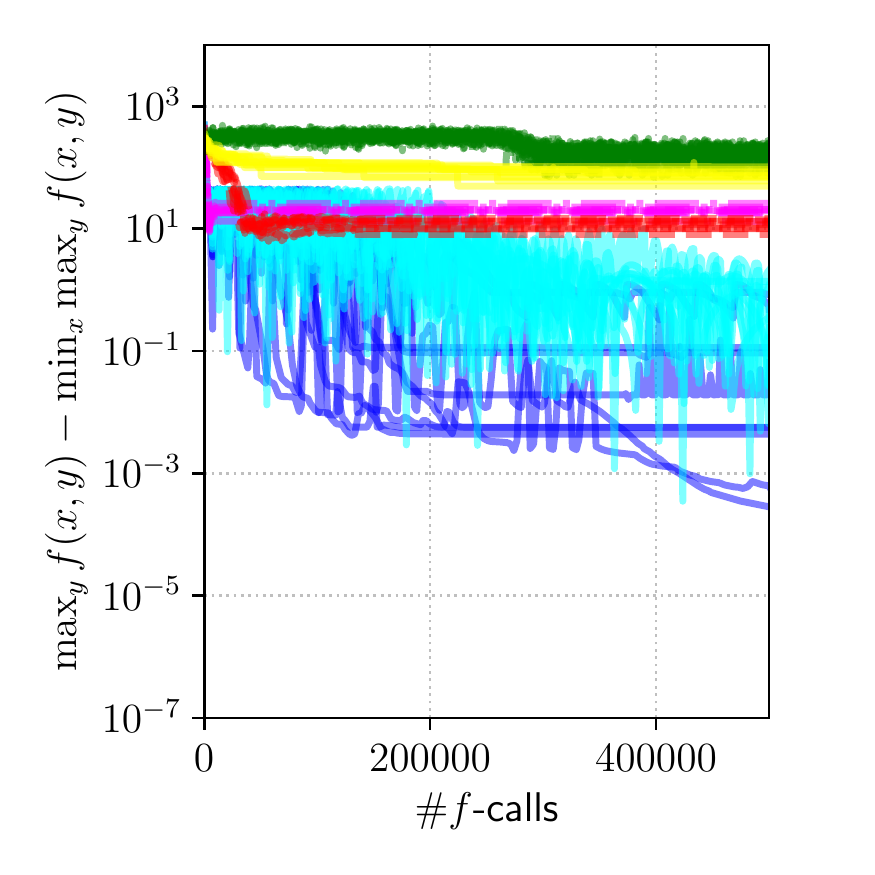}%
    \caption{$f_7$}%
  \end{subfigure}%
\caption{Results of 10 independent runs of Adversarial-CMA-ES with and without sampling distribution (denoted as Adv-CMA-ES(P) and Adv-CMA-ES(no P), respectively), CMA-ES($N_y = 10$), CMA-ES($N_y = 100$), MMDE, and COEVA. The search space dimension is $m = 50$ and $n = 20$ for all cases.}\label{fig:ex2}
\end{figure}

\Cref{fig:ex2} shows the results of 10 independent trials of Adversarial-CMA-ES, CMA-ES($N_y = 10$), CMA-ES($N_y = 100$), MMDE, and COEVA.
Adversarial-CMA-ES succeeds in converging the global min--max saddle point on $f_2$, $f_3$, and $f_6$. The functions $f_2$ and $f_3$ were locally strongly convex--concave functions, and Adversarial-CMA-ES performed well, as expected. The existing coevolutionary approaches, as well as CMA-ES($N_y$), failed to converge on these problems. 
The benchmark problems used to evaluate the performance of existing coevolutionary approaches \cite{Branke2008ppsn,mmde2018,Zhou.2010} are rather low-dimensional problems ($m \leq 2$ and $n \leq 2$). These approaches do not work well on higher-dimensional problems and perform worse than the simple baseline, CMA-ES($N_y$).
CMA-ES($N_y$) tends to the global optimal point on $f_5$. This is because the optimal $x^*$ is optimal for approximate worst-case functions provided that there exists $y$ in $N_y$ samples such that $\sum_{i=1}^{n} y_i > 0$ holds. 
In contrast, no approach succeeded in converging toward the global optimum of the worst-case function on $f_4$, $f_6$, and $f_7$. From these results, we conclude that local strong convexity--concavity is an important factor for the convergence of Adversarial-CMA-ES. These results reveal the limitations of Adversarial-CMA-ES and the difficulty of locating the solution to the worst-case objective if it is not a min--max saddle point.

\section{Application to Robust Berthing Control}\label{sec:berth}

In this section, we analyze the application of Adversarial-CMA-ES to robust berthing control tasks under model uncertainty. 

\subsection{Problem Description}

\begin{figure}[t]
    \centering
    \begin{subfigure}{0.5\hsize}
    \centering
    \includegraphics[width=\hsize]{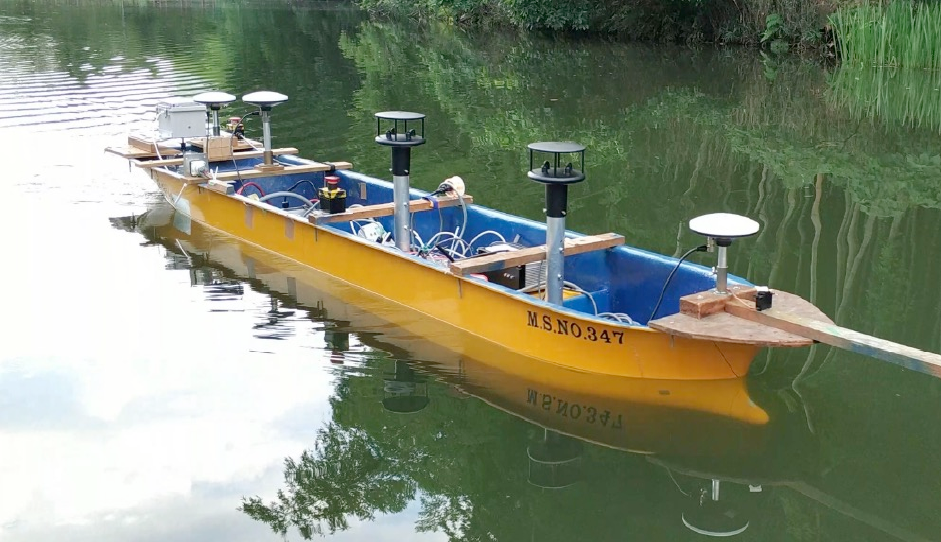}%
    \caption{Photo of 3 m model ship of ESSO OSAKA}%
    \label{fig:esso_photo}%
    \end{subfigure}%
    \begin{subfigure}{0.5\hsize}%
    \includegraphics[width=\hsize]{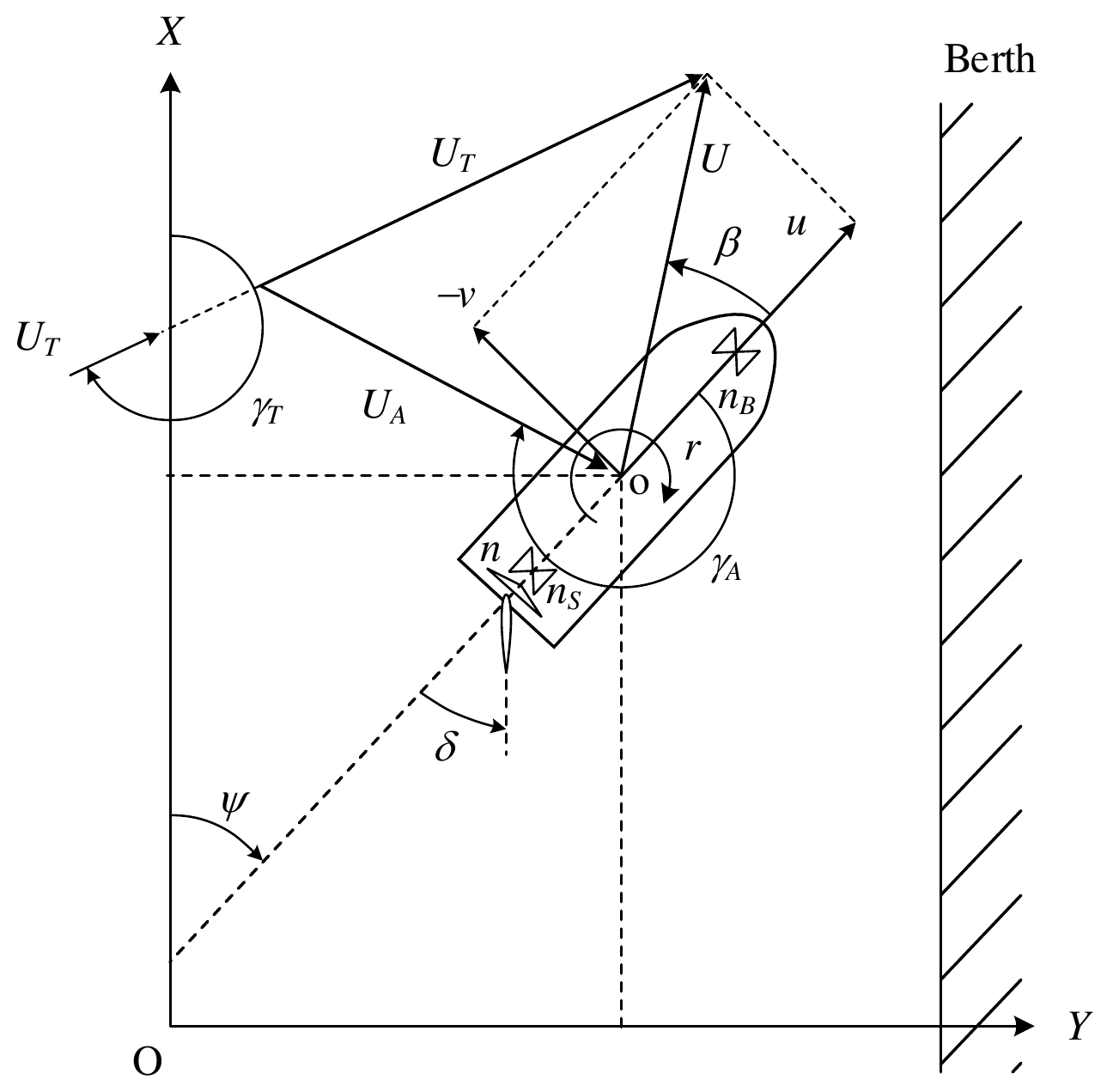}%
    \caption{Coordinate systems}%
    \label{fig:coordinate}%
    \end{subfigure}%
    \caption{Control target: 3 m model ship of ESSO OSAKA}
    \label{fig:esso_osaka}
\end{figure}

The objective of our robust berthing control task is to obtain a controller that controls a subject ship toward the target position near a berth while avoiding collision to the berth under the worst situation in the predefined uncertainty set.
The control target was a 3 m model ship representing the vessel MV ESSO OSAKA (\Cref{fig:esso_photo}).
The state variables $s = (X, u, Y, v_m, \psi, r) \in \R^6$ and the control signal $a = (\delta, n_\mathrm{p}, n_\mathrm{BT}, n_\mathrm{ST}) \in \R^4$ were as described in \Cref{fig:coordinate}.
The controller $u_x : s \mapsto a$ was modeled by a neural network with $m = 99$ dimensional parameter vector $x$. 
The objective function $f(x, y)$ measures the distance between the target position and the final position of the subject ship after a pre-defined control time using the controller $u_x$, penalized by the risk of the collision to the berth, under an uncertainty parameter $y \in \Y$ described below. 
The details of the controller and the objective function are explained in \Cref{apdx:berth:problem}.

The wind conditions $y^{(A)}$ and the model coefficients $y^{(B)}$ with respect to the wind forces are treated as the uncertain factors $y = (y^{(A)}, y^{(B)})$. 
The following three situations are considered. 
(A) The state equation model is accurately modeled, but the wind conditions are uncertain.
In this situation, the uncertainty parameters $y^{(A)} \in \Y_{A}$ represent the wind velocity $U_T$ [m/s] and the wind direction $\gamma_T$ [rad], and their feasible values are in $\Y_{A} = [0, 0.5] \times [0, 2\pi]$. The model coefficients $y^{(B)}$ are set to the same values as in \cite{miyauchi2021}, denoted by $y_\mathrm{est}^{(B)}$.
(B) Wind conditions are known, but the state equation model is uncertain. The coefficients in the state equation model for the effect of the wind force were derived in \cite{Fujiwara1998} using regression of wind tunnel experiment data, and we consider them to be relatively inaccurate. The uncertainty parameters $y^{(B)}$ consist of $10$ coefficients for the wind force. The feasible domain $\Y_{B}$ is constructed to include the coefficient used in \cite{miyauchi2021}, that is, $y_\mathrm{est}^{(B)} \in \Y_{B}$. For each variable, the feasible values are defined by the interval. The interval of the $i$th component of $y^{(B)}$, denoted by $[y_\mathrm{est}^{(B)}]_i$, is set to $[0.9 \cdot [y_\mathrm{est}^{(B)}]_{i}, 1.1 \cdot [y_\mathrm{est}^{(B)}]_{i}]$ for all $i = 1,\dots, 10$. The other model coefficients are set to the same values as \cite{miyauchi2021} and the wind condition is set to $y_\mathrm{est}^{(A)} = (1.5\pi, 0.5)$, meaning that the velocity of wind blowing orthogonally from the sea to the berth is $0.5$ [m/s]. (C) Wind conditions are unknown, and the model coefficients are uncertain. In this situation, $y$ is composed of the uncertainty parameters $y^{(A)}$ and $y^{(B)}$, and the feasible values are set to $\Y_{C} = \Y_{A} \times \Y_{B}$. 

\subsection{Experiment Setting}

We ran Adversarial-CMA-ES and CMA-ES($N_y = 100$) on $\Y_{A}$, $\Y_{B}$, and $\Y_{C}$.\footnote{CMA-ES($N_y = 10$), MMDE and COEVA tested in \Cref{sec:comp} were omitted from the comparison based on our preliminary experiments. The worst-case performance of CMA-ES($N_y = 10$) were worse than the worst-case performance of CMA-ES($N_y = 100$) on our problems. The worst-case performance of MMDE and COEVA were not competitive to the other approaches as demonstrated in \Cref{fig:ex2}.} 
As baselines, we run (1+1)-CMA-ES on $f(x, y_\text{fix})$ under two different situations for $y_\text{fix} \in \Y$.
The first situation was $y_\text{fix} = (y_\mathrm{no}^{(A)}, y_\mathrm{est}^{(B)})$, where $y_\mathrm{no}^{(A)} = (0, 0)$ corresponds to no wind disturbance, and the second situation was $y_\text{fix} = (y_\mathrm{est}^{(A)}, y_\mathrm{est}^{(B)})$, where $y_\mathrm{est}^{(A)} = (1.5\pi, 0.5)$ reflects our prior knowledge that such a wind is difficult to handle for avoiding collision with the berth. Note that (1+1)-CMA-ES on $f(x, y_\text{fix})$ corresponds to CMA-ES($N_y = 1$) with $y^1 = y_\text{fix}$.
Each algorithm was run 20 times independently with random initialization of $x$ and $y$.
The search space for $x$ and $y$ was scaled to $\X = [-1, 1]^{m}$ and $\Y = [-1, 1]^{n}$. 
The box constraint was treated using the mirroring technique described in \Cref{sec:comp}. 
The initial solution $(x, y)$ was drawn uniform-randomly from $\X \times \Y$.
For CMA-ES($N_y = 100$), $y^k$ for $k = 1, \dots, N_y$ were uniform-randomly generated. 
The step sizes $\sigma^x$ and $\sigma^y$ are initialized as one-fourth of the length of the initialization interval. 
The factors $A^x$ and $A^y$ are initialized by the identity matrix. 
The minimal step size is $\bar{\sigma}_{\min} = 10^{-8}$ for both $\sigma^x$ and $\sigma^y$. 
We set $G_{\mathrm{tol}} = 10^{-6}$ and $d^y_{\min} = \bar{\sigma}_{\min} \times \sqrt{n}$ for Adversarial-CMA-ES.
The $f$-call budget was $10^{6}$.

For Adversarial-CMA-ES, we used the restart strategy proposed in \Cref{alg:lr}.
The output of Adversarial-CMA-ES follows \Cref{alg:lr}.
For CMA-ES($N_y$), when the termination condition $\sigma < \bar{\sigma}_{\min}$ was satisfied, the candidate solution was recorded and the algorithm was re-started until it exhausted the $f$-call budget. 
Note that $y^k$ ($k = 1, \dots,N_y$) were not resampled. 
The output of CMA-ES($N_y$) is determined as follows: Let $\{x^1, \dots, x^r\}$ be the set of recorded candidate solutions and the solution obtained at the end of the run. We then selected $x = \argmin_{i=1,\dots,r} \max_{k=1,\dots,N_y} f(x^i, y^k)$ as the output of CMA-ES($N_y$). 

The obtained solutions were evaluated as follows. 
Because the ground truth worst-case objective function value $F(x) = \max_{y \in \Y} f(x, y)$ for a given $x$ is unknown, we perform numerical optimization to approximate $F(x)$. We ran (1+1)-CMA-ES for $500 \times n$ iterations to obtain a local maximal point $y$ of $f(x, y)$. As the objective is expected to have multiple local optima, we repeat it 100 times with different initial search points $y$. The initialization of (1+1)-CMA-ES is as described above. 

\subsection{Results and Discussion}

\begin{figure}
    \centering
    \begin{subfigure}{\hsize}%
    \centering
    \includegraphics[height=\hsize,angle=270]{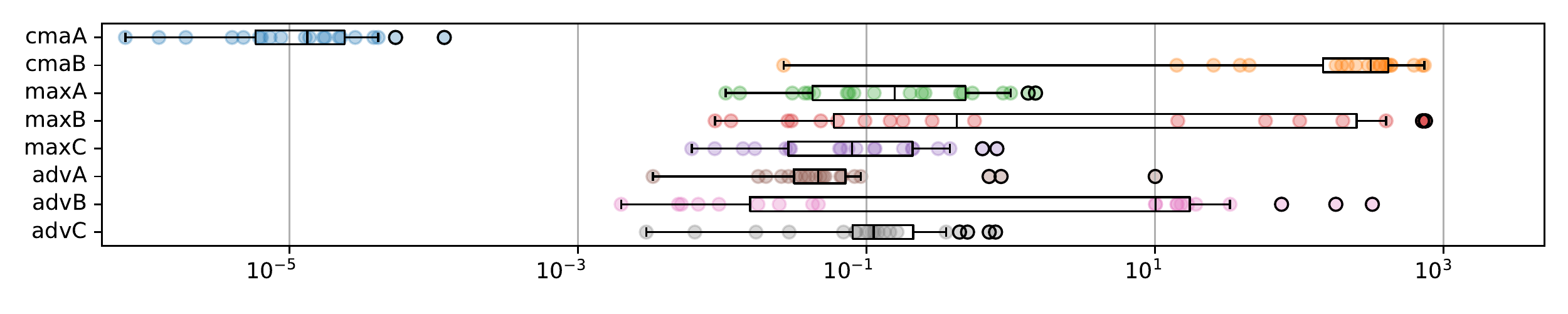}%
    \caption{$f(x, (y_\text{no}^{(A)}, y_\text{est}^{(B)}))$}\label{fig:berth:a0}%
    \end{subfigure}%
    \\
    \begin{subfigure}{\hsize}%
    \centering
    \includegraphics[height=\hsize,angle=270]{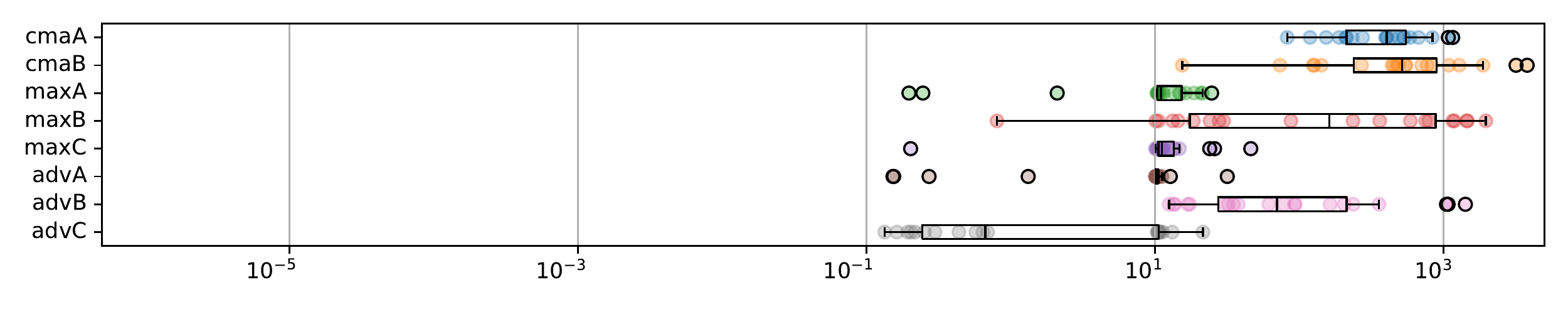}%
    \caption{$\max_{y^{(A)} \in \Y_{A}} f(x, (y^{(A)}, y_\text{est}^{(B)}))$}\label{fig:berth:a}%
    \end{subfigure}%
    \\
    \begin{subfigure}{\hsize}%
    \centering
    \includegraphics[height=\hsize,angle=270]{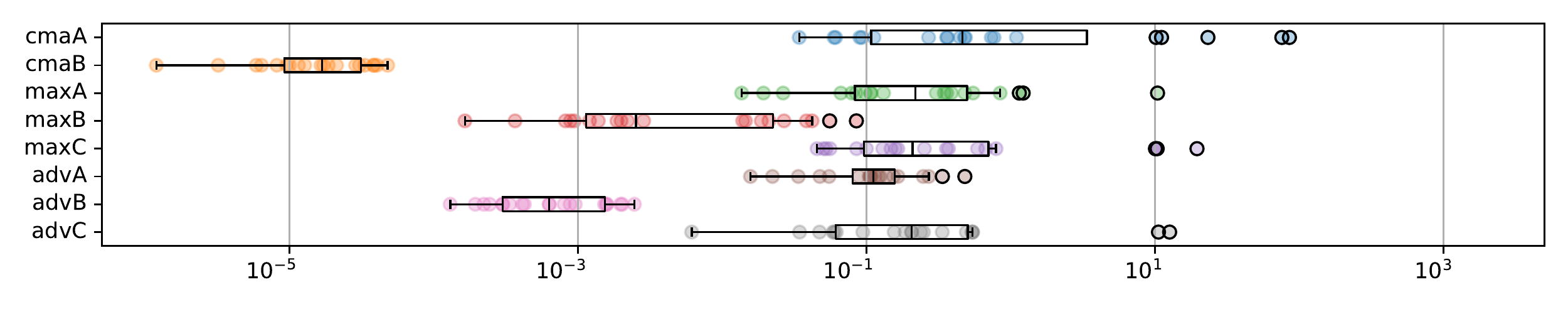}%
    \caption{$f(x, (y_\text{est}^{(A)}, y_\text{est}^{(B)}))$}\label{fig:berth:b0}%
    \end{subfigure}%
    \\
    \begin{subfigure}{\hsize}%
    \centering
    \includegraphics[height=\hsize,angle=270]{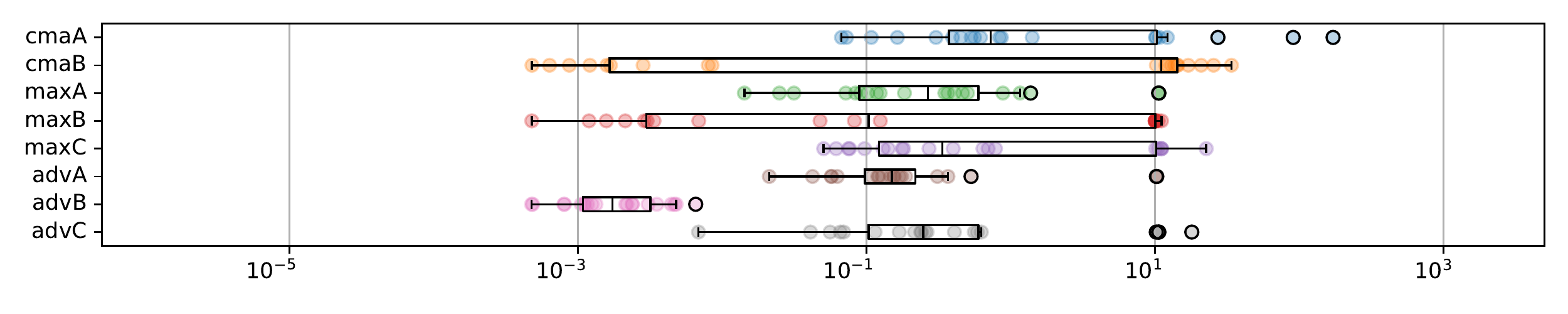}%
    \caption{$\max_{y^{(B)} \in \Y_{B}} f(x, (y_\text{est}^{(A)}, y^{(B)}))$}\label{fig:berth:b}%
    \end{subfigure}%
    \\
    \begin{subfigure}{\hsize}%
    \centering
    \includegraphics[height=\hsize,angle=270]{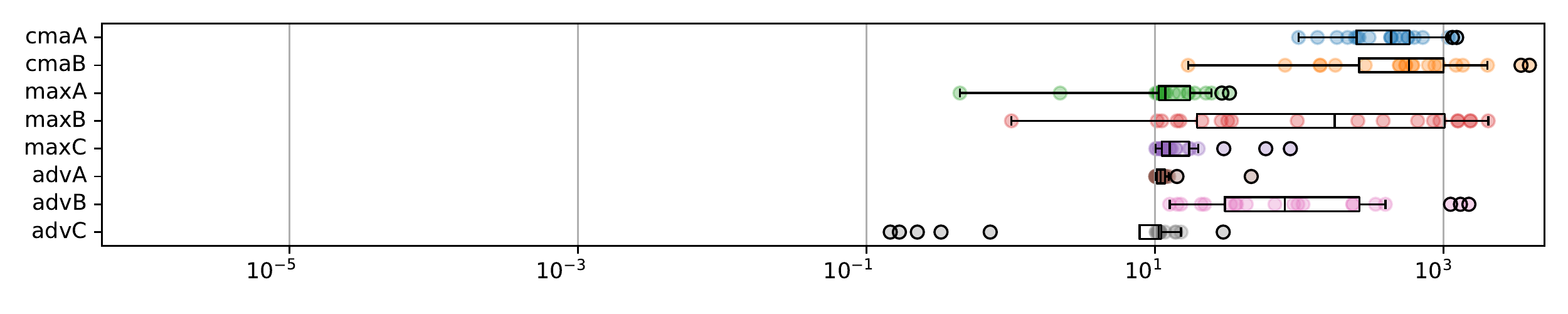}%
    \caption{$\max_{y \in \Y_{C}} f(x, y)$}\label{fig:berth:c}%
    \end{subfigure}%
    \caption{Performance of the controllers obtained in 20 independent trials of (1+1)-CMA-ES on $y = (y_\text{no}^{(A)}, y_\text{est}^{(B)})$ and $y = (y_\text{est}^{(A)}, y_\text{est}^{(B)})$; CMA-ES($N_y = 100$) on $\Y_{A}$, $\Y_{B}$, $\Y_{C}$; and Adversarial-CMA-ES on $\Y_{A}$, $\Y_{B}$, $\Y_{C}$, denoted by \texttt{cmaA}, \texttt{cmaB}, \texttt{maxA}, \texttt{maxB}, \texttt{maxC}, \texttt{advA}, \texttt{advB}, and \texttt{advC}, respectively. Each box indicates the lower quartile $Q1$ and the upper quartile $Q3$, with the line indicating the median $Q2$. The lower and upper whiskers are the lowest datum above $Q1 - 1.5(Q3-Q1)$ and the highest datum below $Q3 + 1.5(Q3-Q1)$.}
    \label{fig:berthing}
\end{figure}

\Cref{fig:berthing} shows the performance of the resulting controllers of 20 independent trials of each algorithm under different situations. 
Some of the trajectories observed for the obtained controllers are discussed in \Cref{apdx:berth}.

The (1+1)-CMA-ES on $y = (y_\text{no}^{(A)}, y_\text{est}^{(B)})$ achieves the best performance under no wind disturbance (\Cref{fig:berth:a0}), while
 (1+1)-CMA-ES on $y = (y_\text{est}^{(A)}, y_\text{est}^{(B)})$ achieved the best performance under the certain wind condition, $y^{(A)} = y_\mathrm{est}^{(A)}$ (\Cref{fig:berth:b0}).
In all trials, they achieved a cost of $< 10^{-3}$. 
However, their performances significantly were degraded under the worst case, particularly when the wind condition was unknown (\Cref{fig:berth:a,fig:berth:c}), where the ship collided with the berth and the cost was $> 10$. 
The performance was less affected by the uncertainty in the model coefficients in this experiment, but the effect would be enhanced if we consider a wider uncertainty set $\Y_{(B)}$. Nonetheless, these results demonstrate the importance of considering model uncertainty to obtain robust berthing control.

The controllers obtained by Adversarial-CMA-ES and CMA-ES($N_y = 100$) on $\Y_{A}$ and $\Y_{C}$ achieved better performance under the worst situation in $\Y_{A}$ (\Cref{fig:berth:a}) than those obtained by the other approaches. 
The numbers of trials that \texttt{maxA}, \texttt{maxC}, \texttt{advA}, and \texttt{advC} succeed in avoiding collision with the berth under the worst case in $\Y_{A}$ were $3$, $1$, $4$, and $11$ out of $20$. Interestingly, Adversarial-CMA-ES on $\Y_{C}$ (\texttt{advC}) achieved better worst-case performance in $\Y_{A}$ than Adversarial-CMA-ES on $\Y_{A}$ (\texttt{advA}), even though \texttt{advC} considered a wider uncertainty set and hence was expected to show more conservative performance. The reason for the superior performance of \texttt{advC} may be explained as follows. Both $y^{(A)}$ and $y^{(B)}$ represent the uncertainty regarding the wind force. Considering the worst case in $\Y_{C}$ results in being more conservative to the wind force. This may help to obtain a controller that can avoid the collision with the berth, while possibly losing the accuracy of the final position. 

The advantage of Adversarial-CMA-ES over CMA-ES($N_y = 100$) was more pronounced in the worst-case performance on $\Y_{B}$ (\Cref{fig:berth:b}). 
The median of \texttt{advB} and that of \texttt{maxB} were better than the median of the other results. 
All 20 trials of \texttt{advB} achieved berthing without collision with the berth in the worst situation. 
In contrast, 9 out of 20 trials failed in \texttt{maxB}. This may be because $N_y = 100$ was not sufficiently large to represent the uncertainty in the 10-dimensional space $\Y_{B}$. 

In the worst-case performance on $\Y_{C}$, 
5 out of 20 trials of \texttt{advC} succeed in avoiding a collision with the berth, whereas all the trials of \texttt{maxC} failed in avoiding a collision with the berth.
Because of the similarity between the results on the worst-case performances on $\Y_{A}$ and $\Y_{C}$, the difficulties in obtaining a robust controller under the worst-case in $\Y_{C}$ mainly comes from the difficulty in treating the worst case in $\Y_{A}$. 
The results may be improved by running the optimization process longer and performing more restarts to locate more local min--max saddle points.


\section{Conclusion}

We have proposed a framework for saddle point optimization with approximate minimization oracles.
Our theoretical analysis has shown the conditions under which the learning rate for the approach converges linearly (i.e., geometrically) toward the min--max saddle point on strongly convex--concave functions. 
Numerical evaluation have shown the tightness of the theoretical results. 
We have also proposed a learning rate adaptation mechanism for practical use. 
Numerical evaluation on convex-concave quadratic problems has demonstrated that the proposed approach with the learning rate adaptation successfully converges linearly toward the min--max saddle point, with the compromise of $f$-calls being no more than three times that of $f$-calls with the best tuned fixed learning rate. 
Comparison with other baseline approaches on several test problems revealed the limitations of existing coevolutionary approaches as well as of the proposed approach on problems with the optimal solution that is not a min--max saddle point. 
The application of the proposed approach to a robust berthing control task demonstrated the usefulness of the proposed approach, and the results imply the importance of considering modeling errors to achieve a reliable and safe solution. 

We conclude the present work with some suggestions for possible avenues for future research.

The main limitation of the proposed approach as a solver to \eqref{eq:spo} is that it may fail to converge to a local minimal solution of the worst-case objective $\max_{y \in \Y} f(x, y)$ if it is not a strict local min--max saddle point of $f$.
Such failure cases were observed in \Cref{fig:ex2}, not only for the proposed approach but also for existing coevolutionary approaches. 
Addressing this difficulty is an important direction for future work. 
For the GDA approach \eqref{eq:grad}, \citet{liang2019interaction} have shown that the GDA failed to converge to the optimal solution on a bi-linear function $f(x, y) = x^\T C y$ and some improved gradient-based approaches \cite{daskalakis2017training,yadav2017stabilizing,Mescheder2017nips} successfully converged. 
We expect that these gradient-based approaches would be useful in improving the here-proposed approach.
The other limitation is that the best possible runtime $\Omega(1/\gamma^*)$ in \eqref{eq:crstar} scales as the interaction term; more precisely, $\beta_G / \alpha_H$, increases. We intend to address this limitation in future work.

The main limitation of our theoretical result (\Cref{thm:conv}) is that approximate minimization oracles are required to satisfy \Cref{def:amo}.
In practice, it is often impossible to guarantee \Cref{def:amo} as we do not know the global minimum/maximum of the objective functions. For the design of Adversarial-CMA-ES and Adversarial-SLSQP, we expect that it is approximately satisfied by running a linearly convergence approximate minimization oracle until a fixed number of improvements are observed. See \Cref{sec:acmaes} for details. In such cases, we have condition \eqref{eq:amo-cond} not with a constant $\epsilon$ but rather with a possibly stochastic and time-dependent sequence $\{\epsilon_t\}$, which is not covered by \Cref{thm:conv}. Dealing with such $\epsilon_t$ will enlarge the scope of \Cref{thm:conv} and bridge the gap between theory and practice.\footnote{Another possible approach is to replace condition $h(\tilde{z}) - \min_{z \in \Z} h(z) \leq \epsilon \cdot (h(\bar{z}) - \min_{z \in \Z} h(z))$ in \Cref{def:amo} with condition $\norm{\nabla h(\tilde{z})} < \epsilon \cdot \norm{\nabla h(\bar{z})}$ under the additional assumption that the objective function is continuously differentiable. An advantage of this approach is that this condition can be approximately verified by estimating the gradients of the objective function \cite{Nesterov2017}.}

The uncertainty parameter is typically constrained in a bounded set $\Y \subset \R^n$ in practice, however, the effect of $\Y$ on the convergence rate has not been investigated in this work. Our theoretical result was developed for unconstrained situation and the proof does not immediately generalize to the constrained situation. The effect of the bound $\Y$ on the convergence rate is to be investigated theoretically and empirically.

The experimental results on the robust berthing control task have demonstrated the usefulness of the proposed approach and the importance of considering model uncertainties. Simultaneously, they revealed the difficulty of obtaining a robust solution with satisfactory utility. Regarding the wind condition uncertainty, it is possible to decompose $\Y_{A}$ into disjoint subsets (e.g., based on the wind direction), train the robust feedback controller for each subset, and switch the controller based on the wind condition measured at the time of operation. Such an approach is not available for the uncertainty in the model coefficients. To improve the worst-case performance, it is important to reduce the set of uncertain parameter values $\Y$ as much as possible. In our experiments, we defined the interval for each uncertain coefficient to form $\Y$, but the corner case may be unrealistic and will degrade the worst-case performance unnecessarily. Designing more intelligent $\Y$ remains as a very important task for practical applications.

\begin{acks}
The authors would like to thank the anonymous reviewers for their valuable comments and suggestions. 
This work is partially supported by JSPS KAKENHI Grant Number 19H04179 and 19K04858.
\end{acks}


\appendix
\section{Proofs}\label{apdx:proof}

\subsection{Proof of \Cref{prop:saddle-error}}
\begin{proof}
  Assume that $(x^*, y^*)$ is a local min--max saddle point of $f$. Then, by definition, there exists a neighborhood $\mathcal{E}_x\times\mathcal{E}_y$ of $(x^*, y^*)$ such that $f(x, y^*) \geq f(x^*, y^*) \geq f(x^*, y)$ holds for any $(x, y) \in \mathcal{E}_x\times\mathcal{E}_y$. Let $(x, y) \in \mathcal{E}_x\times\mathcal{E}_y$. Then, $G_x(x, y^*) = f(x, y^*) - \min_{x' \in \X} f(x', y^*) \geq f(x^*, y^*) - \min_{x' \in \X} f(x', y^*) = G_x(x^*, y^*)$ and $G_y(x^*, y) = \max_{y' \in \Y} f(x^*, y') - f(x^*, y) \geq \max_{y' \in \Y} f(x^*, y') - f(x^*, y^*) = G_y(x^*, y^*)$. This implies that $x^*$ and $y^*$ are local minimal points of $G_x(x, y^*)$ and $G_y(x^*, y)$, respectively.

  Conversely, assume that $x^*$ and $y^*$ are local minimal points of $G_x(x, y^*)$ and $G_y(x^*, y)$, respectively. Then, there exists a neighborhood $\mathcal{E}_x\times\mathcal{E}_y$ of $(x^*, y^*)$ such that $G_x(x, y^*) \geq G_x(x^*, y^*)$ and $G_y(x^*, y) \geq G_y(x^*, y^*)$ for any $(x, y) \in \mathcal{E}_x\times\mathcal{E}_y$. They read $f(x, y^*) \geq f(x^*, y^*)$ and $f(x^*, y^*) \geq f(x^*, y)$, which implies that $(x^*, y^*)$ is a local min--max saddle point of $f$.

  If $(x^*, y^*)$ is the global min--max saddle point of $f$, then $(x^*, y^*)$ is a local minimal point of $G$. Moreover, we have $G_x(x^*, y^*) = G_y(x^*, y^*) = 0$, implying that it is the global minimal point of $G$.
  Conversely, if $(x^*, y^*)$ is the global minimal point of $G$, then it is a local min--max saddle point.
  Moreover, because the global minimum of $G$ is zero, we have $G_x(x^*, y^*) = G_y(x^*, y^*) = 0$.
  Then, we can take $\mathcal{E}_x = \X$ and $\mathcal{E}_y = \Y$ in the above proof, which implies that $(x^*, y^*)$ is the global min--max saddle point.
\end{proof}

\subsection{Proof of \Cref{lem:hess-gap}}
\begin{proof}
Noting that $(\nabla_x f)(\optx(y), y) = 0$ and $(\nabla_y f)(x, \opty(x)) = 0$, we obtain
\begin{equation}
\begin{split}
  \nabla_x G_x(x, y) &= (\nabla_x f)(x, y) , \qquad
  \nabla_y G_x(x, y) = (\nabla_y f)(x, y) - (\nabla_y f)(\optx(y), y ) ,\\
  \nabla_x G_y(x, y) &= (\nabla_x f)(x, \opty(x) ) - (\nabla_x f)(x, y) , \qquad
  \nabla_y G_y(x, y) = - (\nabla_y f)(x, y) .
\end{split}
\label{eq:nabla_g}
\end{equation}
Moreover, we have
\begin{align*}
  \nabla^2 G_x(x, y)
  &= \begin{bmatrix}
    H_{x,x}(x, y) & H_{x,y}(x, y) \\
    H_{y,x}(x, y) & H_{y,y}(x, y) - H_{y,y}(\optx(y), y ) - [J_{\optx}(\optx(y), y)]^\T H_{x,y}(\optx(y), y) 
  \end{bmatrix},
  \\
  \nabla^2 G_y(x, y) 
  &= \begin{bmatrix}
    - H_{x,x}(x, y) + H_{x,x}(x, \opty(x) ) + [J_{\opty}(x, \opty(x))]^\T H_{y,x}(x, \opty(x))  & -H_{x,y}(x, y) \\
    - H_{y,x}(x, y) & - H_{y,y}(x, y)
  \end{bmatrix} .
\end{align*}
In light of \Cref{prop:implicit} and the symmetry $H_{x,y} = H_{y,x}^\T$, we have $[J_{\optx}(\optx(y), y)]^\T = - H_{y,x}(\optx(y), y)(H_{x,x}(\optx(y), y))^{-1}$ and $[J_{\opty}(x, \opty(x))]^\T = - H_{x,y}(x, \opty(x))(H_{y,y}(x, \opty(x)))^{-1}$. 
Then, because $\nabla^2 G = \nabla^2(G_x + G_y) = \nabla^2 G_x + \nabla^2 G_y$, we obtain
\begin{align*}
  \nabla^2 G(x, y)
  &= \begin{bmatrix}
    G_{x,x}(x, \opty(x)) & 0 \\
    0 & G_{y,y}(\optx(y), y )
  \end{bmatrix}.
\end{align*}
The symmetry of $G_{x,x}$ and $G_{y,y}$ are clear from their definitions.
The positivity of $G_{x,x}$ and $G_{y,y}$ follows that $H_{x,x} \succ 0$, $-H_{y,y} \succ 0$, $H_{x,y}(- H_{y,y})^{-1} H_{y,x} \succcurlyeq \sigma_{\min}(H_{x,y})^2 / \sigma_{\max}(- H_{y,y}) \succ 0$ and $H_{y,x}(H_{x,x})^{-1} H_{x,y} \succcurlyeq \sigma_{\min}(H_{x,y})^2 / \sigma_{\max}(H_{x,x}) \succ 0$.
This completes the proof.
\end{proof}

\subsection{Proof of \Cref{thm:conv}}
\begin{proof}
  Let $v_x = \tilde{x}_t - x_t$, $v_y = \tilde{y}_t - y_t$, and $v = (v_x, v_y)$.
  Define $\bx(\tau) = x_t + \tau \cdot v_x$ and $\by(\tau) = y_t + \tau \cdot v_y$.
  Let $w_x = \optx(y_t) - x_t$ and $w_y = \opty(x_t) - y_t$.
  Define $\bbx(\tau) = x_t + \tau \cdot w_x$ and $\bby(\tau) = y_t + \tau \cdot w_y$.
  Then, $\bx(0) = \bbx(0) = x_t$ and $\by(0) = \bby(0) = y_t$.
  Moreover, $\bx(\eta) = x_{t+1}$ and $\by(\eta) = y_{t+1}$. 
  
  We first derive several inequalities on the norms of $v$ and $w$. 
  Note that $\nabla_x f(\optx(y_t), y_t) = 0$ and $\nabla_y f(x_t, \opty(x_t)) = 0$. 
  In light of conditions (1)--(2) in the theorem statement, we have
  \begin{align}
    (\alpha_H / 2) \norm{w_x}_{H_{x,x}^*}^2 \leq G_{x}(x_t, y_t) &= f(x_t,y_t) - f(\optx(y_t), y_t) \leq (\beta_H / 2) \norm{w_x}_{H_{x,x}^*}^2,\label{eq:11}
    \\
    (\alpha_H / 2) \norm{w_y}_{-H_{y,y}^*}^2 \leq G_{y}(x_t, y_t) &= f(x_t,\opty(x_t)) - f(x_t, y_t)\leq (\beta_H / 2) \norm{w_y}_{-H_{y,y}^*}^2 .\label{eq:12}
  \end{align}
  Because of condition~\eqref{eq:amo-cond}, we have
  \begin{align}
    (\alpha_H / 2) \norm{w_x - v_x}_{H_{x,x}^*}^2 &\leq f(\tilde{x}_t, y_t) - f(\optx(y_t), y_t) \leq \epsilon \cdot G_{x}(x_t, y_t),\label{eq:13}
    \\
    (\alpha_H / 2) \norm{w_y - v_y}_{-H_{y,y}^*}^2 &\leq f(x_t, \opty(x_t)) - f(x_t, \tilde{y}_t) \leq \epsilon \cdot G_y(x_t, y_t) .\label{eq:14}
  \end{align}
  Then, from \Cref{eq:11,eq:12,eq:13,eq:14}, we have
  \begin{align}
    \norm{v_x}_{H_{x,x}^*}^2 &\leq (\norm{w_x}_{H_{x,x}^*} + \norm{v_x - w_x}_{H_{x,x}^*} )^2 \leq 2 [ (1 + \sqrt{\epsilon})^2 / \alpha_H] \cdot G_{x}(x_t, y_t), \label{eq:15}
    \\
    \norm{v_y}_{H_{y,y}^*}^2 &\leq (\norm{w_y}_{-H_{y,y}^*} + \norm{v_y - w_y}_{-H_{y,y}^*} )^2 \leq 2 [(1 + \sqrt{\epsilon})^2 / \alpha_H] \cdot G_{y}(x_t, y_t). \label{eq:16}
  \end{align}
 Because $\epsilon < \alpha_H / \beta_H$, we also have
  \begin{align}
    \norm{v_x}_{H_{x,x}^*}^2 &\geq (\norm{w_x}_{H_{x,x}^*} - \norm{v_x - w_x}_{H_{x,x}^*} )^2 
    \geq (2 / \beta_H) [ 1 - \sqrt{(\beta_H / \alpha_H) \cdot \epsilon} ]^2 G_{x}(x_t, y_t), \label{eq:17}
    \\
    \norm{v_y}_{-H_{y,y}^*}^2 &\geq (\norm{w_y}_{-H_{y,y}^*} - \norm{v_y - w_y}_{-H_{y,y}^*} )^2 
    \geq (2 / \beta_H) [ 1 - \sqrt{(\beta_H / \alpha_H) \cdot \epsilon} ]^2 G_{y}(x_t, y_t). \label{eq:18}
  \end{align}

  By applying the mean value theorem repeatedly, we have
  \begin{equation}\begin{split}
    G(x_{t+1}, y_{t+1}) - G(x_{t}, y_{t})
    &= G(\bx(\eta), \by(\eta)) - G(\bx(0), \by(0))
    \\
    &= \int_{0}^{\eta} \nabla G(\bx(\tau), \by(\tau))^\T \mathrm{d}\tau \cdot v
    \\
    &= \int_{0}^{\eta} \left[\nabla G(\bx(0), \by(0)) + \int_{0}^{\tau} \nabla^2 G(\bx(s), \by(s)) \mathrm{d}s \cdot v\right]^\T \mathrm{d}\tau \cdot v
    \\
    &= \eta \cdot \nabla G(x_t, y_t)^\T \cdot v + v^\T \int_{0}^{\eta} \int_{0}^{\tau} \nabla^2 G(\bx(s), \by(s))^\T \mathrm{d}s \mathrm{d}\tau \cdot v
    .
    \end{split}\label{eq:1}
  \end{equation}
  First, we evaluate the first term on the right-most side of \eqref{eq:1}. Note first that we have $\nabla G(x_t, y_t)^\T v = \nabla_x G(x_t, y_t)^\T v_x + \nabla_y G(x_t, y_t)^\T v_y$. In light of \eqref{eq:nabla_g}, we have $\nabla_x G(x_t, y_t) = (\nabla_x f)(x_t, \opty(x_t)) = (\nabla_x f)(\bbx(0), \bby(1))$. Noting that $(\nabla_x f)(\bbx(1), \bby(0)) = (\nabla_x f)(\optx(y_t), y_t) = 0$, by the mean value theorem, we obtain
  \begin{equation}
    \nabla_x G(x_t, y_t)^\T v_x
    = -w_x^\T \cdot \int_{0}^{1} H_{x,x}(\bbx(1-\tau), \bby(\tau)) \mathrm{d}\tau \cdot v_x 
    + w_y^\T \cdot \int_{0}^{1} H_{y,x}(\bbx(1-\tau), \bby(\tau)) \mathrm{d}\tau \cdot v_x 
    .\label{eq:nabla_gxt}
  \end{equation}
  Analogously, we obtain
  \begin{equation}
    \nabla_y G(x_t, y_t)^\T v_y
    = w_y^\T \cdot \int_{0}^{1} H_{y,y}(\bbx(1 - \tau), \bby(\tau)) \mathrm{d}\tau \cdot v_y 
    - w_x^\T \cdot \int_{0}^{1} H_{x,y}(\bbx(1 - \tau), \bby(\tau)) \mathrm{d}\tau \cdot v_y 
    .\label{eq:nabla_gyt}
  \end{equation}
  Inserting $v_x = w_x + (v_x - w_x)$ and $v_x = w_y + (v_y - w_y)$ into \eqref{eq:nabla_gxt} and \eqref{eq:nabla_gyt}, using $H_{x,y} = H_{y,x}^{\T}$, and rearranging them, we obtain
  \begin{equation}
  \begin{split}
    \nabla G(x_t, y_t)^\T v
    &= \int_{0}^{1} \left( 
    - \begin{bmatrix}
    w_x \\
    w_y
    \end{bmatrix}^\T
    \begin{bmatrix}
     H_{x,x} & 0 \\
    0 & - H_{y,y}
    \end{bmatrix}
    \begin{bmatrix}
    w_x \\
    w_y
    \end{bmatrix}
    - \begin{bmatrix}
    v_x - w_x \\
    v_y - w_y
    \end{bmatrix}^\T
    \begin{bmatrix}
     H_{x,x} & - H_{x,y} \\
     H_{y,x} & - H_{y,y}
    \end{bmatrix}
    \begin{bmatrix}
    w_x \\
    w_y
    \end{bmatrix}
    \right)\mathrm{d}\tau,
  \end{split}
  \label{eq:nabla_gv}
  \end{equation}
where we drop $(\bbx(1-\tau), \bby(\tau))$ from $H_{x,x}$, $H_{x,y}$, $H_{y,x}$ and $H_{y,y}$ for compact expressions.

We aim to bound each term of \eqref{eq:1} and \eqref{eq:nabla_gv}. 
The second term on \eqref{eq:1} is bounded by using conditions (3) and (4) of the theorem statement as well as the fact that $\nabla^2 G$ is block-diagonal (\Cref{lem:hess-gap}) as
\begin{equation}
    \alpha_{G} \norm{v}_{\diag(H_{x,x}^*, -H_{y,y}^*)}^2  \leq v^\T \nabla^2 G(\bx(s), \by(s))^\T \cdot v \leq \beta_{G} \norm{v}_{\diag(H_{x,x}^*, -H_{y,y}^*)}^2 . \label{eq:111}
\end{equation}
The first term on \eqref{eq:nabla_gv} is bounded by using conditions (1) and (2) of the theorem statement as
\begin{equation}
    -\beta_{H} \norm{w}_{\diag(H_{x,x}^*, -H_{y,y}^*)}^2 \leq
    - \begin{bmatrix}
    w_x \\
    w_y
    \end{bmatrix}^\T
    \begin{bmatrix}
     H_{x,x} & 0 \\
    0 & - H_{y,y}
    \end{bmatrix}
    \begin{bmatrix}
    w_x \\
    w_y
    \end{bmatrix}
    \leq - \alpha_{H} \norm{w}_{\diag(H_{x,x}^*, -H_{y,y}^*)}^2
    .\label{eq:112}
\end{equation}
The second term on \eqref{eq:nabla_gv} is bounded as
\begin{multline}
\abs*{\begin{bmatrix}
    v_x - w_x \\
    v_y - w_y
    \end{bmatrix}^\T
    \begin{bmatrix}
     H_{x,x} & - H_{x,y} \\
     H_{y,x} & - H_{y,y}
    \end{bmatrix}
    \begin{bmatrix}
    w_x \\
    w_y
    \end{bmatrix}}
    \leq     \norm{w - v}_{\diag(H_{x,x}^*, -H_{y,y}^*)}
    \cdot \norm{w}_{\diag(H_{x,x}^*, -H_{y,y}^*)}
    \\
    \cdot \sigma_{\max}\left(\begin{bmatrix}
     \sqrt{H_{x,x}^{*}}^{-1} H_{x,x} \sqrt{H_{x,x}^{*}}^{-1} & -\sqrt{H_{x,x}^{*}}^{-1} H_{x,y} \sqrt{-H_{y,y}^{*}}^{-1} \\
     \sqrt{-H_{y,y}^{*}}^{-1} H_{y,x} \sqrt{H_{x,x}^{*}}^{-1} & \sqrt{-H_{y,y}^{*}}^{-1} (-H_{y,y})
\sqrt{-H_{y,y}^{*}}^{-1}
\end{bmatrix}\right) ,\label{eq:113}
\end{multline}%
where the greatest singular value is bounded by using conditions (1)--(4) in the theorem statement as
\begin{equation}
    \sigma_{\max}\left(\begin{bmatrix}
     \sqrt{H_{x,x}^{*}}^{-1} H_{x,x} \sqrt{H_{x,x}^{*}}^{-1} & - \sqrt{H_{x,x}^{*}}^{-1} H_{x,y} \sqrt{-H_{y,y}^{*}}^{-1} \\
    \sqrt{-H_{y,y}^{*}}^{-1} H_{y,x} \sqrt{H_{x,x}^{*}}^{-1} & \sqrt{-H_{y,y}^{*}}^{-1} (-H_{y,y})
\sqrt{-H_{y,y}^{*}}^{-1}
\end{bmatrix}\right)
\leq 
\beta_H \cdot \beta_G^{1/2} / \alpha_{H}^{1/2} . \label{eq:nabla_gv2}
\end{equation}

\Cref{eq:1,eq:nabla_gv,eq:nabla_gv2,eq:111,eq:112,eq:113,eq:11,eq:12,eq:13,eq:14,eq:15,eq:16} lead to 
\begin{equation}
    G(x_{t+1}, y_{t+1}) - G(x_{t}, y_{t})
    \leq 
    \left[- 2 \eta \frac{\alpha_{H}}{\beta_{H}} \left( 1 - \frac{\beta_H^2}{\alpha_H^2} \sqrt{\frac{\beta_G}{\alpha_{H}}\cdot \epsilon} \right)
    + \eta^2 \cdot (1 + \sqrt{\epsilon})^2 \cdot \frac{\beta_{G}}{\alpha_H}
    \right] G(x_t, y_t).
\end{equation}
Here, the right-hand side is $\gamma \cdot G(x_{t}, y_{t})$ with $\gamma$ defined in \eqref{eq:cr}. Hence, $G(x_{t+1},y_{t+1}) \leq (1+\gamma) \cdot G(x_t, y_t)$. 
Note that $\log(1 + \gamma) < \gamma$ for all $\gamma \in (-1, 0)$, we thus obtain $\log\left( G(x_{t+1}, y_{t+1}) \right) - \log\left( G(x_t, y_t) \right) < \gamma$. Because $\log\left( G(x_{t}, y_{t}) \right) - \log\left( G(x_0, y_0) \right) < \gamma \cdot t$,
the minimal $t$ that $\log\left( G(x_{t}, y_{t}) \right) - \log\left( G(x_0, y_0) \right) \leq \log(\zeta)$ is no greater than $\left\lceil\frac{1}{\gamma}\log\left(\frac{1}{\zeta}\right)\right\rceil = T_\zeta$. 
Similarly, \Cref{eq:1,eq:nabla_gv,eq:nabla_gv2,eq:111,eq:112,eq:113,eq:11,eq:12,eq:13,eq:14,eq:17,eq:18} lead to
\begin{equation}
    G(x_{t+1}, y_{t+1}) - G(x_t, y_t) 
    \geq 
    \left[ - 2 \eta \frac{\beta_{H}}{\alpha_{H}}\left( 1 + \sqrt{\frac{\beta_G}{\alpha_{H}}\cdot \epsilon} \right)  
    + \eta^2 \frac{\alpha_{G}}{\beta_H} \left( 1 - \sqrt{\frac{\beta_H}{\alpha_{H}}\cdot \epsilon} \right)^2
    \right] G(x_t, y_t) \enspace.
    \label{eq:9}
\end{equation}
The right-hand side of \Cref{eq:9} is positive if $\eta > 2 \cdot \bar{\eta}$. This completes the proof.
\end{proof}

\section{Details of Automatic Berthing Control Problem}\label{apdx:berth:problem}

\paragraph{Subject Ship}
The control target is a 3 m model ship of MV ESSO OSAKA (\Cref{fig:esso_osaka}), following a related study \cite{maki2020,Maki2020b}.
The state variables $s = (X, u, Y, v_m, \psi, r) \in \R^6$ are the $X$ [m] and $Y$ [m] coordinates of the Earth-fixed coordinate system, the longitudinal velocity $u$ [m/s] and the lateral velocity $v_m$ [m/s] at the mid-ship, and the yaw direction $\psi$ [rad] as seen from the $X$ coordinates and the yaw angular velocity $r$ [rad/s].
The control signal $a = (\delta, n_\mathrm{p}, n_\mathrm{BT}, n_\mathrm{ST}) \in \R^4$ consists of the rudder angle $\delta$ [rad], propeller revolution number $n_{\mathrm{p}}$ [rps], the bow thruster revolution number $n_{\mathrm{BT}}$ [rps], and the stan-thruster revolution number $ n_{\mathrm{ST}}$ [rps].
Their feasible values are in $U = \big[-\frac{35}{180}\pi, \frac{35}{180}\pi\big] \times [-20, 20] \times [-20, 20] \times [-20, 20]$.
We employ the state equation model $\dot s = \phi(s, a; y)$ proposed in \cite{miyauchi2021}, where $y \in \Y$ represents the model uncertainty described in \Cref{sec:berth}.

\paragraph{Feedback Controller}
The feedback controller $u_x: \R^6 \to U$ is modeled by the following neural network parameterized by $x = (B, W, V)$: 
\begin {equation}
  u_x (s) = V \cdot \texttt {softmax} (\alpha \cdot (B + W \cdot s)) , \label {eq:nn}
\end {equation}
where $W \in [-1, 1]^{K \times 6} $ and $ B \in [-1, 1]^{K} $ define a linear map $z = \alpha \cdot (B + W \cdot s)$ from the state vector $s$ to the $K$ dimensional latent space, and $V \in U^{K} \subset \mathbb {R}^{m \times K}$ is a matrix consisting of $K$ feasible control vectors as its columns.
The softmax function
\begin{equation*}
    \texttt {softmax}: z = (z_1, \dots, z_K) \mapsto \frac{(\exp (z_1), \dots, \exp (z_K))}{\sum_ {k = 1} ^ {K} \exp (\exp (z_1), \dots, \exp (z_K))} \in \Delta^{K-1}
\end{equation*}
outputs a point in the $K-1$ dimensional standard simplex $ \Delta ^ {K-1} = \{z \in \mathbb {R} ^ K: z_1 \geq 0, \dots, z_K \geq 0, \ \text {and} \ \sum_ {k = 1} ^ {K} z_k = 1 \}$.
The output is a combination of the columns of $V$ weighted by the softmax output. 
$ \alpha> 0 $ is a parameter that determines whether the output of \texttt {softmax} is close to the one-hot vector.

The architecture of this neural network is interpreted as follows.
First, $ z = \texttt {softmax} (\alpha \cdot (B + W \cdot s)) $ on the first layer divides the state space into $K$ regions.
For example, if the greatest element of the vector $B + W \cdot s$ is the $k$th coordinate, then $z$ is approximated by the one-hot vector $e_k$ with $1$ on the $k$th coordinate and $0$ on the other coordinates if $\alpha$ is sufficiently large. 
In such a situation, $ u (s) = V \cdot z \approx V \cdot e_k = v_k$, where $v_k$ is the $k$th column of $V$.
In other words, this neural network approximates the control law that divides the state space using a Voronoi diagram with respect to the Euclidean metric and outputs the corresponding column of $V$ as a control signal in each region.
If we set $\alpha$ to be greater, $z$ is more likely to be close to a one-hot vector, which makes it easier to express the bang-bang type control.
If we set $\alpha$ to be smaller, $z$ is more likely to take a value in the middle of $\Delta^{K-1}$, which makes it easier to express a continuous control. 

Based on our preliminary experiments, we set $\alpha = 4$ and $K = 9$ in the following experiments.
Then, $x$ is of $m = 99$ dimension.

\paragraph{Objective Function}
The objective is to find the parameter $x := (B, W, V) \in \X$ of the controller $u_{x}$ that minimizes the cost $C$ of the trajectory $(s_{t \in [0, t_{\max}]}, a_{t \in [0, t_{\max}]})$ in the worst environment $y \in \Y$ for $u_{x}$. It is modeled as
\begin{align*}
  \min_{x \in \X} \max_{y \in \Y} f(x, y) &= \min_{x \in \X} \max_{y \in \Y} C( s_{t \in [0, t_{\max}]}, a_{t \in [0, t_{\max}]})\\
 \text{subject to}\  s_t &= s_0 + \int_{0}^{t} \phi(s_\tau, a_\tau; y) \mathrm{d}\tau \quad \text{and} \quad
  a_t = u_{x}(s_{\lfloor t / \mathrm{d}t \rfloor \cdot \mathrm{d}t}) ,
\end{align*}
where $\mathrm{d}t$ [seconds] is the control time span, that is, the control signal $a_t$ changes every $\mathrm{d}t$, and $s_0$ is the initial state. 

We define the cost of the trajectory as
\begin{align*}
  C( s_{t \in [0, t_{\max}]}, a_{t \in [0, t_{\max}]}) = C_1  + w \cdot (C_2 + \mathbb{I}\{C_2 > 0\}).
\end{align*}
where $w > 0$ is the hyperparameter that determines the trade-off between utility and safety, 
\begin{align*}
  C_1 = \frac16\sum_{i=1}^{6} (s_{t_{\max}, i} - s_{\mathrm{fin},i} )^2,
\end{align*}
evaluates the deviation of the final ship state from the target state $s_{\mathrm{fin}}$, and 
\begin{align*}
  C_2 = \frac14 \sum_{i=1}^{4} \int_{0}^{t_{\max}} \mathrm{dist}(P_{\tau,i}, C_\mathrm{berth}) \mathrm{d}\tau,
\end{align*}
measures the collision risk, where $P_{\tau,1}, \dots, P_{\tau,4}$ represents the coordinates of the four vertices of the rectangle surrounding the ship at time $\tau$ and $\mathrm{dist}(P, C_\mathrm{berth})$ measures the distance from a point $P$ to the closest point on the berth boundary. 
Refer to \cite{maki2020,Maki2020b} for the definitions of $C_1$ and $C_2$. 


Following \cite{maki2020}, we set $t_{\max} = 200$ [seconds] and $\mathrm{d}t = 10$ [seconds].
The initial state is $s_0 = \left( 15.0, 0.01,  6.0,  0.0, \pi, 0.0 \right)$ and
the target state is $s_\mathrm{fin} = \left( 3.0,  0.0,  9.5,  0.0, \pi, 0.0 \right)$.
The boundary of the berth is $C_\mathrm{berth} = \{Y = 9.994625\}$. 
The trade-off coefficient is set to $w = 10$. 
That is, the cost $f(x, y) < 10$ implies that the controller $u_x$ produces a trajectory without collision with the berth under the uncertainty parameter $y$. 

\paragraph{Differences from Previous Works}

Our problem formulation mostly follows previous studies \cite{maki2020,Maki2020b} but with certain differences. 
First, we optimize the feedback controller, whereas the control signals for each time period as well as the total control time are directly optimized in \cite{maki2020,Maki2020b}, which we believe is not suitable for obtaining robust control. 
Second, we modify the objective function. Previous studies include the term penalizing the control time as they formulate the problem as minimization of the control time. Because we did not optimize the control time, it is excluded from our objective function definition. Moreover, for better collision avoidance, we replaced $w \cdot C_2$ with $w \cdot (C_2 + \mathbb{I}\{C_2 > 0\})$. 
Third, following \cite{miyauchi2021}, we implement thrusters to realize robust control under external disturbances and adopt the state equation model proposed in \cite{miyauchi2021}.

\section{Additional Results for Automatic Berthing Control Problem}\label{apdx:berth}

\begin{figure}
    \centering
    \begin{subfigure}{0.5\hsize}%
    \centering
    \includegraphics[width=\hsize]{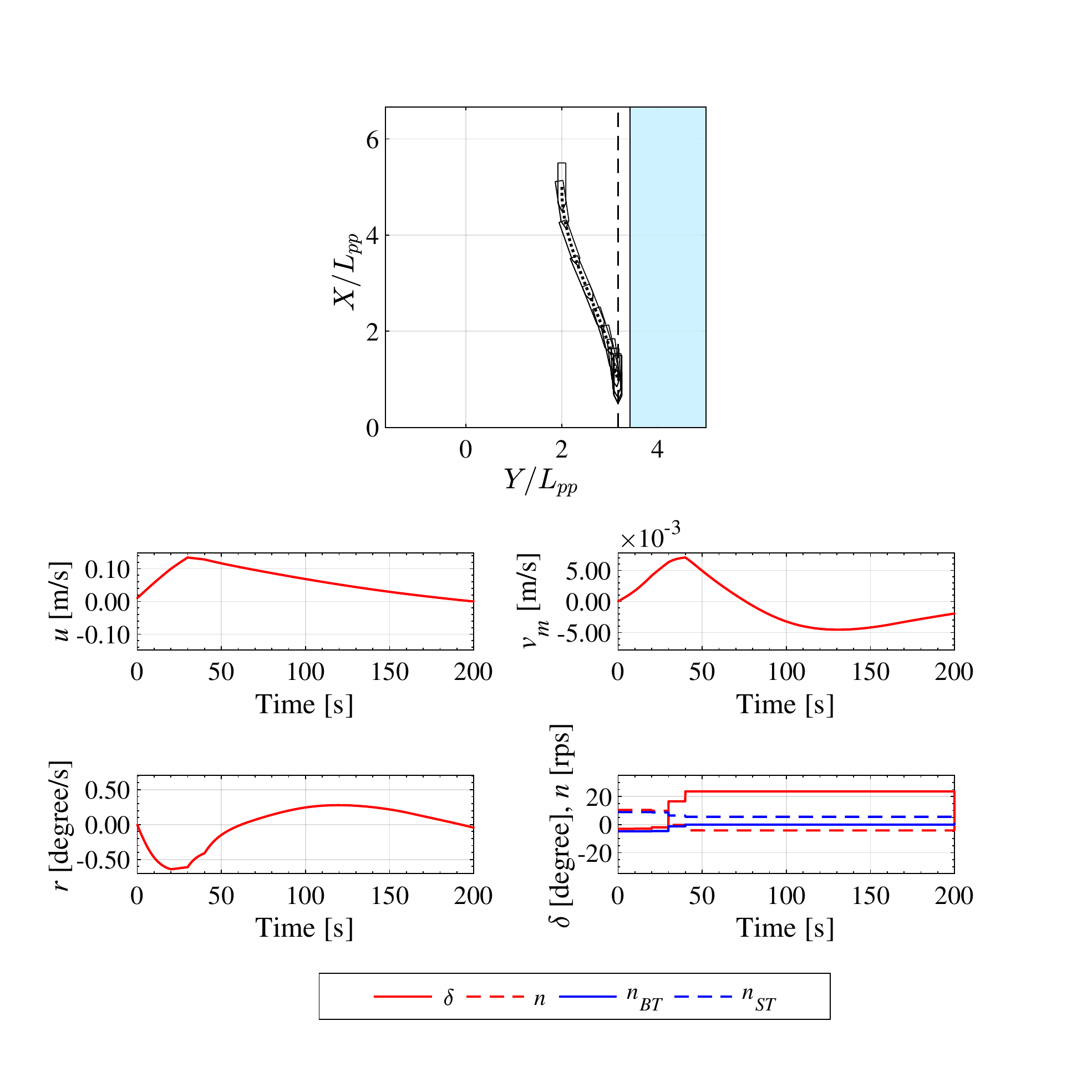}\\%
    $f(x, y) = 7.30 \times 10^{-7}$
    \caption{$y^{(A)} = y_\mathrm{no}^{(A)}$, $y^{(B)} = y_\mathrm{est}^{(B)}$}\label{fig:berth:cmaA:A0}%
    \end{subfigure}%
    \begin{subfigure}{0.5\hsize}%
    \centering
    \includegraphics[width=\hsize]{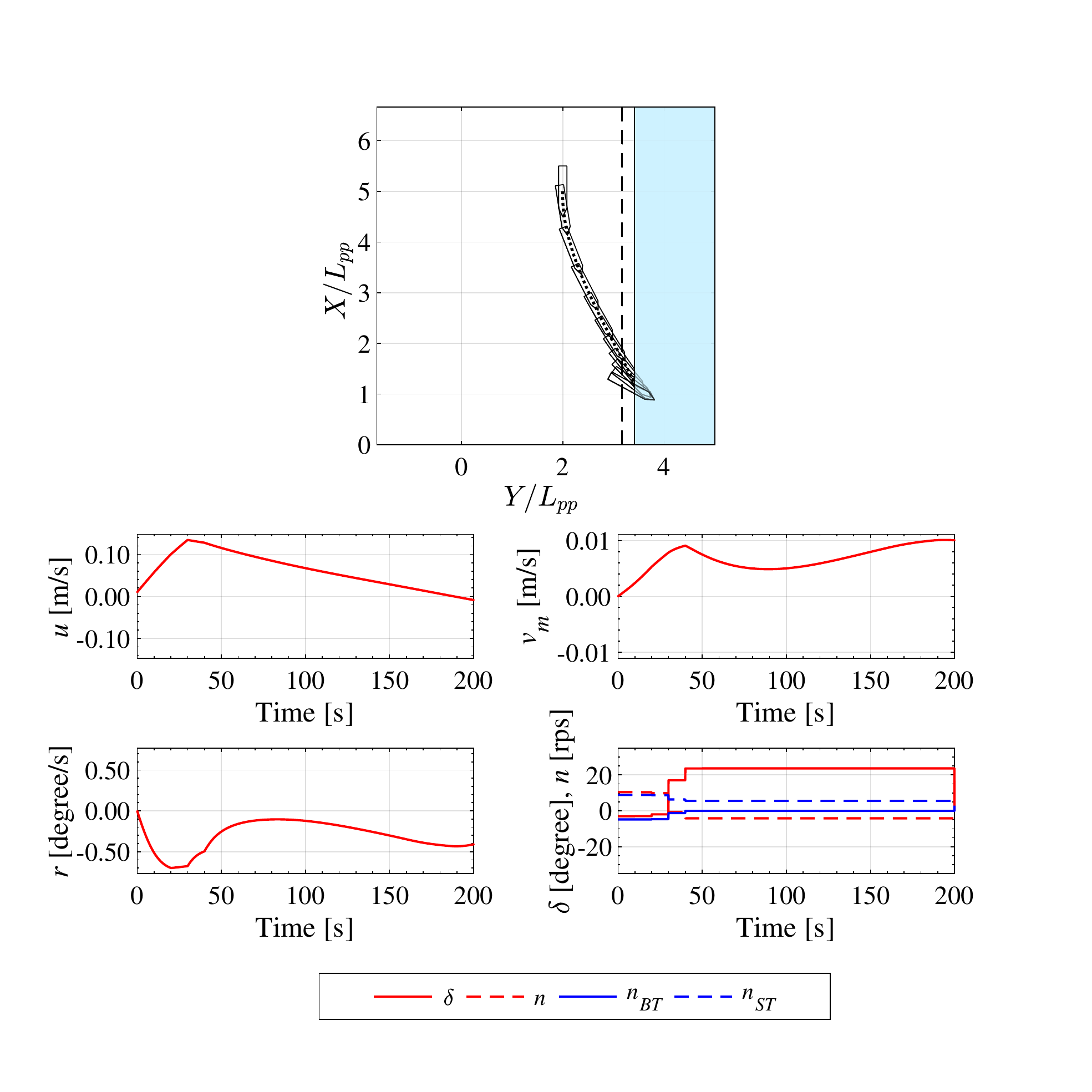}\\%
    $f(x, y) = 4.89 \times 10^{+2}$
    \caption{Worst case in $y^{(A)} \in \Y_{A}$, $y^{(B)} = y_\mathrm{est}^{(B)}$}\label{fig:berth:cmaA:A}%
    \end{subfigure}%
    \caption{Trajectories of the best controller obtained by CMA-ES($y_\mathrm{no}^{(A)}$)}
    \label{fig:berth:A:cma}
\end{figure}
\begin{figure}
    \begin{subfigure}{0.5\hsize}%
    \centering
    \includegraphics[width=\hsize]{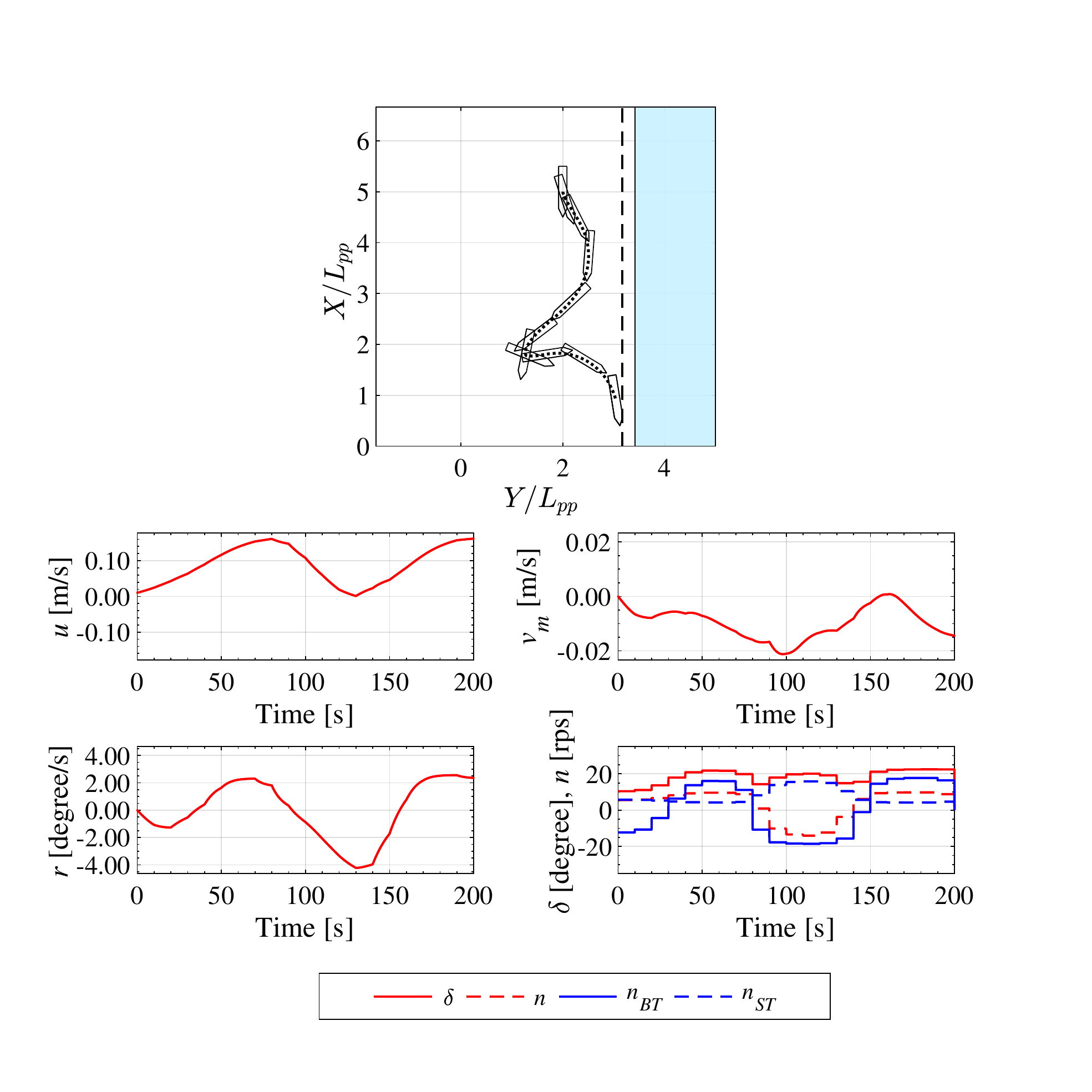}\\%
    $f(x, y) = 4.79 \times 10^{-2}$%
    \caption{$y^{(A)} = y_\mathrm{no}^{(A)}$, $y^{(B)} = y_\mathrm{est}^{(B)}$}\label{fig:berth:advA:A0}%
    \end{subfigure}%
    \begin{subfigure}{0.5\hsize}%
    \centering
    \includegraphics[width=\hsize]{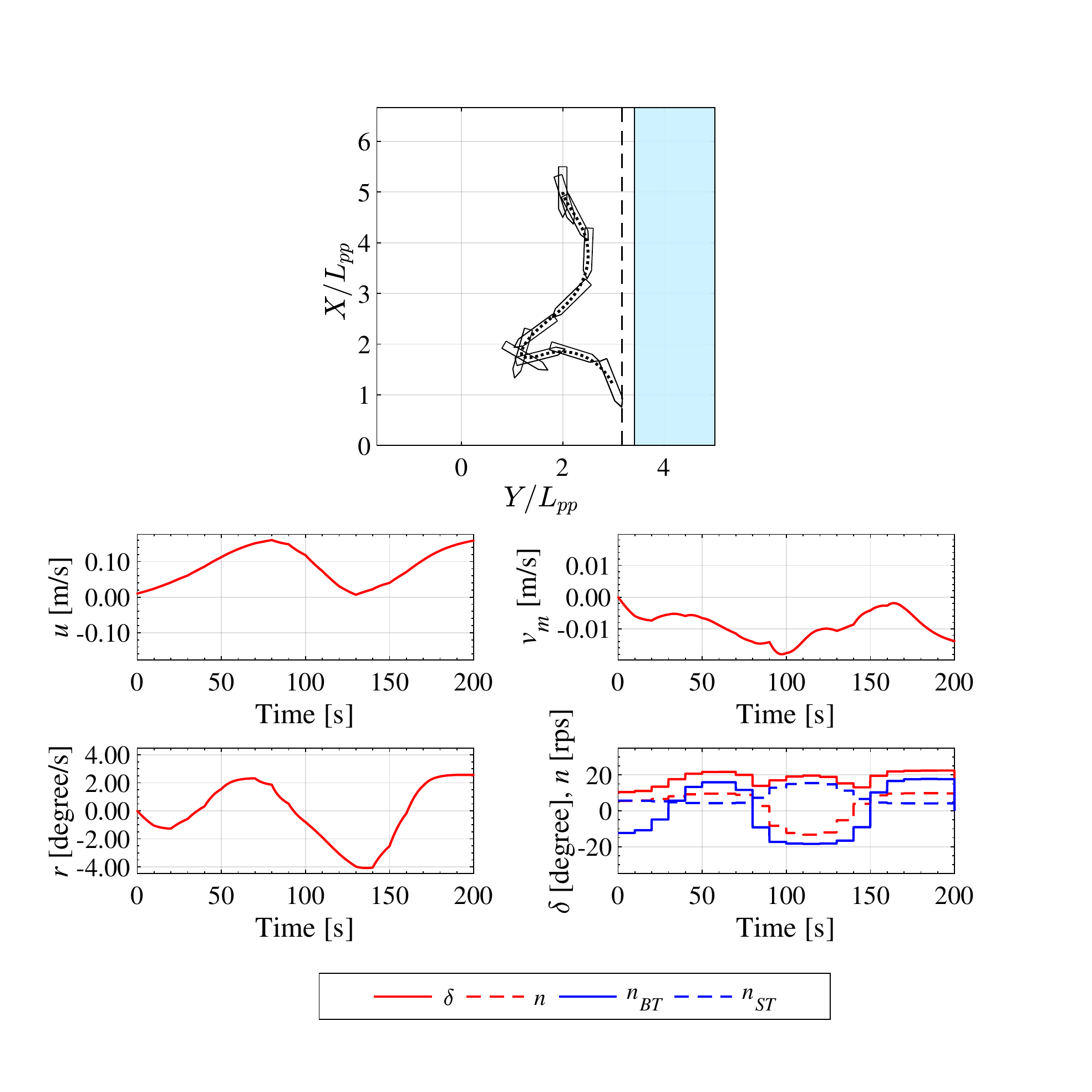}\\%
    $f(x, y) = 1.53 \times 10^{-1}$%
    \caption{Worst case in $y^{(A)} \in \Y_{A}$, $y^{(B)} = y_\mathrm{est}^{(B)}$}\label{fig:berth:advA:A}%
    \end{subfigure}%
    \caption{Trajectories of the best controller obtained by Adversarial-CMA-ES on $\Y_A$}
    \label{fig:berth:A:adv}
\end{figure}
\begin{figure}
    \centering
    \begin{subfigure}{0.5\hsize}%
    \centering
    \includegraphics[width=\hsize]{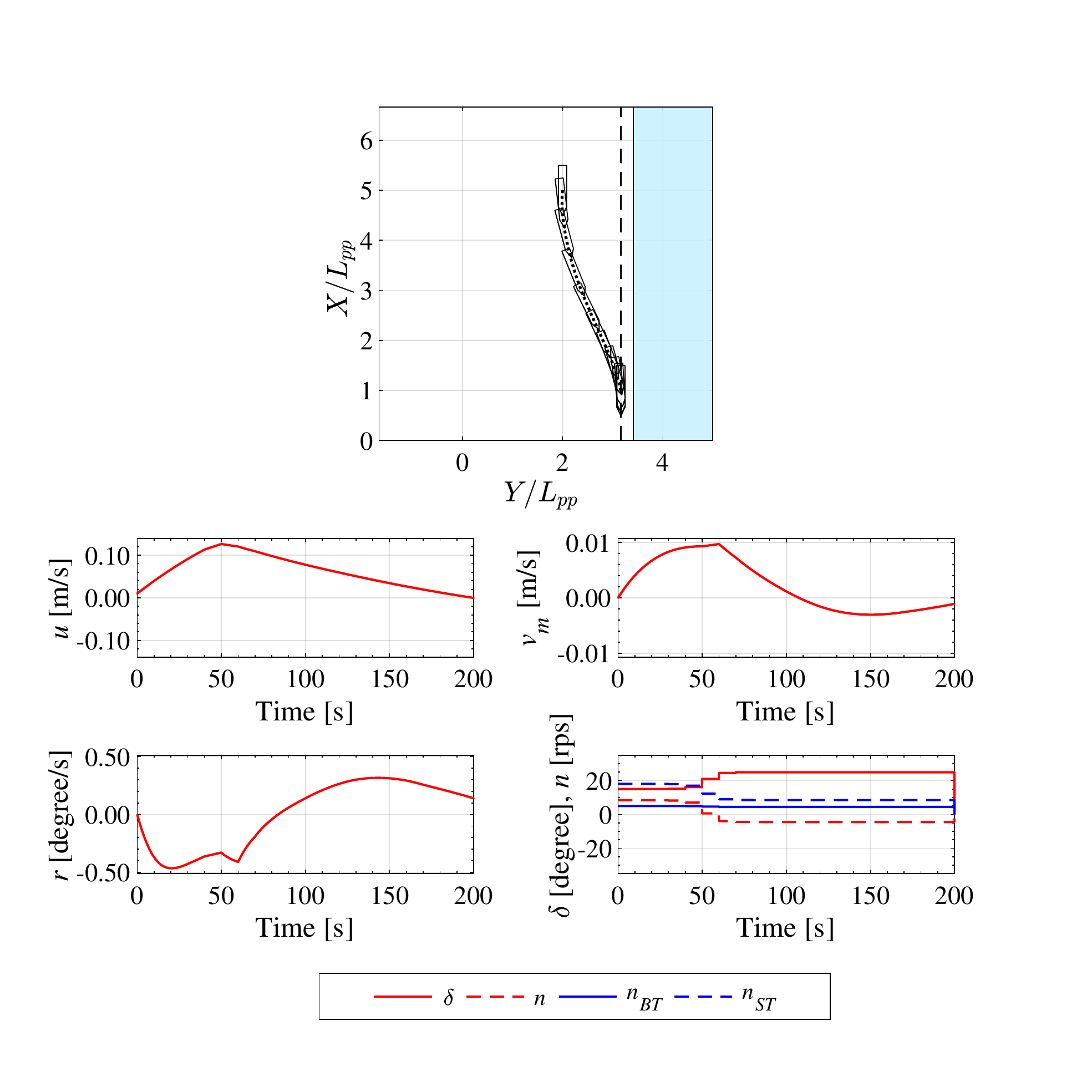}\\%
        $f(x, y) = 1.20 \times 10^{-6}$%
    \caption{$y^{(A)} = y_\mathrm{est}^{(A)}$, $y^{(B)} = y_\mathrm{est}^{(B)}$}\label{fig:berth:cmaB:B0}%
        \end{subfigure}%
    \begin{subfigure}{0.5\hsize}%
    \centering
    \includegraphics[width=\hsize]{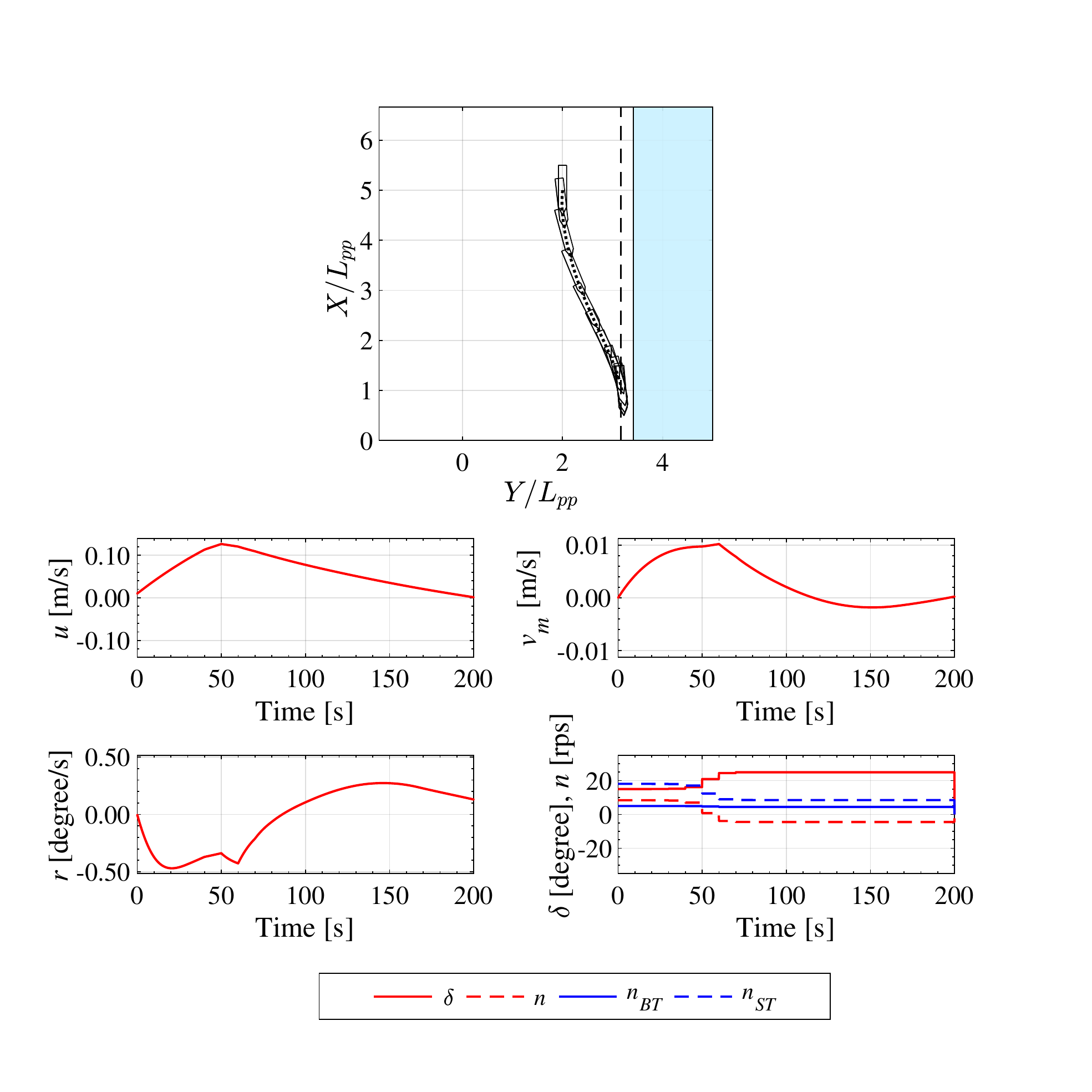}\\%
        $f(x, y) = 1.44 \times 10^{+1}$%
    \caption{Worst case in $y^{(B)} \in \Y_{B}$, $y^{(A)} = y_\mathrm{est}^{(A)}$}\label{fig:berth:cmaB:B}%
        \end{subfigure}%
    \caption{Trajectories of the best controller obtained by CMA-ES($y_\mathrm{est}^{(B)}$)}
    \label{fig:berth:B:cma}
\end{figure}
\begin{figure}
    \begin{subfigure}{0.5\hsize}%
    \centering
    \includegraphics[width=\hsize]{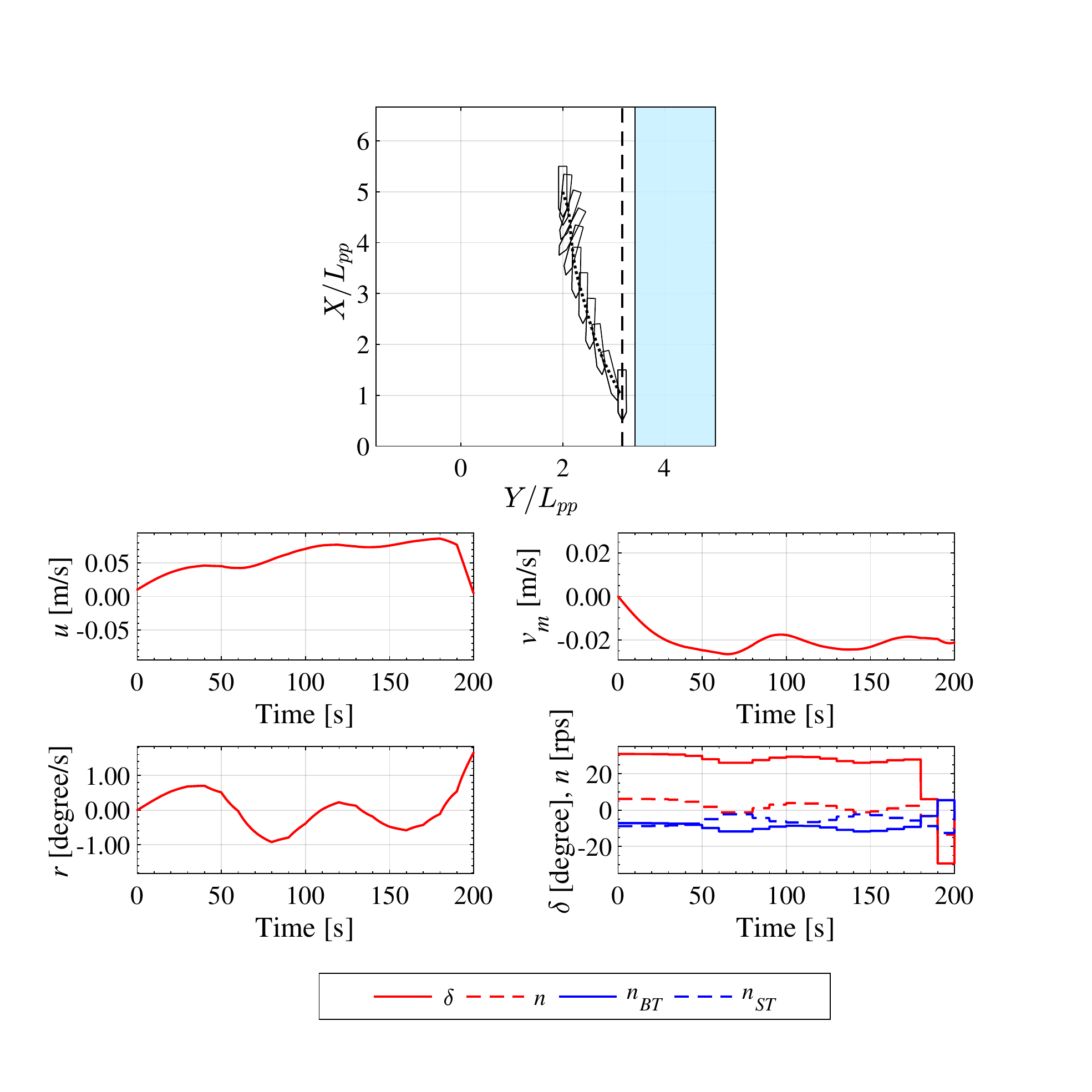}\\%
    $f(x, y) = 2.23 \times 10^{-4}$%
    \caption{$y^{(A)} = y_\mathrm{est}^{(A)}$, $y^{(B)} = y_\mathrm{est}^{(B)}$}\label{fig:berth:advB:B0}%
    \end{subfigure}%
    \begin{subfigure}{0.5\hsize}%
    \centering
    \includegraphics[width=\hsize]{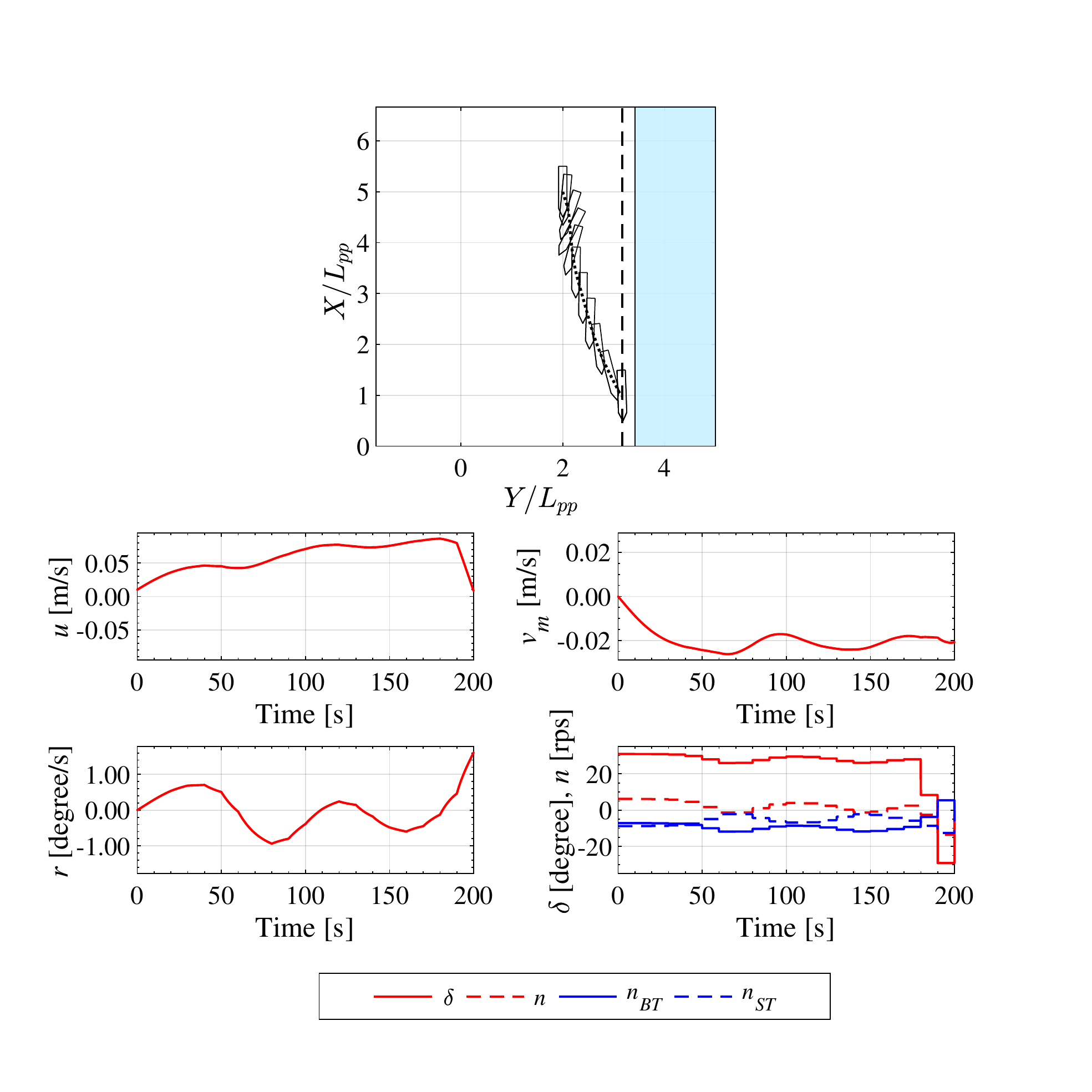}\\%
    $f(x, y) = 4.77 \times 10^{-4}$%
    \caption{Worst case in $y^{(B)} \in \Y_{B}$, $y^{(A)} = y_\mathrm{est}^{(A)}$}\label{fig:berth:advB:B}%
    \end{subfigure}%
    \caption{Trajectories of the best controller obtained by Adversarial-CMA-ES on $\Y_B$}
    \label{fig:berth:B:adv}
\end{figure}

\Cref{fig:berth:A:cma,fig:berth:A:adv,fig:berth:B:cma,fig:berth:B:adv} visualize the trajectories obtained in the experiments in \Cref{sec:berth}. The route of the ship, that is, $(X, Y, \psi)$ at each time, is displayed in the top figure. The $X$ and $Y$ axes are scaled by $L_{pp} = 3~[\mathrm{m}]$. The changes in the velocities, $(u, v_m, r)$, as well as the changes in the control signals, $(\delta, n_p, n_{BT}, n_{ST})$, are plotted at the bottom. Note that $r$ and $\delta$ are plotted on a degree basis for better intuition.
\Cref{fig:berth:A:cma} shows the trajectories observed for the best controller obtained by CMA-ES($y_\mathrm{no}^{(A)}$), which is the controller optimized under $y = (y_\mathrm{no}^{(A)}, y_\mathrm{est}^{(B)})$, that is, no wind $y^{(A)} = y_\mathrm{no}^{(A)}$ and model parameter $y^{(B)} = y_\mathrm{est}^{(B)}$ used in the previous study.
\Cref{fig:berth:A:adv} shows the trajectories observed for the best controller obtained by Adversarial-CMA-ES on $\Y_A$, which is the controller optimized under the worst wind condition $y^{(A)} \in \Y_{A}$ with $y^{(B)} = y_\mathrm{est}^{(B)}$. 
For \Cref{fig:berth:A:cma,fig:berth:A:adv}, the left figure is the trajectory under $y = (y_\mathrm{no}^{(A)}, y_\mathrm{est}^{(B)})$ and the right figure is the trajectory under the worst wind condition $y^{(A)} \in \Y_{A}$ with $y^{(B)} = y_\mathrm{est}^{(B)}$. 
\Cref{fig:berth:B:cma} shows the trajectories observed for the best controller obtained by CMA-ES($y_\mathrm{est}^{(B)}$), which is the controller optimized under $y = (y_\mathrm{est}^{(A)}, y_\mathrm{est}^{(B)})$.
\Cref{fig:berth:B:adv} shows the trajectories observed for the best controller obtained by Adversarial-CMA-ES on $\Y_B$, which is the controller optimized under the worst model parameter $y^{(B)} \in \Y_{B}$ with wind condition $y^{(A)} = y_\mathrm{est}^{(A)}$. 
For \Cref{fig:berth:B:cma,fig:berth:B:adv}, the left figure is the trajectory under $y = (y_\mathrm{est}^{(A)}, y_\mathrm{est}^{(B)})$ and the right figure is the trajectory under the worst model parameter $y^{(B)} \in \Y_{B}$ with $y^{(A)} = y_\mathrm{est}^{(A)}$. 

\end{document}